\documentclass{article}
\usepackage{amsmath,amsxtra,amssymb,latexsym,amsfonts,amscd}
\usepackage{multicol,color}
\usepackage{float}
\usepackage{soul}
\usepackage{graphicx}
\definecolor{labelkey}{rgb}{0,0.08,0.45}
\definecolor{refkey}{rgb}{0,0.6,0.0}
\definecolor{Brown}{rgb}{0.45,0.0,0.05}
\definecolor{lime}{rgb}{0.00,0.8,0.0}
\definecolor{lblue}{rgb}{0.5,0.5,0.99}
\usepackage{pstricks-add}
\usepackage{amssymb}
\usepackage{theorem}
\usepackage{fancyhdr}
\usepackage{tikz}
\usepackage{enumerate}
\usepackage[margin=1.1in]{geometry}
\usetikzlibrary{arrows}
\usepackage{hyperref}
\definecolor{labelkey}{rgb}{0.6,0.6,0.6}
\definecolor{refkey}{rgb}{0,0.6,0.0}
\def\disp{\displaystyle}
\def\e{\epsilon}
\def\ve{\varepsilon}
\def\dd{\delta}

\def\lm{\lambda}
\def\O{\Omega}
\def\Tilde{\widetilde}

\def\({\left(}
\def\){\right)}
\def\[{\left[}
\def\]{\right]}
\def\n{\left \|}
\def\en{\right \|}
\def\nn {\left \{ }
\def\hnn {\right \}}
\def\l{\left}
\def\r{\right}

\def\oa{\bar a}
\def\ob{\bar b}

\def\ox{\bar{x}}
\def\oy{\bar{y}}
\def\oz{\bar{z}}

\def\ou{\bar{u}}

\def\oA{\bar A}
\def\oB{\bar B}
\def\oX{\bar X}

\def\gph{\hbox{}}

\def\gg{\gamma}
\def\dn{\downarrow}

\def\sr{\Longrightarrow }
\def\tto{\rightrightarrows}
\def\st{\stackrel}

\def\Limsup{\mathop{{\rm Lim}\,{\rm sup}}}

\def\hat{\widehat}
\def\Hat{\widehat}

\def\Tilde{\widetilde}
\def\Bar{\overline}
\def\la{\langle}
\def\ra{\rangle}
\def\ve{\varepsilon}
\def\B{I\!\!B}
\def\h{\hfill\Box}
\def\R{\mathbb{R}}
\def\N{\mathbb{N}}

\def\gra{\bigtriangledown }
\def\co{\mbox{\rm co}\,}
\def\gph{\mbox{\rm gph}\,}
\def\rep{\mbox{\rm rep}\,}
\def\epi{\mbox{\rm epi}\,}

\def\bd{\mbox{\rm bd}\,}

\def\dn{\downarrow}
\def\O{\Omega}

\def\vph{\varphi}
\def\emp{\emptyset}
\def\st{\stackrel}
\def\oR{\Bar{\R}}
\def\lm{\lambda}

\def\gg{\gamma}

\def\dd{\delta}

\def\al{\alpha}
\def\vth{\vartheta}

\def\be{\beta}
\def\Th{\Theta}

\def\N{I\!\!N}
\def\th{\theta}

\def\vth{\vartheta}

\newtheorem{theorem}{Theorem}[section]

\newtheorem{corollary}[theorem]{Corollary}
\newtheorem{proposition}[theorem]{Proposition}

\theoremstyle{plain}{\theorembodyfont{\rmfamily}
}
\theoremstyle{plain}{\theorembodyfont{\rmfamily}
}
\theoremstyle{plain}{\theorembodyfont{\rmfamily}
}
\theoremstyle{plain}{\theorembodyfont{\rmfamily}
\newtheorem{example}[theorem]{Example}}

\theoremstyle{plain}{\theorembodyfont{\rmfamily}
}

\def\eq{\begin{equation}}
\def\eeq{\end{equation}}

\begin{document}
\begin{center}
{\bf OPTIMIZATION OF FULLY CONTROLLED SWEEPING PROCESSES}\\[3ex]
TAN H. CAO\footnote{Department of Applied Mathematics and Statistics, State University of New York--Korea, Yeonsu-Gu, Incheon, Republic of Korea (tan.cao@stonybrook.edu). Research of this author is supported by the National Research Foundation of Korea grant funded by the Korea Government (MIST) NRF-2020R1F1A1A01071015.}
\quad GIOVANNI COLOMBO\footnote{Dipartimento di Matematica ``Tullio Levi-Civita'', Universit$\grave{\textrm{a}}$ di Padova, via Trieste 63, 35121 Padua, Italy (colombo@math.unipd.it). Research of this author is partially supported by the University of Padova grant SID 2018 ``Controllability, stabilizability and infimum gaps for control systems," BIRD 187147, and is affiliated to Istituto Nazionale di Alta Matematica (GNAMPA).}
\quad BORIS S. MORDUKHOVICH\footnote{Department of Mathematics, Wayne State University, Detroit, Michigan 48202, USA (boris@math.wayne.edu). Research of this author was partially supported by the USA National Science Foundation under grants DMS-1007132 and DMS-1512846, by the USA Air Force Office of Scientific Research grant \#15RT0462, and by the Australian Research Council under Discovery Project DP-190100555.} \quad DAO NGUYEN\footnote{Department of Mathematics, Wayne State University, Detroit, Michigan 48202, USA (dao.nguyen2@wayne.edu). Research of this author was partially supported by the USA National Science Foundation under grant DMS-1512846 and by the USA Air Force Office of Scientific Research grant \#15RT0462.}
\end{center}
\small{\bf Abstract.} The paper is devoted to deriving necessary optimality conditions in a general optimal control problem for dynamical systems governed by controlled sweeping processes with hard-constrained control actions entering both polyhedral moving sets and additive perturbations. By using the first-order and mainly second-order tools of variational analysis and generalized differentiation, we develop a well-posed method of discrete approximations, obtain optimality conditions for solutions to discrete-time control systems, and then establish by passing to the limit verifiable necessary optimality conditions for local minimizers of the original controlled sweeping process that are expressed entirely in terms of its given data. The efficiency of the obtained necessary optimality conditions for the sweeping dynamics is illustrated by solving three nontrivial examples of their own interest.\\[1ex]
\vspace{-0.2in}

\section{Introduction and Problem Formulation}\label{intro}
\setcounter{equation}{0}\vspace*{-0.1in}

The paper continues recent developments on necessary optimality conditions for controlled sweeping processes. The {\em sweeping dynamics} was originally described by Jean-Jacques Moreau \cite{moreau} in the form
\begin{equation}\label{SP}
\dot x(t)\in-N\big(x(t);C(t)\big)\;\mbox{ a.e. }\;t\in[0,T]\;\mbox{ with }\;x(0):=x_0\in C(0),
\end{equation}
where $N(x;C)$ stands for the normal cone of convex analysis defined by
\begin{equation}\label{NC}
N(x;C):=\big\{v\in\R^n\big|\;\la v,y-x\ra\le 0,\;y\in C\big\}\;\textrm{ if }\;x\in C\;\textrm{ and }\;N(x;C):=\emp\textrm{ if }\;x\notin C
\end{equation}
for the continuously moving convex set $C=C(t)$ at the point $x=x(t)$. Moreau's sweeping process \eqref{SP} and its modifications have been developed in dynamical system theory with many applications to various areas of mechanics, economics, traffic equilibria, robotics, etc.; see, e.g., \cite{bt,CT,HB,KMM,ve} with the references therein. However, optimization and control problems for sweeping processes were formulated much later. A primal reason for this situation is that the Cauchy problem in \eqref{SP} has a {\em unique} solution as shown in \cite{moreau} (see also the survey in \cite{CT,KMM} for subsequent developments), and so there is nothing to optimize. This is very different from the standard setting of control systems governed by ODE and Lipschitzian differential inclusions as in \cite{pont} and \cite{cl,m-book,v}. First optimal control problems for sweeping processes were formulated with controls functions acting in additive perturbations as in \cite{et}. Nevertheless, necessary optimality conditions for controlled sweeping processes were first established only in \cite{chhm1} (see also \cite{chhm3}) for a new class of problems with control functions acting in the moving set formalized as $C(t)=C(u(t))$ on $[0,T]$. Another type of controlled sweeping processes was introduced in \cite{bk}, with deriving necessary optimality conditions, where control functions entered a linear ODE system adjacent to the sweeping dynamics. The recent years have witnessed a rapidly growing interest to the derivation of necessary optimality conditions for various types of controlled sweeping processes with their broad applications to practical models; see, e.g., \cite{ao,ac,cm2,cm3,cg,chhm1,chhm3,cmn1,cmn,hm,pfs,mn,zeidan} and the references therein.

This paper addresses a general class of optimal control problems governed by a perturbed sweeping process over controlled polyhedral sets. Namely, we consider the optimal control problem $(P)$  described as follows: minimize
\begin{equation}\label{minimize}
J[x,a,b,u]:=\vph\big(x(T)\big)+\int_0^T\ell\(t,x(t),a(t),b(t),u(t),\dot{x}(t),\dot{a}(t),\dot{b}(t)\)dt
\end{equation}
subject to the perturbed sweeping dynamics
\begin{equation}\label{Problem}
\left\{\begin{array}{ll}
\dot{x}(t)\in-N\big(x(t);C(t)\big)+g\big(x(t),u(t)\big)\;\textrm{ a.e. }\;t\in[0,T],\\
x(0)=x_0\in C(0)\subset\R^n
\end{array}\right.
\end{equation}
with trajectories $x(\cdot)\in W^{1,2}([0,T];\R^n)$ generated by measurable controls $u(\cdot)\in L^2([0,T];\R^d)$ in the additive perturbations of \eqref{Problem} satisfying the constraint
\begin{equation}\label{u}
u(t)\in U\subset\R^d\;\textrm{ a.e. }\;t\in[0,T]
\end{equation}
as well as absolutely continuous controls $a(\cdot)=\(a_1(\cdot),\ldots,a_m(\cdot)\)\in W^{1,2}([0,T];\R^{mn})$ and $b(\cdot)=\(b_1(\cdot),\ldots, b_m(\cdot)\)\in W^{1,2}([0,T];{\R^m})$ acting in the moving set
\begin{equation}\label{SPC}
C(t):=\nn x\in\R^n\;\left|\;\la a_i(t),x\ra\le b_i(t),\;i=1,\ldots,m\right.\hnn\;\mbox{ for all }\;t\in[0,T]
\end{equation}
under the the pointwise constraints on the entire time interval given by
\begin{equation}\label{a}
\|a_i(t)\|=1\;\textrm{ for all }\;t\in[0,T],\;i=1,\ldots,m,
\end{equation}
\begin{equation}\label{ab}
\dot a_i(t)\in A_i\subset\R^n\;\mbox{ and }\;\dot b_i(t)\in B_i\subset\R\;\mbox{ a.e. }\;t\in[0,T],\;i=1,\ldots,m,
\end{equation}
where the points $x_0\in C(0)$, $a_{i0}:=a_i(0)\in\R^n$, $b_{i0}:=b_i(0)\in\R$ are fixed together with the final time $T>0$.

Problem $(P)$ was formulated in our previous paper \cite{ccmn19a} without imposing the pointwise constraints \eqref{ab} on the functions $a_i(t)$ and $b_i(t)$, which control the evolution of the moving set $C(t)$ in \eqref{SPC} and hence of the sweeping dynamics \eqref{Problem}. However, necessary optimality conditions for (either local or global) minimizers of $(P)$ were not obtained in \cite{ccmn19a} even for the mentioned particular case of the problem. In \cite{ccmn19a} we proved the existence of feasible and optimal solutions to the problem formulated therein under natural assumptions on the given data, constructed a well-posed sequence of discrete approximations of a designated local minimizer of the continuous-time problem with the $W^{1,2}\times L^2$-strong convergence of their optimal solutions, and derived necessary conditions for optimal solutions to discrete approximations. The first two issues (existence and convergence of discrete approximations) are valid, with minimal changes in the statements and the proofs, for the constrained control problem $(P)$ formulated above, while the derivation of necessary optimality conditions for the obtained versions of the discrete counterparts of the sweeping systems requires some work, which is done in what follows. However, the {\em main thrust} of this paper is on establishing {\em necessary optimality conditions} for local minimizers of the general problem $(P)$, which has never been accomplished earlier, except some particular cases when controls are acting separately in either perturbations, or moving sets; see \cite{ccmn19a} for more discussions on previously known results. Note that the enhanced technique developed in this paper allows us to improve and clarify some important optimality conditions even for the special cases of $(P)$ studied before. Observe also that the presence of the velocity constraints in \eqref{ab} significantly complicates the derivation of necessary optimality and leads us to new conditions of the maximum principle type, which are important in applications; see, e.g., illustration in Section~\ref{sec:exa}.

Observe further that, besides the hard/pointwise constraints on control functions in both the additive perturbations $g(x,u)$ and the moving set $C(t)$ imposed in \eqref{u}, \eqref{a}, and \eqref{ab}, we automatically have the mixed pointwise state-control constraints given by the bilinear inequalities
\begin{equation}\label{state}
\la a_i(t),x(t)\ra\le b_i(t)\;\mbox{ for all }\;t\in[0,T]\;\mbox{ and }\;i=1,\ldots,m,
\end{equation}
which follow from the sweeping inclusion \eqref{Problem} due to the normal cone definition \eqref{NC} ensuring that $x(t)\in C(t)$ on $[0,T]$. Together with the non-Lipschitzian/highly discontinuous sweeping dynamics in \eqref{Problem}, this makes deriving necessary optimality conditions for local minimizers of $(P)$ to be a very challenging task, which does not allow us to employ conventional techniques of the calculus of variations and optimal control theory.

We rely here on the {\em method of discrete approximations}, which was developed in \cite{m95} (see also \cite[Chapter~6]{m-book}) to derive necessary optimality conditions for control problems governed by Lipschitzian differential inclusions. The controlled sweeping dynamics and optimization problems for it are dramatically different---for any type of the control entering discussed above---from the Lipschitzian framework. Thus the method of discrete approximations requires significant modifications to be useful in the derivation of necessary optimality conditions. This has been done in \cite{cm2,cm3,cg,chhm1,chhm3,cmn1,cmn,hm,mn} for various types of sweeping control problems, while not for the problem $(P)$ formulated in \eqref{minimize}--\eqref{state}. Now we intend to establish necessary optimality conditions for local minimizers (in the sense precisely defined in Section~\ref{prel}) of the problem $(P)$ under consideration by employing the discrete approximation scheme and convergence results developed in \cite{ccmn19a}.

The rest of the paper is organized as follows. In Section~\ref{prel} we first formulate and discuss the standing assumptions and also the type of ``relaxed intermediate local minimizers" studied below for the original problem $(P)$. Then we construct a well-posed sequence of discrete approximation problems $(P^k)$ such that these problems admit optimal solutions and that each sequence of their optimal solutions strongly converges as $k\to\infty$ to the prescribed local minimizer of $(P)$. This preliminary section follows, without proofs, our previous paper \cite{ccmn19a} with rather small adaptations related to incorporating the new control constraints imposed in \eqref{ab}.

The next Section~\ref{sec:disc} is devoted to the derivation of necessary conditions for optimal solutions to each problem $(P^k)$ with any fixed $k\in\N$ by reducing it to a constrained problem of mathematical programming and using appropriate tools of variational analysis and generalized differentiation, mainly of the second-order. We implement here the same approach as in \cite{ccmn19a} while carefully incorporating the additional pointwise control constraints, which lead us to new optimality conditions in comparison with \cite{ccmn19a}.

Section~\ref{sec:sweep} is the culmination of this paper. It establishes new necessary optimality conditions for relaxed intermediate local minimizers of the original sweeping control problem $(P)$ by passing to the limit from the optimality conditions for the discrete approximation problems derived in Section~\ref{sec:disc}. The proof is quite involved being heavily based on advanced tools of variational analysis and second-order computations. The obtained results are expressed entirely in terms of the given data of $(P)$ and include, among various primal and dual relationships, novel local and global maximization conditions that can be treated as far-going sweeping counterparts of the Pontryagin Maximum Principle for the hard-constrained problems under consideration.

The concluding Section~\ref{sec:exa} presents three examples of sweeping optimal control problems of their own interest. The provided arguments and calculations illustrate how the main results of Section~\ref{sec:sweep} can be applied for all the three types of controls under consideration in $(P)$. Using the obtained necessary optimality conditions for the corresponding problems $(P)$ in these examples, together with imposing and justifying some additional assumptions on the class of controls considered therein, allows us to explicitly determine a unique solutions satisfying them.

Throughout the paper we use the standard notation of variational analysis; see, e.g., \cite{m18,rw}. Recall that, depending on the context, the sign $^*$ indicates the matrix transposition and the set duality/polarity defined by $\O^*:=\{v\in\R^n\;|\;\la v,x\ra\le 0\;\mbox{ for all }\;x\in\O\}$, where $\O\subset\R^n$. As usual, $\N:=\{1,2,\ldots\}$.\vspace*{-0.2in}

\section{Standing Assumptions and Preliminaries}\label{prel}
\setcounter{equation}{0}\vspace*{-0.1in}

According to the descriptions in Section~\ref{intro}, by {\em feasible solutions} to problem $(P)$ we understand the collections of controls $u(\cdot)\in L^2([0,T];\R^d)$, $a(\cdot)=(a_1(\cdot),\ldots,a_m(\cdot))\in W^{1,2}([0,T];\R^{mn})$, and $b(\cdot)=(b_1(\cdot),\ldots,b_m(\cdot))\in W^{1,2}([0,T];\R^m)$ together with the corresponding trajectories $x(\cdot)\in W^{1,2}([0,T];\R^n)$ of the sweeping dynamics, which satisfy all the relationships in \eqref{minimize}--\eqref{state}. Let us formulate the {\em standing assumptions} of the paper that ensure, in particular, the {\em existence} of feasible solutions:

{\bf(H1)} The control set $U$ in \eqref{u} is closed and bounded in $\R^d$.

{\bf(H2)} The perturbation mapping $g\colon\R^n\times\R^d\to\R^n$ in \eqref{Problem} is uniformly Lipschitz continuous with respect to both variables $x$ and $u\in U$, i.e., there exists a constant $L_g>0$ such that
\begin{equation}\label{e:gL}
\left\|g(x_1,u_1)-g(x_2,u_2)\right\|\le L_g\left(\left\|x_1-x_2\right\|+\left\|u_1-u_2\right\|\right)\;\;\mbox{for all}\;\;(x_1,u_1),(x_2,u_2)\in\R^n\times U.
\end{equation}
Moreover, the mapping $g$ satisfies the sublinear growth condition
\begin{equation*}
\|g(x,u)\|\le M\(1+\|x\|\)\mbox{ for all }\;u\in U\;\mbox{ with some }\;M>0.
\end{equation*}

{\bf(H3)} There exists a continuous function $\vartheta\colon[0,T]\to\R$ such that $\disp\sup_{t\in[0,T]}\vartheta(t)<0$ and
\begin{equation*}
C^0(t):=\big\{x\in\R^n\big|\;\la a_i(t),x\ra-b_i(t)<\vartheta(t),\;i=1,\ldots,m\big\}\ne\emp\;\mbox{ for all }\;t\in[0,T].
\end{equation*}

{\bf(H4)} The constraint sets $A_i$ and $B_i$ in \eqref{ab} are closed and bounded in $\R^n$ and $\R$, respectively.

{\bf (H5)} The terminal cost $\vph\colon\R^n\to\oR:=(-\infty,\infty]$ is lower semicontinuous (l.s.c.), while the running cost $\ell\colon\R^{2(n+nm+m)+d}\to\oR$ is bounded from below and l.s.c.\ around the reference feasible solution to $(P)$ for a.e.\ $t\in[0,T]$. Furthermore, $\ell$ is a.e.\ continuous in $t$ and is uniformly majorized by a summable function on $[0,T]$.\vspace*{0.03in}

The aforementioned existence of feasible solutions to $(P)$ under the assumptions in (H1)--(H3) follows from \cite[Theorem~2.1]{ccmn19a} the proof of which is based on the constructions of \cite{Tol}. More precisely, Theorem~2.1 from \cite{ccmn19a} states that any control triples $(u(\cdot),a(\cdot),b(\cdot))\in L^2([0,T];\R^d)\times W^{1,2}([0,T];\R^{mn})\times W^{1,2}([0,T];\R^m)$ generates a unique trajectory $x(\cdot)\in W^{1,2}([0,T];\R^n)$ of the sweeping inclusion \eqref{Problem}.

Since this paper deals with deriving necessary optimality conditions in $(P)$, it is natural to define an appropriate notion of {\em local} minimizers for our study. Adapting the concept of intermediate (between weak and strong) local minimizers introduced in \cite{m95} for Lipschitzian differential inclusions, recall that a feasible solution $(\ox(\cdot),\oa(\cdot),\ob(\cdot),\ou(\cdot))$ to $(P)$ is an {\em intermediate local minimizer} (i.l.m.) for $(P)$ if $(\ox(\cdot),\oa(\cdot),\ob(\cdot),\ou(\cdot))\in W^{1,2}([0,T];\R^n)\times W^{1,2}([0,T];\R^{mn})\times W^{1,2}([0,T];\R^m)\times L^2([0,T];\R^d)$ and there exists $\ve>0$ such that $J[\ox,\oa,\ob,\ou]\le J[x,a,b,u]$ for any feasible solutions $(x(\cdot),a(\cdot),b(\cdot),u(\cdot))$ to $(P)$ satisfying
\begin{equation}\label{ilm}
\big\|x(\cdot)-\ox(\cdot)\big\|_{W^{1,2}}+\big\|\big(a(\cdot),b(\cdot)\big)-\big(\oa(\cdot),\ob(\cdot)\big)\big\|_{W^{1,2}}+\|u(\cdot)-\ou(\cdot)\|_{L^2}\le\ve.
\end{equation}
If we replace $\|x(\cdot)-\ox(\cdot)\|_{W^{1,2}}$ in \eqref{ilm} by the ${\cal C}$-norm $\disp\max_{t\in[0,T]}\|x(t)-\ox(t)\|$, then we speak about a {\em strong local minimizer} of $(P)$, which is properly stronger than the notion of intermediate local minimizers for this problem.

However, in the general setting of $(P)$, without any a priori convexity assumptions, we need some stability of local minimizers with respect to performing limiting procedures. It can be achieved by the following relaxation of $(P)$ performed in the conventional line of the calculus of variations and optimal control; cf.\ \cite{m-book,v} for the case of ODE control systems and Lipschitzian differential inclusions and also \cite{dfm,et,Tol1} for non-Lipschitzian ones.

To proceed, define the set-valued mapping $F\colon\R^n\times\R^{mn}\times\R^m\times\R^d\tto\R^n$ by
\begin{equation}\label{F}
F(x,a,b,u):=N\big(x;C(a,b)\big)-g(x,u)\;\mbox{ with }\;C(a,b):=\big\{x\in\R^n\;\big|\;\la a_i,x\ra\le b_i,\;i=1,\ldots,m\big\}
\end{equation}
and deduce from the constructions in \eqref{NC} and \eqref{F} the explicit representation of $F$ in the form
\begin{equation}\label{F-rep}
F(x,a,b,u)=\bigg\{\sum_{i\in I(x,a,b)}\eta_ia_i\;\bigg|\;\eta_i\ge 0\bigg\}-g(x,u)\;\mbox{ for }\;x\in C(a,b),
\end{equation}
where the collection of active constraint indices of the polyhedron $C(a,b)$ is given by
\begin{equation}\label{e:AI}
I(x,a,b):=\big\{i\in\{1,\ldots,m\}\;\big|\;\la a_i,x\ra=b_i\big\}.
\end{equation}
Label by $\ell_F(t,a,b,u,\dot x,\dot a,\dot b)$ the restriction of the integrand $\ell$ on the set $F(x,a,b,u)$ with $\ell_F(t,a,b,u,\dot x,\dot a,\dot b):=\infty$ if $\dot x\notin
F(x,a,b,u)$, and then denote by $\Hat\ell_F$ the convexification of $\ell_F(t,a,b,u,\dot x,\dot a,\dot b)$ (i.e., the largest l.s.c.\ convex function majorized by $\ell(t,x,a,b,\cdot,\cdot,\cdot,\cdot)$) with respect to $(u,\dot x,\dot a,\dot b)$ when $u$ takes values from the the convex hull $\co U$ of the constraint set \eqref{u}. The {\em relaxed optimal control problem} $(R)$ associated with $(P)$ is defined as follows:
\begin{equation}\label{R}
\textrm{minimize}\hspace{0.1in}\hat
J[x,a,b,u]:=\vph\big(x(T)\big)+\int_0^T\hat\ell_F\big(t,x(t),a(t),b(t),u(t),\dot{x}(t),\dot{a}(t),\dot{b}(t)\big)dt
\end{equation}
over feasible solutions $(x(\cdot),a(\cdot),b(\cdot),u(\cdot))\in W^{1,2}([0,T];\R^n)\times W^{1,2}([0,T];\R^{mn})\times W^{1,2}([0,T];\R^m)\times L^2([0,T];\R^d)$ satisfying the constraints in \eqref{a} and \eqref{ab} and giving a finite value of the integrand $\hat\ell_F$ in \eqref{R}.

A quadruple $(\ox(\cdot),\oa(\cdot),\ob(\cdot),\ou(\cdot))$ is said to be a {\em relaxed intermediate local minimizer} (r.i.l.m.) for problem $(P)$ if it is feasible to this problem with $\hat J[\ox,\oa,\ob,\ou]=J[\ox,\oa,\ob,\ou]$ and if there exists a number $\ve>0$ such that $J[\ox,\oa,\ob,\ou]\le\hat J[x,a,b,u]$ for any feasible solution $(x(\cdot),a(\cdot),b(\cdot),u(\cdot))$ to $(R)$. satisfying \eqref{ilm}. Observing that there is no difference between intermediate local minimizers and their relaxed counterparts in the case where the sets $U$ and $g(x;U)$ are convex and the integrand $\ell(t,x,a,b,\cdot,\cdot,\cdot,\cdot)$ is convex with respect to the control and velocity variables, we refer the reader to \cite{ccmn19a} for a detailed discussion on the fulfillment of this phenomenon without imposing any convexity assumptions due the ``hidden convexity" inherent in such systems.

The main goal of this paper is to establish necessary optimality conditions for the given relaxed intermediate local minimizer $(\ox(\cdot),\oa(\cdot),\ob(\cdot),\ou(\cdot))$ of problem $(P)$ using the method of discrete approximations. Proceeding similarly to \cite{ccmn19a} with additional approximations of the new constraints in \eqref{ab}, for each $k\in\N$ consider the discrete partition/mesh of $[0,T]$ defined by
\begin{equation*}
\Delta_k:=\left\{0=t^k_0<t^k_1<\ldots<t^k_{\nu(k)-1}<t^k_{\nu(k)}=T\right\}\;\;\mbox{with}\;\;h^k_j:=t^k_{j+1}-t^k_j\le\dfrac{\Tilde\nu}{\nu(k)}\;\;\mbox{for}\;\;j=0,\ldots,\nu(k)-1,
\end{equation*}
where $\nu(\cdot)\in\N$ depends on $k$ such that $\nu(k)\ge k$, and where $\Tilde\nu>0$ is some constant. Then we approximate the original problem $(P)$ by the discrete-time optimization problems $(P^k)$ as follows:
\begin{equation*}
\begin{aligned}
&\textrm{minimize }J_k[x^k,a^k,b^k,u^k]:=\vph(x^k_{\nu(k)})+\sum^{\nu(k)-1}_{j=0}h^k_j\ell\bigg(t^k_j,x^k_j,a^k_j,b^k_j,u^k_j,\dfrac{x^k_{j+1}-x^k_i}{h^k_j},\dfrac{a^k_{j+1}-a^k_j}{h^k_j},\dfrac{b^k_{j+1}-b^k_j}{h^k_j}
\bigg)\\
&+\frac{1}{2}\sum_{j=0}^{\nu(k)-1}\int_{t^k_j}^{t^k_{j+1}}\n\(\dfrac{x^k_{j+1}-x^k_j}{h^k_j},\dfrac{a^k_{j+1}-a^k_j}{h^k_j},\dfrac{b^k_{j+1}-b^k_j}{h^k_j},u^k_j\)-\(\dot{\ox}(t),\dot{\oa}(t),
\dot{\ob}(t),\ou(t)\)\en^2dt
\end{aligned}
\end{equation*}
over the discrete quadruples represented by
\begin{eqnarray*}
(x^k,a^k,b^k,u^k)=(x^k_0,x^k_1,\ldots,x^k_{\nu(k)},a^k_{i0},a^k_{i1},\ldots,a^k_{i{\nu(k)}},b^k_{i0},b^k_{i1},\ldots,\;b^k_{i{\nu(k)}},u^k_0,u^k_1,\ldots,u^k_{\nu(k)-1})
\end{eqnarray*}
with $i=1,\ldots,m$ subject to the constraints
\begin{equation*}
x^k_{j+1}\in x^k_j-h^k_jF(x^k_j,a^k_j,b^k_j,u^k_j),\;j=0,\ldots,\nu(k)-1,
\end{equation*}
\begin{equation*}
\la a^{k}_{i\nu(k)},x^k_{\nu(k)}\ra\le b^{k}_{i\nu(k)},\;i=1,\ldots,m,
\end{equation*}
\begin{equation*}
x^k_0=x_0\in C(0),\;a^k_0=a_0,\;b^k_0=b_0,\;u^k_0=\ou(0),
\end{equation*}
\begin{equation*}
\sum_{j=0}^{\nu(k)-1}\int_{t^k_j}^{t^k_{j+1}}\n\(x^k_j,a^k_j,b^k_j,u^k_j\)-\(\ox(t),\oa(t),\ob(t),\ou(t)\)\en^2\;dt\le \dfrac{\ve}{2},
\end{equation*}
\begin{equation*}
\sum_{j=0}^{\nu(k)-1}\int_{t^k_j}^{t^k_{j+1}}\n\(\dfrac{x^k_{j+1}-x^k_j}{h^k_j},\dfrac{a^k_{j+1}-a^k_j}{h^k_j},\dfrac{b^k_{j+1}-b^k_j}{h^k_j}\) - \(\dot{\ox}(t),\dot{\oa}(t),\dot{\ob}(t)\)\en^2\;dt \le \dfrac{\ve}{2},
\end{equation*}
\begin{equation}\label{disc-ab}
u^k_j\in U,\;\alpha_{ij}^k:=\dfrac{a^k_{i\,(j+1)}-a^k_{ij}}{h^k_j}\in A_i,\;\beta_{ij}^k:=\dfrac{b^k_{i\,(j+1)}-b^k_{ij}}{h^k_j}\in B_i,\;j=0,\ldots,\nu(k)-1,\;i=1,\ldots,m,
\end{equation}
\begin{equation}\label{const_a}
1-\ve_k\le\|a^k_{ij}\|\le 1+\ve_k,\;i=1,\ldots,m,\;j=0,\ldots,\nu(k),
\end{equation}
where $\ve_k>0$ are fixed numbers that are sufficiently small.

Further, we say that the {\em positive linear independence constraint qualification} (PLICQ) holds at $x\in C(t)$ if
\begin{equation}\label{PLICQ}
\bigg[\sum_{i\in I(x,a(t),b(t))}\al_ia_i=0,\,\al_i\ge 0\bigg]\Longrightarrow\big[\al_i=0\;\;\mbox{for all}\;\;i\in I(x,a(t),b(t))\big],\quad t\in[0,T].
\end{equation}
This constraint qualification ensures the fulfillment of the local version of the estimate
\begin{equation}\label{A5'}
\sum_{i\in I(x,a(t),b(t))}\lm_i\n a_i(t)\en\le\sigma\bigg\|\sum _{i\in I(x,a(t),b(t))}\lm_i a_i(t)\bigg\|\;\mbox{ for all }\;x\in C(t)\;\mbox{ and }\;\lm_i\ge 0
\end{equation}
with some constant $\sigma>0$ that is known as the {\em inverse triangle inequality}. Conversely, \eqref{A5'} yields PLICQ at all $x\in C(t)$; see \cite{ccmn19a,ve} for more discussions. It is obvious that PLICQ is essentially weaker that the classical {\em linear independence constraint qualification} (LICQ) at $x\in C(t)$, which corresponds to \eqref{PLICQ} with arbitrary $\al_i\in\R$ instead of the nonnegative ones. While PLICQ plays a crucial role in deriving the main necessary optimality conditions of this paper, LICQ is used only to get the dynamic complementary slackness conditions in Theorem~\ref{Th2}.

First we present the following important result of the method of discrete approximations ensuring the strong convergence of optimal solutions of discrete problems $(P^k)$ to the prescribed r.i.l.m.\ of problem $(P)$. It is established in \cite[Theorem~5.2]{ccmn19a} for problem $(P)$ without the constraints in \eqref{ab} and their discrete approximations in \eqref{disc-ab}. However, the presence of the new constraints can be easily incorporated into the proof, and so we omit it here while referring the reader to \cite{ccmn19a}.\vspace*{-0.1in}

\begin{theorem}[\bf strong convergence of discrete optimal solutions]\label{str-conver} Let $(\ox(\cdot),\oa(\cdot),\ob(\cdot),\ou(\cdot))$ be an r.i.l.m.\ for problem
$(P)$  such that the functions $\dot\ox(\cdot),\dot\oa(\cdot),\dot\ob(\cdot)$ and $\ou(\cdot)$ are of bounded variation ${\rm(BV)}$ on $[0,T]$. Suppose in addition to the standing assumptions {\rm(H1)--(H5)} that the terminal cost $\vph$ is continuous around $\ox(T)$, that the running cost $\ell$ is continuous at $\big(t,\ox(t),\oa(t),\ob(t),\ou(t),$ $\dot\ox(t),\dot\oa(t),\dot\ob(t)\big)$
for a.e.\ $t\in[0,T]$, and that $\ell$ is uniformly majorized around the given r.i.l.m.\ by a summable function on $[0,T]$, i.e., there exists $\mu(\cdot)\in L^1([0,T];\R_+)$ with
\begin{equation*}
\big|\ell\big(t,x,a,b,u,\dot x,\dot a,\dot b\big)\big|\le\mu(t)\;\mbox{ for all }\;(x,a,b,u,\dot x,\dot u,\dot a,\dot b)\;\mbox{ near }\;\big(\ox(t),\oa(t),\ob(t),\ou(t),\dot\ox(t),\dot\oa(t),\dot\ob(t)\big)\;\mbox{ a.e. }\;t\in[0,T].
\end{equation*}
Take any sequence of optimal solutions $\(\ox^k(\cdot),\oa^k(\cdot),\ob^k(\cdot),\ou^k(\cdot)\)$ to the discrete problems $(P_k)$ and extend them to the entire interval $[0,T]$ piecewise
linearly for $\(\ox^k(\cdot),\oa^k(\cdot),\ob^k(\cdot)\)$ and piecewise constantly for $\ou^k(\cdot)$. Then the extended sequence of $\(\ox^k(\cdot),\oa^k(\cdot),\ob^k(\cdot),\ou^k(\cdot)\)$ converges to $\(\ox(\cdot),\oa(\cdot),\ob(\cdot),\ou(\cdot)\)$ as $k\to\infty$ in the norm topology of $W^{1,2}([0,T];\R^n)\times W^{1,2}([0,T];\R^{mn})\times W^{1,2}([0,T];\R^m )\times L^2([0,T];\R^d)$.
\end{theorem}
\vspace*{-0.25in}

\section{Necessary Optimality Conditions for Discrete Approximations}\label{sec:disc}
\setcounter{equation}{0}\vspace*{-0.1in}

In this section we derive necessary optimality conditions for discrete approximation problems $(P^k)$ constructed in Section~\ref{prel} for each $k\in\N$. Since the presence of the pointwise constraints on $\alpha_{ij}^k$ and $\beta_{ij}^k$ in \eqref{disc-ab} is essential for the optimality conditions in $(P^k)$ and for the subsequent limiting procedure from $(P^k)$ to $(P)$ as $k\to\infty$, we give a detailed proof of the next theorem in the line of \cite[Theorem~7.2]{ccmn19a} for discrete-time problems without those constraints while mainly concentrating on the handling of the new constraints. Recall the notation
\begin{equation*}
\[\alpha,y\]:=(\alpha_1y_1,\ldots,\alpha_my_m)\in\R^{nm}\;\mbox{ and }\;{\rm rep}_m(x):=(x,\ldots,x)\in\R^{nm}
\end{equation*}
for any vectors $\alpha=(\alpha_1,\ldots,\alpha_m)\in\R^m$ and $y=(y_1,\ldots,y_m)\in\R^{nm}$ with $y_i\in\R^n$ as $i=1,\ldots,m$.

To formulate and prove the necessary optimality conditions for problems $(P^k)$ and $(P)$ that are given in this and next sections, respectively, we first define and briefly discuss the constructions of generalized differentiation used in what follows; see the books \cite{m-book,m18,rw} for more details.

Let $\O\subset\R^n$ be locally closed around $\ox\in\O$. The (Mordukhovich limiting) {\em normal cone} to $\O$ at $\ox$ is
\begin{equation}\label{nor_con}
N(\ox;\O)=N_\O(\ox):=\big\{v\in\R^n\;\big|\;\exists\,x_k\to\ox,\;w_k\in\Pi(x_k;\O),\;\al_k\ge 0\;\mbox{ with }\;\al_k(x_k-w_k)\to v\big\},
\end{equation}
where $\Pi(x;\O):=\{w\in\O\;|\;\|x-w\|={\rm dist}(x;\O)\}$ stands for the Euclidean projection of $x\in\R^n$ to $\O$. If $\O$ is convex, then \eqref{nor_con} reduces to \eqref{NC}, but in general this normal cone is nonconvex while satisfying comprehensive calculus rules together with the associated subdifferential and coderivative constructions for functions and multifunctions.

Given a set-valued mapping (multifunction) $S\colon\R^n\tto\R^m$, its {\em coderivative} at $(\ox,\oy)\in\gph S$ is given by
\begin{equation}\label{coderivative}
D^*S(\ox,\oy)(u):=\big\{v\in\R^n\;\big|\;(v,-u)\in N\big((\ox,\oy);\gph S\big)\big\},\quad u\in\R^m,
\end{equation}
via the normal cone \eqref{nor_con} to the graph $\gph S:=\{(x,y)\in\R^n\times\R^m\;|\;y\in S(x)\}$. For single-valued smooth mappings $S\colon\R^n\to\R^m$ we have $D^*S(\ox)(u)=\{\nabla S(\ox)^*u\}$, when $Q^*$ signifies the adjoint/transposed matrix of $Q$.

The {\em subdifferential} of an l.s.c.\ function $\phi\colon\R^n\to\oR$ at $\ox$ with $\phi(\ox)<\infty$ is defined by
\begin{equation}\label{sub}
\partial\phi(\ox):=\big\{v\in\R^n\big|\;(v,-1)\in
N\big((\ox,\phi(\ox));\epi\phi\big)\big\}
\end{equation}
via the normal cone \eqref{nor_con} of the epigraph $\epi\phi:=\{(x,\al)\in\R^{n+1}\;|\;\al\ge\phi(x)\}$.

We have the following upper estimate of the coderivative \eqref{coderivative} of the velocity mapping \eqref{F} via the problem data taken from \cite[Theorem~6.2]{ccmn19a}, where $A:=[a_{ij}]$ as $i=1,\ldots,m$ and $j=1,\ldots,n$ with the vector columns $a_i$.\vspace*{-0.1in}

\begin{proposition}[\bf coderivative evaluations]\label{Th7} Consider the multifunction $F$ from \eqref{F}, where $g(x,u)$ is ${\cal C}^1$ around the reference points. Suppose that the vectors $\nn a_i\;|\;i\in I(x,a,b)\hnn$ are positively linearly independent at any triple $(x,a,b)\in\R^n\times\R^{mn}\times\R^m$. Then for all such triples and all $(w,u)\in\R^m\times U$ with $w+g(x,u)\in N(x;C(a,b))$ we have the coderivative upper estimate
\begin{equation*}
\begin{aligned}
D^*F(x,a,b,u,w)(y)\subset\bigcup_{\substack{p\in N_{\R^m_-}(Ax-b),\;A^*p=w+g(x,u)\\q\in D^*N_{\R^m_-}(Ax-b,p)(Ay)}}
\(\begin{array}{c}
A^*q-\nabla_x g(x,u)^*y\\
\hline
p_1y+q_1x\\
\vdots\\
p_my+q_mx\\
\hline
-q\\
-\nabla_u g(x,u)^*y
\end{array}\)
\end{aligned}
\end{equation*}
for any $y\in\disp\bigcap_{\nn i\;|\;p_i>0\hnn}a^\perp_i$, where the vector
$q\in\R^m$ satisfies the conditions
\begin{equation*}
\left\{\begin{array}{ll}
q_i=0\;\mbox{ for all }\;i\;\mbox{ such that either }\;\la a_i,x\ra<b_i\;\mbox{ or }\;p_i=0,\;\mbox{ or}\;\la a_i,y\ra<0,\\
q_i\ge 0\;\mbox{ for all }\;i\;\mbox{ such that }\;\la a_i,x\ra=b_i,\;p_i=0,\;\mbox{ and }\;\la a_i,y\ra>0.
\end{array}
\right.
\end{equation*}
\end{proposition}

Now we are ready to derive necessary conditions for optimal solutions to problem $(P_k)$ for each fixed $k\in\N$. The obtained optimality conditions extend those in \cite[Theorem~7.2]{ccmn19a} and are expressed entirely via the given data with the usage of the generalized differential constructions from \eqref{nor_con}--\eqref{sub}, where the subdifferential of $\ell$ is taken with respect to all but time variables, and where we skip indicating the dependence of $\ell$ on $t$.\vspace*{-0.1in}

\begin{theorem}[\bf necessary conditions for discrete optimal solutions]\label{Th1} Let
\begin{equation}\label{disc-opt}
(\ox^k,\oa^k,\ob^k,\ou^k)=(\ox_0^k,\ldots,\ox_{\nu(k)}^k,\oa_0^k,\ldots,\oa_{\nu(k)}^k,\ob_0^k,\ldots,\ob_{\nu(k)}^k,\ou_0^k,\ldots,\ou_{\nu(k)-1}^k)
\end{equation}
be an optimal solution to $(P^k)$ for a fixed index $k\in\N$, where $g$ is ${\cal C}^1$-smooth while $\vph$ and $\ell$ are locally Lipschitzian around the corresponding components of \eqref{disc-opt}. Assume also that the vectors $\nn\oa^k_i\;|\;i\in I(\ox^k,\oa^k,\ob^k)\hnn$ are positively linearly independent and define the auxiliary quadruples by
\begin{equation*}
\begin{aligned}
&\left(\th_{j}^{uk},\th_{j}^{xk},\th_{j}^{ak},\th_{j}^{bk}\right):=\\
&\left(\int_{t_j^k}^{t_{j+1}^k}\(\ou_j^k-\ou(t)\)dt,\int_{t_j^k}^{t_{j+1}^k}\(\frac{\ox_{j+1}^k-\ox_j^k}{h^k_j}-\dot{\ox}(t)\)dt,\int_{t_j^k}^{t_{j+1}^k}
\(\frac{\oa_{j+1}^k-\oa_j^k}{h^k_j}-\dot{\oa}(t)\)dt, \int_{t_j^k}^{t_{j+1}^k}\(\frac{\ob_{j+1}^k-\ob_j^k}{h^k_j}-\dot{\ob}(t)\)dt\right).
\end{aligned}
\end{equation*}
Then there exist dual elements
$\lm^k\ge 0$, $\psi^{uk}=(\psi^{uk}_0,\ldots,\psi^{uk}_{\nu(k)-1})\in\R^{d\nu(k)}$,
$\psi^{ak}=(\psi^{ak}_0,\ldots,\psi^{ak}_{\nu(k)-1})\in\R^{mn\nu(k)}$, $\psi^{bk}=(\psi^{bk}_0,\ldots,\psi^{bk}_{\nu(k)-1})\in\R^{m\nu(k)}$, $\xi^k=(\xi_1^k,\ldots,\xi_{m}^k)\in\R^{m}_+$, $p_j^k=\(p_{j}^{xk},p_{j}^{ak},p_{j}^{bk}\)\in\R^{n+mn+m}$ as
$j=0,\ldots,\nu(k)$ together with
$\gg^k=(\gg_0^k,\ldots\gg^k_{\nu(k)-1})\in\R^{m\nu(k)}$, $\al^{1k}=\(\al_{0}^{1k},\ldots,\al_{\nu(k)}^{1k}\)\in\R^{m(\nu(k)+1)}_+$,
$\al^{2k}=(\al_{0}^{2k},\ldots,\al_{\nu(k)}^{2k})\in\R^{m(\nu(k)+1)}_{-}$, and $\eta^k=(\eta_0^k,\ldots,\eta^k_{\nu(k})\in\R^{m(\nu(k)+1}_+$ as well as subgradients
\begin{equation}\label{subcol}
\(w_{j}^{xk},w_{j}^{ak},w_{j}^{bk},w_{j}^{uk},v_{j}^{xk},v_{j}^{ak},v_{j}^{bk}\)\in\partial\ell\Big(\ox_j^k,\oa_j^k,\ob_j^k,\ou_j^k,\frac{\ox_{j+1}^k-\ox_j^k}{h^k_j},
\frac{\oa_{j+1}^k-\oa_j^k}{h^k_j},\frac{\ob_{j+1}^k-\ob_j^k}{h^k_j}\Big)
\end{equation}
as $j=0,\ldots,\nu(k)-1$ such that the following necessary optimality conditions are satisfied:\\[1ex]
{\sc $\bullet$ primal-dual dynamic equations} for all $j=0,\ldots,\nu(k)-1$ and $i=1,\ldots,m:$
\begin{equation}\label{87}
-\frac{\ox_{j+1}^k-\ox_j^k}{h^k_j}+g(\ox_j^k,\ou_j^k)=\sum_{i=1}^{m}\eta_{ij}^k\oa_{ij}^k,
\end{equation}
\begin{equation}\label{conx}
\frac{p_{j+1}^{xk}-p_{j}^{xk}}{h^k_j}-\lm^kw_{j}^{xk}=-\nabla_xg(\ox_j^k,\ou_j^k)^*\(-\frac{1}{h^k_j}\lm^k\th_{j}^{xk}-\lm^k v_{j}^{xk}+p_{j+1}^{xk}\)+\sum_{i=1}^{m}\gg_{ij}^k\oa_{ij}^k,
\end{equation}
\begin{equation}\label{cona}
\begin{array}{ll}
&\dfrac{p_{j+1}^{ak}-p_{j}^{ak}}{h^k_j}-\lm^kw_{j}^{ak}-\dfrac{2}{h^k_j}\[\al_{j}^{1k}+\al_{j}^{2k},\oa_{j}^k\]\\
&\quad\qquad =\[\gg_{j}^k,\rep_m(\ox_j^k)\]+\[\eta_j^k,\rep_m\(-\dfrac{1}{h^k_j}\lm^k\th_{j}^{xk}-\lm^k v_{j}^{xk}+p_{j+1}^{xk}\)\],
\end{array}
\end{equation}
\begin{equation}\label{conb}
\frac{p_{j+1}^{bk}-p_{j}^{bk}}{h^k_j}-\lm^kw_{j}^{bk}=-\gg_{j}^k,
\end{equation}
\begin{equation}\label{congg}
\gg^k_{ij}\in D^*N_{\R_-}\(\la\oa_{ij}^k,\ox_j^k\ra -\ob_{ij}^k,\eta^k_{ij}\)\Big(\la\oa_{ij}^k,-\dfrac{1}{h^k_j}\lm^k\th_{j}^{xk}-\lm^k v_{j}^{xk}+p_{j+1}^{xk}\ra\Big),
\end{equation}
which implies that $\gg^k_{ij}$ satisfies the following relations:
\eq\label{congg1}
\begin{cases}
\gg^k_{ij}=0\;\;\mbox{if}\;\;\la\oa^k_{ij},\ox^k_j\ra<\ob^k_{ij}\;\;\mbox{or}\;\;\eta^k_{ij}=0,\;\la\oa_{ij}^k,-\dfrac{1}{h^k_j}\lm^k\th_{j}^{xk}-\lm^k v_{j}^{xk}+p_{j+1}^{xk}\ra<0,\\
\gg^k_{ij}\ge 0\;\;\mbox{if}\;\;\la\oa^k_{ij},\ox^k_j\ra=\ob^k_{ij},\;\eta^k_{ij}=0,\;\la\oa_{ij}^k,-\dfrac{1}{h^k_j}\lm^k\th_{j}^{xk}-\lm^k v_{j}^{xk}+p_{j+1}^{xk}\ra>0,\\
\gg^k_{ij}\in\R \;\;\mbox{if}\;\;\eta^k_{ij}>0,\;\la\oa_{ij}^k,-\dfrac{1}{h^k_j}\lm^k\th_{j}^{xk}-\lm^k v_{j}^{xk}+p_{j+1}^{xk}\ra=0,
\end{cases}
\eeq
\begin{equation}\label{cony}
-\dfrac{1}{h^k_j}\lm^k\th_{j}^{uk}-\lm^kw_{j}^{uk}-\dfrac{1}{h^k_j}\psi_j^{uk}=-\nabla_ug(\ox_j^k,\ou_j^k)^*\Big(-\frac{1}{h^k_j}\lm^k\th_{j}^{xk}-\lm^k v_{j}^{xk}+p_{j+1}^{xk}\Big),
\end{equation}
\eq\label{dloc-max}
\begin{cases}
\psi_j^{uk}\in N(\ou_j^k;U),\\[1ex]
\disp\psi_{j}^{Ak}:=-\frac{1}{h^k_j}\lm^k\th_{j}^{ak}+p_{j+1}^{ak}-\lm^kv_{j}^{ak}\in N(\dot{\oa}^k_{j};A),\\[1ex]
\disp\psi_{j}^{Bk}:=-\frac{1}{h^k_j}\lm^k\th_{j}^{bk}+p_{j+1}^{bk}-\lm^kv_{j}^{bk}\in N(\dot{\ob}^k_{j};B)
\end{cases}
\eeq
with the notation $A:=\prod_{i=1}^{m}A_i$ and $B:=\prod_{i=1}^{m}B_i$.\\[1ex]
{$\bullet$ \sc transversality conditions} for all $i=1,\ldots,m$:
\begin{equation}\label{nmutx}
-p_{\nu(k)}^{xk}\in\lm^k\partial\vph(\ox_{\nu(k)}^k)+\sum_{i=1}^{m}\eta_{i\nu(k)}^k\oa_{i\nu(k)}^k,
\end{equation}
\begin{equation}\label{nmuta}
 p_{\nu(k)}^{ak}=-2\[\al_{\nu(k)}^{1k}+\al_{\nu(k)}^{2k},\oa_{\nu(k)}^k\]-\[\eta_{\nu(k)}^k,\rep_m(\ox_{\nu(k)}^k)\],
\end{equation}
\begin{equation}\label{nmutb}
p_{i\nu(k)}^{bk}=\eta^k_{i\nu(k)}\ge 0,\;\la \oa^{k}_{i\nu(k)},\ox^k_{\nu(k)}\ra<\ob^k_{i\nu(k)}\sr p_{i\nu(k)}^{bk}=0,
\end{equation}
with dual vectors $\al_{\nu(k)}^{1k}$ and $\al_{\nu(k)}^{2k}$ taken from
\begin{equation}\label{t:7.20}
\begin{cases}
\al_{\nu(k)}^{1k}\in N_{[0,1+\ve_k]}\big(\|\oa_{\nu(k)}^k\|\big),\\[1ex]
\al_{\nu(k)}^{2k}\in N_{[1-\ve_k,\infty]}\big(\|\oa_{\nu(k)}^k\|\big).
\end{cases}
\end{equation}
{$\bullet$ \sc complementarity slackness conditions}:
\begin{equation}\label{eta}
\[\la a_{ij}^k,\ox_j^k\ra<\ob_{ij}^k\]\sr\eta_{ij}^k=0\;\mbox{ for all }j=0,\ldots,\nu(k)-1\;\mbox{ and }\;i=1,\ldots,m,
\end{equation}
where in addition we have the implications
\begin{equation}\label{94}
\begin{cases}
\[\la\oa_{ij}^k,\ox_j^k\ra<\ob_{ij}^k\]\sr\gg_{ij}^k=0,\\[1ex]
\[\la\oa_{i\nu(k)}^k,\ox_{\nu(k)}^k\ra<\ob_{i\nu(k)}^k\]\sr\eta_{i\nu(k)}^k=0,
\end{cases}
\end{equation}
\begin{equation}\label{96}
\eta_{ij}^k>0\sr\[\Big\la\oa_{ij}^k,-\frac{1}{h^k_j}\lm^k\th_{j}^{xk}-\lm^kv_{j}^{xk}+p_{j+1}^{xk}\Big\ra=0\].
\end{equation}
{$\bullet$ \sc nontriviality conditions} with $\psi^k=\(\psi^{Ak},\psi^{Bk},\psi^{uk}\)$:
\begin{equation}\label{ntc0}
\lm^k+\n\al^{1k}+\al^{2k}\en+\disp\n\eta_{\nu(k)}^k\en+\sum_{j=0}^{\nu(k)-1}\n p_{j}^{xk}\en+\n p_{0}^{ak}\en+\n p_{0}^{bk}\en+\sum_{j=0}^{\nu(k)-1}\|\psi^k_j\|\ne 0,
\end{equation}
\begin{equation}\label{ntc1}
\lm^k+\n\al^{1k}+\al^{2k}\en+\n\gg^k\en+\n p^{ak}_{\nu(k)}\en+\n p^{bk}_{\nu(k)}\en\ne 0.
\end{equation}
\end{theorem}\vspace*{-0.1in}
{\bf Proof.} Fix $\ve>0$ from the definition of the given relaxed intermediate local minimizer $(\ox(\cdot),\oa(\cdot),\ob(\cdot),\ou(\cdot))$ for the original problem $(P)$ and form the cumulating vector
\begin{equation*}
z:=\(x_0^k,\ldots,x_{\nu(k)}^k,a_0^k,\ldots,a_{\nu(k)}^k,b_0^k,\ldots,b_{\nu(k)}^k,u_0^k,\ldots,u_{\nu(k)-1}^k,X_0^k,\ldots,X_{\nu(k)-1}^k,A_0^k,\ldots,A_{\nu(k)-1}^k,B_0^k,\ldots,B_{\nu(k)-1}^k\)
\end{equation*}
of feasible solutions to $(P^k)$. We clearly rewrite $(P^k)$ in the form of mathematical programming $(MP)$:
\begin{equation*}
\begin{aligned}
&\textrm{minimize }\phi_0(z):=\vph(x^k_{\nu(k)})+\sum_{j=0}^{\nu(k)-1}h^k_j\ell(x_j^k,a_j^k,b_j^k,u_j^k,X_j^k,A_j^k,B_j^k)\\
&+\frac{1}{2}\sum_{j=0}^{\nu(k)-1}\int_{t_j^k}^{t_{j+1}^k}\n\(X_j^k-\dot{\ox}(t),A_j^k-\dot{\oa}(t),B_j^k-\dot{\ob}(t),u_j^k-\ou(t)\)\en^2dt
\end{aligned}
\end{equation*}
subject to finitely many equality, inequality, and geometric constraints
\begin{equation}\label{e:5.13*}
\kappa(z):=\sum_{j=0}^{\nu(k)-1}\int_{t^k_j}^{t^k_{j+1}}\n\(x^k_j,a^k_j,b^k_j,
u^k_j\)-\(\ox(t),\oa(t),\ob(t),\ou(t)\)\en^2dt\le\dfrac{\ve}{2},
\end{equation}
\begin{equation}\label{e:5.14*}
\phi(z):=\sum_{j=0}^{\nu(k)-1}\int_{t_j^k}^{t_{j+1}^k}\n\(X_j^k,A_j^k,B_j^k,u_j^k\)-\big(\dot{\ox}(t),\dot{\oa}(t),\dot{\ob}(t),\ou(t)\big)\en^2dt-\frac{\e}{2}\le 0,
\end{equation}
\begin{equation*}
g^x_j(z):=x_{j+1}^k-x_j^k-h^k_jX_j^k=0,\;j=0,\ldots,\nu(k)-1,
\end{equation*}
\begin{equation*}
g^a_j(z):=a_{j+1}^k-a_j^k-h^k_jA_j^k=0,\;j=0,\ldots,\nu(k)-1,
\end{equation*}
\begin{equation*}
g^b_j(z):=b_{j+1}^k-b_j^k-h^k_jB_j^k=0,\;j=0,\ldots,\nu(k)-1,
\end{equation*}
\begin{equation*}
q_i(z):=\la a_{i\nu(k)}^k,x_{\nu(k)}^k\ra-b_{i\nu(k)}^k\leq 0,\;i=1,\ldots,m,
\end{equation*}
\begin{equation*}
l^1_{ij}(z):=\|a_{ij}^k\|^2-(1+\dd_k)^2\le 0,\;i=1,\ldots,m,\;j=0,\ldots,\nu(k),
\end{equation*}
\begin{equation*}
l^2_{ij}(z):=\|a_{ij}^k\|^2-(1-\dd_k)^2\ge 0,\;i=1,\ldots,m,\;j=0,\ldots,\nu(k),
\end{equation*}
\begin{equation}\label{e:5.16*}
z\in\Xi_j:=\big\{z\big|\;-X_j^k\in F(x_j^k,a_j^k,b_j^k,u_j^k)\big\},\;j=0,\ldots,\nu(k)-1,
\end{equation}
\begin{equation}\label{uinU}
z\in\O_j^u=\big\{z\big|\;u^k_j\in U\big\},\;\;j=0,\ldots,\nu(k)-1,
\end{equation}
\begin{equation}\label{uinU1}
z\in\O_{j}^a=\big\{z\big|\; \dot a^k_{j}\in A\big\},\;i=1,\ldots,m,\;j=0,\ldots,\nu(k)-1,
\end{equation}
\begin{equation}\label{uinU2}
z\in\O_{j}^b=\big\{z\big|\; \dot b^k_{j}\in B\big\},\;i=1,\ldots,m,\;j=0,\ldots,\nu(k)-1.
\end{equation}
The necessary optimality conditions for the solution $\oz$ to $(MP)$ corresponding to \eqref{disc-opt} in $(P^k)$ follow from \cite[Proposition~6.4(ii) and
Theorem~6.5(ii)]{m-book}. Theorem~\ref{str-conver} tells us that the inequality constraints in \eqref{e:5.13*} and \eqref{e:5.14*} are inactive for $k$ sufficiently large (suppose that for all $k\in\N$ without loss of generality), and hence the associated multipliers collapse to zero. In this way we find $\lm^k\ge 0,\;\xi^k=(\xi_1^k,\ldots,\xi_m^k)\in\R^m_+,\;\al^{1k}=(\al_{0}^{1k},\ldots,\al_{\nu(k)}^{1k})\in\R^{\nu(k)+1}_+,\;\al^{2k}=(\al_{0}^{2k},\ldots,\al_{\nu(k)}^{2k})\in\R^{\nu(k)+1}_{-},\;p_j^k=
\(p_{j}^{xk},p_{j}^{ak},p_{j}^{bk}\)\in\R^{n+mn+m}$ for $j=1,\ldots,\nu(k)$, and $z^*_j=(x^*_{0j},\ldots,x^*_{\nu(k)j},a^*_{0j},\ldots,a^*_{\nu(k)j},b^*_{0j},\ldots,b^*_{\nu(k)j},u^*_{0j},\ldots,u^*_{(\nu(k)-1)j},X^*_{0j},\ldots,X^*_{(\nu(k)-1)j},A^*_{0j},\ldots,
A^*_{(\nu(k)-1)j},B^*_{0j},\ldots,\\B^*_{(\nu(k)-1)j})$ for $j=0,\ldots,\nu(k)$, not zero all together, such that
\begin{equation}\label{69}
z^*_j\in N(\oz,\Xi_j)+N(\oz;\O_j^u)+\Big[N(\oz;\O_{j}^a)+N(\oz;\O_{j}^b)\Big]\;\mbox{ for }\;j=0,\ldots,\nu(k)-1,
\end{equation}
\begin{equation}\label{70}
-z^*_0-\ldots-z^*_{\nu(k)}\in\lm^k\partial\phi_0(\oz)+\sum_{i=1}^{m}
\xi_i^k\nabla
q_i(\oz)+\sum_{j=0}^{\nu(k)}\sum^m_{i=1}\al_{ij}^{1k}\nabla
l^1_{ij}(\oz)+\sum_{j=0}^{\nu(k)}\sum^m_{i=1}\al_{ij}^{2k}\nabla
l^2_{ij}(\oz)+\sum_{j=0}^{\nu(k)-1}\nabla g_j(\oz)^*p_{j+1}^{k},
\end{equation}
\begin{equation}\label{71+}
\xi_i^k q_i(\oz)=0,\;\;i=1,\ldots,m,
\end{equation}
\begin{equation}\label{71l1}
\al_{ij}^{1k}\(\n a_{ij}^k\en-(1+\dd_k)\)=0,\;i=1,\ldots,m,\;\;j=0,\ldots,\nu(k),
\end{equation}
\begin{equation}\label{71l2}
\al_{ij}^{2k}\(\n a_{ij}^k\en-(1-\dd_k)\)=0,\;i=1,\ldots,m,\;j=0,\ldots,\nu(k).
\end{equation}

The inclusion in \eqref{69} requires clarification. Indeed, it follows from the aforementioned necessary optimality conditions for $(MP)$ taken from \cite{m-book} that
\begin{equation}\label{69a}
z^*_j\in N\big(\oz;\Xi_j\cap\Th_j\cap\O^u_j\big),\;\mbox{ where }\;
\Th_j:= \O^a_{j}\cap\O^b_{j} \;\mbox{ for }\;j=0,\ldots,\nu(k)-1.
\end{equation}
To pass from \eqref{69a} to \eqref{69}, we need to use the intersection formula for the normal cone \eqref{nor_con}. Fix any $j\in\{0,\ldots,\nu(k)-1\}$ and denote $\Lambda_j:=\Th_j\cap\O^u_j$. The intersection rule from \cite[Theorem~2.16]{m18} reads as
\begin{equation}\label{69b}
N\big(\oz;\Xi_j\cap\Lambda_j\big)\subset N(\oz;\Xi_j)+N(\oz;\Lambda_j)\;\mbox{ if }\;N(\oz;\Xi_j)\cap\big(-N(\oz;\Lambda_j)\big)=\{0\},
\end{equation}
and thus we should check the fulfillment of the qualification condition therein in order to use \eqref{69b} in \eqref{69a}. To proceed, pick any $\Tilde z^*_j\in N(\oz;\Xi_j)\cap(-N(\oz;\Lambda_j))$ and observe by \eqref{e:5.16*} that its diagonal components satisfy
\begin{equation}\label{intersec1}
\begin{array}{ll}
&(\Tilde x^*_{jj},\Tilde a^*_{jj},\Tilde b^*_{jj}, \Tilde u^*_{jj},-\Tilde X^*_{jj}, \Tilde A^*_{jj},\Tilde B^*_{jj})\in N\Big(\Big(\ox_j^k,\oa^k_j,\ob^k_k,\ou^k_j,\disp-\frac{\ox^{k+1}_j-\ox^{k+1}_j}{h_k}\Big);\gph F\Big) \times \{0\} \times \{0\},\\\\
&-(\Tilde x^*_{jj},\Tilde a^*_{jj},\Tilde b^*_{jj}, \Tilde X^*_{jj},\Tilde u^*_{jj}, \Tilde A^*_{jj},\Tilde B^*_{jj})\in \{0\}\times\{0\}\times\{0\}\times\{0\}\times N\big((\ou^k_j, \dot\oa^k_j, \dot\ob^k_j); U\times A\times B\big)
\end{array}
\end{equation}
for all  $j=0,\ldots,\nu(k)-1$ with its other components equal zero. We get from \eqref{intersec1} that $\Tilde x^*_{jj}=0, \Tilde a^*_{jj =0}$, $\Tilde b^*_{jj} = 0$, $\Tilde X^*_{jj}=0$, $\Tilde A^*_{jj}=0$, and $\Tilde B^*_{jj}=0$. It follows from the first inclusion in \eqref{intersec1} by \eqref{coderivative} that
\begin{equation*}
(0,0,0,\Tilde u^*_{jj})\in D^*F\Big(\ox^k_j,\oa^k_j,\ob^k_j,\ou^k_j,-\frac{\ox^{k+1}_j-\ox^k_j}{h_k}\Big)(0),\quad j=0,\ldots,\nu(k)-1.
\end{equation*}
Applying to the latter inclusion the coderivative evaluation from Proposition~\ref{Th7} under the imposed PLICQ condition tells us that $\Tilde u^*_{jj}$ is zero whenever $j=0,\ldots,\nu(k)-1$. This shows that $\Tilde z^*_j=0$ for all such $j$, and thus the intersection rule in \eqref{69b} holds.\vspace*{-0.03in}

To proceed further with verifying \eqref{69} for $j=0,\ldots,\nu(k)-1$, we first apply to $N(\oz;\Lambda_j)=N(\oz;\Th_j\cap\O^u_j)$ in \eqref{69b} the normal cone intersection rule from \cite[Theorem~2.16]{m18} and then use the result of \cite[Corollary~2.17]{m18} giving us the sum decomposition of the normal cone $N(\ox;\Th_j)$ to the intersection of finitely many sets that define $\Th_j$ in \eqref{69a}. Taking into account the structures of the sets in \eqref{e:5.16*}--\eqref{uinU2} and of the qualification conditions in \cite[Theorem~2.16 and Corollary~2.17]{m18}, for each index $i\in\{1,\ldots,m\}$ and $j\in\{0,\ldots,\nu(k)-1\}$ we find $\psi^{uk}_j$ satisfying the normal cone inclusions \eqref{dloc-max} and such that
\begin{equation}\label{e:5.18*}
\left\{\begin{array}{ll}
\bigg(x^*_{jj},a^*_{jj},b^*_{jj},u^*_{jj}-\psi^{uk}_{j},-X^*_{jj}\bigg)\in
N\Big(\Big(\ox_j^k,\oa_j^k,\ob_j^k,\ou_j^k,-\dfrac{\ox_{j+1}^k-\ox_j^k}{h^k_j}\Big);\gph F\Big),\\
(A^*_{jj}, B^*_{jj})\in N(\dot a^k_j;A)\times N(\dot b^k_j;B)
\end{array}
\right.
\end{equation}
while the other components of $z^*_j$ for $j=0,\ldots,\nu(k)$ are zero. Note that the failure of the qualification conditions ensuring the aforementioned intersection rules leads us to \eqref{e:5.18*} and the subsequent results of this theorem with $\|\psi^{uk}\|\ne 0$ in the first nontriviality condition in \eqref{ntc1}. Otherwise, we arrive at the claimed assertions of the theorem with the ``full" nontriviality \eqref{ntc1} as stated; see below. This gives us the representation
\begin{eqnarray*}
\begin{aligned}
-z^*_0-\ldots-z^*_{\nu(k)}&=\big(-x^*_{0\nu(k)}-x^*_{00},-x^*_{11},\ldots,-x^*_{\nu(k)-1,\nu(k)-1},0,-a^*_{0\nu(k)}-a^*_{00},-a^*_{11},\ldots,-a^*_{\nu(k)-1\,
\nu(k)-1},0,\\
&\quad -b^*_{0\nu(k)}-b^*_{00},-b^*_{11},\ldots,-b^*_{\nu(k)-1\,\nu(k)-1},0,-u^*_{0\nu(k)}-u^*_{00},\ldots,-u^*_{\nu(k)-1\,\nu(k)-1},\\
&\quad -X^*_{00},\ldots,-X^*_{\nu(k)-1\,\nu(k)-1}, -A^*_{00},\ldots,-A^*_{\nu(k)-1\,\nu(k)-1},-B^*_{00},\ldots,-B^*_{\nu(k)-1\,\nu(k)-1}\big).
\end{aligned}
\end{eqnarray*}
Observe further that the right-hand side of the inclusion in \eqref{70} is represented by
\begin{eqnarray*}
\lm^k\partial\phi_0(\oz)+\sum_{i=1}^{m}\xi_i^k\nabla
q_i(\oz)+\sum_{j=0}^{\nu(k)}\sum_{i=1}^m\al_{ij}^{1k}\nabla
l^1_{ij}(\oz)+\sum_{j=0}^{\nu(k)}\sum_{i=1}^m\al_{ij}^{2k}\nabla
l^2_{ij}(\oz)+\sum_{j=0}^{\nu(k)-1}\gra g_j(\oz)^*p_{j+1}^k
\end{eqnarray*}
under the complementary slackness conditions
\begin{eqnarray*}
\xi_i^k\(\la a_{i\nu(k)}^k,x_{\nu(k)}^k\ra-b_{i\nu(k)}^k\)=0,\quad i=1,\ldots,m.
\end{eqnarray*}
Defining now the quantities
\begin{equation*}
\rho_j(\oz):=\int_{t_j^k}^{t_{j+1}^k}\n\(\frac{\ox_{j+1}^k-\ox_j^k}{h^k_j}-\dot{\ox}(t),\frac{\oa_{j+1}^k-\oa_j^k}{h^k_j}-\dot{\oa}(t),\frac{\ob_{j+1}^k-\ob_j^k
}{h^k_j}-\dot{\ob}(t),\ou_j^k(t)-\ou(t)\)\en^2dt,
\end{equation*}
we come up with the following relationships:
$$\(\sum_{i=1}^{m}\xi_i^k\nabla
q_i(\oz)\)_{(x_{\nu(k)},a_{\nu(k)},b_{\nu(k)},u_{\nu(k)})}=\(\sum_{i=1}^{m}
\xi_i^k\oa_{ik}^k,\[\xi^k,\rep_m(\ox_{\nu(k)}^k)\],-\xi^k,0\),$$
$$\(\sum_{j=0}^{\nu(k)}\sum^m_{i=1}\al_{ij}^{1k}\nabla
l^1_{ij}(\oz)\)_{(a_j)}=2\[\al_{j}^{1k},\oa_{j}^k\],\;\;j=0,\ldots,\nu(k),$$
$$\(\sum_{j=0}^{\nu(k)}\sum^m_{i=1}\al_{ij}^{2k}\nabla
l^2_{ij}(\oz)\)_{(a_j)}=2\[\al_{j}^{2k},\oa_{j}^k\],\;\;j=0,\ldots,\nu(k),$$
$$\(\sum_{j=0}^{\nu(k)-1}\nabla g_j(\oz)^*p_{j+1}^k\)_{(x_j,a_j,b_j)}=\left\{
\begin{array}{llll}
-p_{1}^{k}&\textrm{ if }& j=0\\[1ex]
p_{j}^{k}-p_{j+1}^{k}&\textrm{ if }& j=1,\ldots,\nu(k)-1\\[1ex]
p_{\nu(k)}^{k}&\textrm{ if }& j=\nu(k)
\end{array}
\right.,
$$
\begin{eqnarray*}
\begin{aligned}
&\(\sum_{j=0}^{\nu(k)-1}\nabla
g_j(\oz)^*p_{j+1}^k\)_{(X_j,A_j,B_j)}=(-h^k_0p_{1}^{xk},-h^k_1p_{2}^{xk},\ldots,
-h^k_{\nu(k)-1}p_{\nu(k)}^{xk},\\
&\qquad \qquad -h^k_0p_{1}^{ak},-h^k_1p_{2}^{ak},\ldots,-h^k_{\nu(k)-1}p_{\nu(k)}^{ak},
-h^k_0p_{1}^{bk},-h^k_1p_{2}^{bk},\ldots,-h^k_{\nu(k)-1}p_{\nu(k)}^{bk}),
\end{aligned}
\end{eqnarray*}
$$
\partial\phi_0(\oz)\subset\partial\vph(\ox_{\nu(k)}^k)+\sum_{j=0}^{\nu(k)-1}
h^k_j\partial\ell\(\ox_j^k,\oa_j^k,\ob_j^k,\ou_j^k,\oX_j^k,\oA_j^k,
\oB_j^k\)+\frac{1}{2}\sum_{j=0}^{\nu(k)-1}\gra \rho_j(\oz).$$
Moreover, we represent the set $\lm^k\partial\phi_0(\oz)$ as the collection of vectors
$$
\begin{aligned}
&\lm^k(h^k_0w_{0}^{xk},h^k_1w_{1}^{xk},\ldots,h^k_{\nu(k)-1}kw_{\nu(k)-1}^{xk},
v_{\nu(k)}^k,h^k_0w_{0}^{ak},h^k_1w_{1}^{ak},\ldots,h^k_{\nu(k)-1}w_{\nu(k)-1}^{
ak},0,\\
&h^k_0w_{0}^{bk},h^k_1w_{1}^{bk},\ldots,h^k_{\nu(k)-1}w_{\nu(k)-1}^{bk},0,\th_{0
}^{uk}+h^k_0w_{0}^{uk},\th_{1}^{uk}+h^k_1w_{1}^{uk},\ldots,\th_{\nu(k)-1}^{uk}
+h^k_{\nu(k)-1}w_{\nu(k)-1}^{uk},\\
&\th_{0}^{xk}+h^k_0v_{0}^{xk},\th_{1}^{xk}+h^k_1v_{1}^{xk},\ldots,\th_{\nu(k)-1}
^{xk}+h^k_{\nu(k)-1}v_{\nu(k)-1}^{xk},\th_{0}^{ak}+h^k_0v_{0}^{ak},\th_{1}^{ak}
+h^k_1v_{1}^{ak},\ldots,
\th_{\nu(k)-1}^{ak}+h^k_{\nu(k)-1}v_{\nu(k)-1}^{ak},\\
&\th_{0}^{bk}+h^k_0v_{0}^{bk},\th_{1}^{bk}+h^k_1v_{1}^{bk},\ldots,\th_{\nu(k)-1}
^{bk}+h^k_{\nu(k)-1}v_{\nu(k)-1}^{bk}),
\end{aligned}
$$
where $v_{\nu(k)}^k\in\partial\vph(\ox_{\nu(k)}^k)$, where $(\th_{j}^{uk},\th_{j}^{xk},\th_{j}^{ak},\th_{j}^{bk})$ are taken from the formulation of the theorem, and where
\begin{equation*}
\(w_{j}^{xk},w_{j}^{ak},w_{j}^{bk},w_{j}^{uk},v_{j}^{xk},v_{j}^{ak},v_{j}^{bk}\)\in\partial\ell\(\ox_j^k,\oa_j^k,\ob_j^k,\ou_j^k,\frac{\ox_{j+1}^k-\ox_j^k}{h^k_j},\frac{\oa_{j+1}^k-\oa_j^k}{h^k_j},
\frac{\ob_{j+1}^k-\ob_j^k}{h^k_j}\)
\end{equation*}
for $j=0,\ldots,\nu(k)-1$. Unifying all of this yields the relationships
\begin{equation}\label{daux}
-x^*_{00}-x^*_{0\nu(k)}=\lm^kh^k_0w_{0}^{xk}-p_{1}^{xk},
\end{equation}
\begin{equation}\label{e:5.21*x}
-x^*_{jj}=\lm^kh^k_jw_{j}^{xk}+p_{j}^{xk}-p_{j+1}^{xk},\;j=1,\ldots,\nu(k)-1,
\end{equation}
\begin{equation}\label{e:5.24*x}
0=\lm^kv_{\nu(k)}^k+p_{\nu(k)}^{xk}+\sum_{i=1}^{m}\xi_i^k\oa_{ik}^k,\textrm{
where }\;v_{\nu(k)}^k\in\partial\vph (\ox_{\nu(k)}^{k}),
\end{equation}
\begin{equation}\label{daua}
-a^*_{00}-a^*_{0\nu(k)}=\lm^kh^k_0w_{0}^{ak}+2\[\al_{0}^{1k}+\al_{0}^{2k},\oa_{0}^k\]-p_{1}^{ak},
\end{equation}
\begin{equation}\label{e:5.21*a}
-a^*_{jj}=\lm^kh^k_jw_{j}^{ak}+2\[\al_{j}^{1k}+\al_{j}^{2k},\oa_{j}^k\]+p_{j}^{
ak}-p_{j+1}^{ak},\;j=1,\ldots,\nu(k)-1,
\end{equation}
\begin{equation}\label{e:5.24*a}
0=2\[\al_{\nu(k)}^{1k}+\al_{\nu(k)}^{2k},\oa_{\nu(k)}^k\]+p_{\nu(k)}^{ak}+\[
\xi^k,\rep_m(\ox_{\nu(k)}^k)\],
\end{equation}
\begin{equation}\label{daub}
-b^*_{00}-b^*_{0\nu(k)}=\lm^kh^k_0w_{0}^{bk}-p_{1}^{bk},
\end{equation}
\begin{equation}\label{e:5.21*b}
-b^*_{jj}=\lm^kh^k_jw_{j}^{bk}+p_{j}^{bk}-p_{j+1}^{bk},\;\;j=1,\ldots,\nu(k)-1,
\end{equation}
\begin{equation}\label{e:5.24*b}
0=p_{\nu(k)}^{bk}-\xi^k,
\end{equation}
\begin{equation}\label{dauu}
-u^*_{00}=\lm^k\th_{0}^{uk}+\lm^kh^k_0w_{0}^{uk},
\end{equation}
\begin{equation}\label{e:5.22*}
-u^*_{jj}=\lm^k\th_{j}^{uk}+\lm^kh^k_jw_{j}^{uk},\;\;j=1,\ldots,\nu(k)-1,
\end{equation}
\begin{equation}\label{e:5.23*X}
-X^*_{jj}=\lm^k\th_{j}^{xk}+\lm^kh^k_jv_{j}^{xk}-h^k_jp_{j+1}^{xk},\;\;j=0,\ldots,
\nu(k)-1,
\end{equation}
\begin{equation}\label{e:5.23*A}
-A^\ast_{jj}=\lm^k\th_{j}^{ak}+\lm^kh^k_jv_{j}^{ak}-h^k_jp_{j+1}^{ak},\;\;j=0,\ldots,\nu(k)-1,
\end{equation}
\begin{equation}\label{e:5.23*B}
-B^\ast_{jj}=\lm^k\th_{j}^{bk}+\lm^kh^k_jv_{j}^{bk}-h^k_jp_{j+1}^{bk},\;\;j=0,\ldots,
\nu(k)-1.
\end{equation}

To proceed now with verifying the necessary optimality conditions of the theorem, we deduce from \eqref{e:5.24*x}, \eqref{e:5.24*a}, and \eqref{e:5.24*b}, respectively, that
\begin{equation}\label{mutx}
-p_{\nu(k)}^{xk}=\lm^kv_{\nu(k)}^k+\sum_{i=1}^{m}\xi_i^k\oa_{i\nu(k)}
^k\in\lm^k\partial\vph(\ox_{\nu(k)}^k)+\sum_{i=1}^{m}\xi_i^k\oa_{i\nu(k)}^k,
\end{equation}
\begin{equation}\label{muta}
p_{\nu(k)}^{ak}=-2\[\al_{\nu(k)}^{1k}+\al_{\nu(k)}^{2k},\oa_{\nu(k)}^k\]-\[\xi^k,\rep_m(\ox_{\nu(k)}^k)\],
\end{equation}
\begin{equation}\label{mutb}
p_{\nu(k)}^{bk}=\xi^k.
\end{equation}
Next extend each vector $p^k$ by adding the zero component $p_0^k:=\(x^*_{0\nu(k)},a^*_{0\nu(k)},b^*_{0\nu(k)},u^*_{0\nu(k)}\)$. It then follows
from the relationships in \eqref{e:5.21*x}, \eqref{e:5.21*a}, \eqref{e:5.21*b},
\eqref{e:5.23*X}--\eqref{e:5.23*B} that
$$
\begin{aligned}
\frac{x^*_{jj}}{h^k_j}&=\frac{p_{j+1}^{xk}-p_{j}^{xk}}{h^k_j}-\lm^kw_{j}^{xk},\\
\frac{a^*_{jj}}{h^k_j}&=\frac{p_{j+1}^{ak}-p_{j}^{ak}}{h^k_j}-\lm^kw_{j}^{ak}
-\frac{2}{h^k_j}\[\al_{j}^{1k}+\al_{j}^{2k},\oa_{j}^k\],\\
\frac{b^*_{jj}}{h^k_j}&=\frac{p_{j+1}^{bk}-p_{j}^{bk}}{h^k_j}-\lm^kw_{j}^{bk},\\
\frac{u^*_{jj}}{h^k_j}&=-\frac{1}{h^k_j}\lm^k\th_{j}^{uk}-\lm^kw_{j}^{uk},\\
\frac{X^*_{jj}}{h^k_j}&=-\frac{1}{h^k_j}\lm^k\th_{j}^{xk}+p_{j+1}^{xk}-\lm^kv_{j
}^{xk},\\
\frac{A^*_{jj}}{h^k_j}&=-\frac{1}{h^k_j}\lm^k\th_{j}^{ak}+p_{j+1}^{ak}-\lm^kv_{j}^{ak},\\
\frac{B^*_{jj}}{h^k_j}&=-\frac{1}{h^k_j}\lm^k\th_{j}^{bk}+p_{j+1}^{bk}-\lm^kv_{j}^{bk}.
\end{aligned}
$$
Plugging these expressions into the left-hand side of \eqref{e:5.18*} and taking into account the equalities in \eqref{71+}--\eqref{71l2}, \eqref{e:5.24*x}, \eqref{e:5.24*a}, and \eqref{e:5.24*b} give us the conditions
\begin{equation}\label{xi}
\xi_i^k\(\la\oa_{ik}^k,\ox_{k}^k\ra-\ob_{ik}^k\)=0,\;\;i=1,\ldots,m,
\end{equation}
\begin{equation}\label{al1}
\al_{ij}^{1k}\(\|\oa_{ij}^k\|-(1+\dd_k)\)=0,\;\;i=1,\ldots,m,\;\;j=0,\ldots,
\nu(k),
\end{equation}
\begin{equation}\label{al2}
\al_{ij}^{2k}\(\|\oa_{ij}^k\|-(1-\dd_k)\)=0,\;i=1,\ldots,m,\;\;j=0,\ldots,\nu(k),
\end{equation}
\begin{equation}\label{e:5.10*}
\begin{aligned}
&\Bigg(\frac{p_{j+1}^{xk}-p_{j}^{xk}}{h^k_j}-\lm^kw_{j}^{xk},\frac{p_{j+1}^{ak}
-p_{j}^{ak}}{h^k_j}-\lm^kw_{j}^{ak},
\frac{p_{j+1}^{bk}-p_{j}^{bk}}{h^k_j}-\lm^kw_{j}^{bk},-\frac{1}{h^k_j}\lm^k\th_{j}
^{uk}-\lm^kw_{j}^{uk},-p_{j+1}^{xk}\\
&+\lm^k\Big(v_{j}^{xk}+\dfrac{1}{h^k_j}\th_{j}^{xk}\Big)\Bigg)\in\(0,\frac{2}{h^k_j}
\[\al_{j}^{1k}+\al_{j}^{2k},\oa^k_j\],0,\frac{1}{h^k_j}\psi^{uk}_j,0\)\\
&\qquad\qquad\qquad\qquad\qquad\quad +N\(\(\ox_j^k,\oa_j^k,\ob_j^k,\ou_j^k,-\frac{\ox_{j+1}^k-\ox_j^k}{h^k_j}\);\gph F\),\;\;j=0,\ldots,\nu(k)-1,
\end{aligned}
\end{equation}
and the vectors $(\psi_j^{uk},\psi_{j}^{ak},\psi_{j}^{ak})$ satisfying the normal cone inclusions in \eqref{dloc-max}. It readily follows from condition \eqref{e:5.10*} and the coderivative definition that
\begin{equation*}
\begin{aligned}
&\bigg(\dfrac{p_{j+1}^{xk}-p_{j}^{xk}}{h^k_j}-\lm^kw_{j}^{xk},\dfrac{p_{j+1}^{ak
}-p_{j}^{ak}}{h^k_j}-\lm^kw_{j}^{ak}-\dfrac{2}{h^k_j}\[\al_{j}^{1k}+\al_{j}^{2k},\oa_{j}^k\],\dfrac{p_{j+1}^{bk}-p_{j}^{bk}}
{h^k_j}-\lm^kw_{j}^{bk},\\
&\qquad -\dfrac{1}{h^k_j}\lm^k\th_{j}^{uk}-\lm^kw_{j}^{uk}
-\dfrac{1}{h^k_j}\psi^{uk}_{j}\bigg)\in
D^*F\(\ox_j^k,\oa_j^k,\ob_j^k,\ou_j^k,-\frac{\ox_{j+1}^k-\ox_j^k}{h^k_j}
\)\(-\frac{1}{h^k_j}\lm^k\th_{j}^{xk}-\lm^k v_{j}^{xk}+p_{j+1}^{xk}\)
\end{aligned}
\end{equation*}
for all $j=0,\ldots,\nu(k)-1$. Observing the inclusion
\begin{equation*}
-\frac{\ox_{j+1}^k-\ox_j^k}{h^k_j}+g(\ox_j^k,\ou_j^k)\in N\bigg(\ox_j^k;C\(\oa_j^k,\ob_j^k\)\bigg),\mbox{where }j=0,\ldots,\nu(k)-1,
\end{equation*}
and using the PLICQ property of the vectors $\nn\oa^k_i|\;i\in I(\ox^k,\oa^k,\ob^k)\hnn$ give us a unique vector $\eta_j^k\in\R^m_+$ such that for all $i=1,\ldots,m$, we have the conditions
\begin{equation*}\label{h:5.39}
\sum_{i=1}^m\eta_{ij}^k\oa_{ij}^k=-\frac{\ox_{j+1}^k-\ox_j^k}{h^k_j}+g(\ox_j^k,
\ou_j^k)\;\textrm{ with }\;\eta_{ij}^k\in
N_{\R_-}\Big(\la\oa_{ij}^k,\ox_j^k\ra-\ob_{ij}^k\Big),\;j=0,\ldots,\nu(k)-1,
\end{equation*}
which justify the implications in \eqref{87} and \eqref{eta}. Applying now the coderivative upper estimate obtained in Proposition~\ref{Th7} under the PLICQ assumption of the theorem verify the relationships
$$
\(\begin{matrix}
\dfrac{p_{j+1}^{xk}-p_{j}^{xk}}{h^k_j}-\lm^kw_{j}^{xk},\dfrac{p_{j+1}^{ak}-p_{j}^{ak}}{h^k_j}-\lm^kw_{j}^{ak}-\dfrac{2}{h^k_j}\[\al_{j}^{1k}+\al_{j}^{2k},\oa_{j}^k\],\dfrac{p_{j+1}^{bk}-p_{j}^{bk}}{h^k_j}
-\lm^kw_{j}^{bk},\\
-\dfrac{1}{h^k_j}\lm^k\th_{j}^{uk}-\lm^kw_{j}^{uk}-\dfrac{1}{h^k_j}\psi^{uk}_{j}
\end{matrix}\)
$$
$$
\in\(\begin{matrix}
\disp-\nabla g_x(\ox_j^k,\ou_j^k)^*\(-\frac{1}{h^k_j}\lm^k\th_{j}^{xk}-\lm^kv_{j}^{xk}+p_{j+1}^{xk}\)+\sum_{i=1}^{m}\gg_{ij}^k\oa_{ij}^k,\\
\[\gg_{j}^k,\rep_m(\ox_j^k)\]+\[\eta_j^k,\rep_m\(-\dfrac{1}{h^k_j}\lm^k\th_{j}^{xk}-\lm^k v_{j}^{xk}+p_{j+1}^{xk}\)\],\\
-\gg_{j}^k,\;\disp-\nabla g_u(\ox_j^k,\ou_j^k)^*\(-\frac{1}{h^k_j}\lm^k\th_{j}^{xk}-\lm^kv_{j}^{xk}+p_{j+1}^{xk}\)
\end{matrix}\)
$$
for $j=0,\ldots,\nu(k)-1$, where the components $\gg^k_{ij}$ of the vectors $\gg_j^k\in\R^m$ as $i=1,\ldots,m$ are taken from \eqref{congg}. The obtained relationships together with the direct calculation of the coderivative $D^*N_{\R_-}$ in \eqref{congg} ensure the fulfillment of all the conditions in \eqref{conx}--\eqref{conb} as well as the inclusion in \eqref{cony}.

Defining further $\eta_{\nu(k)}^k:=\xi^k$ via $\xi^k$ yields $\eta_j^k\in\R^m_+$ for $j=0,\ldots,\nu(k)$ and allows us to deduce the transversality conditions in \eqref{nmutx}--\eqref{nmutb} from those in \eqref{mutx}--\eqref{mutb}. It follows from \eqref{xi} and the definition of $\eta_{\nu(k)}^k$ that the second complementarity slackness condition in \eqref{94} holds. Observing that \eqref{96}  fulfills due to
\begin{equation*}
-\dfrac{1}{h^k_j}\lm^k\th_{j}^{xk}-\lm^kv_{j}^{xk}+p_{j+1}^{xk}\in\disp\bigcap_{\nn i|\;\eta^k_{ij}>0\hnn}(\oa^k_{ij})^\perp
\end{equation*}
by Proposition~\ref{Th7}, we deduce from \eqref{al1} and \eqref{al2} that both inclusions in \eqref{t:7.20} are satisfied. Furthermore, conditions \eqref{nmutx}, \eqref{nmuta}, and \eqref{nmutb} clearly follow from \eqref{mutx}--\eqref{mutb} due to \eqref{xi}.

Next we show that the general nontriviality condition in \eqref{ntc0} holds. Suppose on the contrary that $\lm^k=0,\,\xi^k=0,\,\al^{1k}+\al^{2k}=0,\,p_{j}^{xk}=0,\,p_{0}^{ak}=0,\,p_{0}^{bk}=0,\,\psi^{ak}_j=0,\,\psi^{bk}_j=0,\,\psi^{uk}_j=0$ for all $j=0,\ldots,\nu(k)-1$, which yields in turn $x^*_{0\nu(k)}=p_{0}^{xk}=0,\,a^*_{0\nu(k)}=p_{0}^{ak}=0$, and $b^*_{0\nu(k)}=p_{0}^{bk}=0$. Then we deduce from \eqref{e:5.24*x}, \eqref{e:5.24*a}, and \eqref{e:5.24*b} that $\(p_{\nu(k)}^{xk},p_{\nu(k)}^{ak},p_{\nu(k)}^{bk}\)=0$, and hence  that $\(p_{j}^{xk},p_{j}^{ak},p_{j}^{bk}\)=0$ for all $j=0,\ldots,\nu(k)$. This contradicts therefore the nontriviality condition in the equivalent problem $(MP)$.

It remains to verify the enhanced nontriviality condition \eqref{ntc1}. Arguing by contraposition, suppose that $\lm^k=0,\;\al^{1k}+\al^{2k}=0$, $\gg^k=0$, $p^{ak}_{\nu(k)}=0$, and $p^{bk}_{\nu(k)}=0$. Hence $\eta^k_{\nu(k)}=p^{bk}_{\nu(k)}=0$ due to \eqref{mutb} and the definition of $\eta^k_{\nu(k)}$. Then it follows from \eqref{nmutx} that $p^{xk}_{\nu(k)}=0$. Thus $\(p^{xk}_j,p^{ak}_j,p^{bk}_j\)=(0,0,0)$ for all $j=0,\ldots,\nu(k)-1$ by \eqref{conx}--\eqref{conb} and by $p^{xk}_{\nu(k)}=0$, $p^{ak}_{\nu(k)}=0$, $p^{bk}_{\nu(k)}=0$. This implies that $\psi^{uk}_j=0,\;\psi^{Ak}_j=0,\;\psi^{Bk}_j=0$, $j=0,\ldots,\nu(k)-1$, by \eqref{cony} and \eqref{dloc-max}. It means that \eqref{ntc0} is violated, which is a contradiction that justifies \eqref{ntc1} and thus completes the proof of the theorem. $\h$

We conclude this section by showing that the normal cone inclusions in \eqref{dloc-max} yield certain {\em maximization conditions} for optimal controls in $(P^k)$ under additional assumptions.\vspace*{-0.07in}

\begin{corollary}[\bf discrete maximization conditions]\label{cor}  In addition to the assumptions of Theorem~{\rm\ref{Th1}}, suppose that the normal cones $N(\ou^k_j;U)$, $N(\dot{\oa}^k_j;A)$ and  $N(\dot{\ob}^k_j;B)$ in \eqref{dloc-max} are tangentially generated, i.e., they are dual/polar to some tangent sets $N(\ou^k_j;U)=T^*(\ou^k_j;U)$, $N(\dot{\oa}^k_j;A)=T^*(\dot{\oa}^k_j;A)$, and $N(\dot{\ob}^k_j;B)$ for all $j=0,\ldots,\nu(k)-1$. Then for such $j$ the following {\sc local maximization conditions} hold:
\begin{equation}\label{lmp}
\langle\psi^{uk}_j,\ou^k_j\rangle=\max_{u\in T(\ou^k_j;U)}\la\psi^{uk}_j,u\ra,\;\;\la\psi^{ak}_{j},\dot{\oa}^k_{j}\ra=\max_{v\in T(\dot{a}^k_{j};A)}\la\psi^{ak}_{j},v\ra,\;\;\la\psi^{bk}_{j},\dot{\ob}^k_{j}\ra=\max_{v\in T(\dot{\ob}^k_{j};B)}\la\psi^{bk}_{j},v\ra.
\end{equation}
If furthermore the sets $U$, $A$, and $B$ are convex, then we have the {\sc global maximization conditions}
\begin{equation}\label{gmp}
\langle\psi^{uk}_j,\ou^k_j\rangle=\max_{u\in U}\la\psi^{uk}_j,u\ra,\;\;\la\psi^{ak}_{j},\dot{\oa}^k_{j}\ra=\max_{v\in A}\la\psi^{ak}_{j},v\ra,\;\;\la\psi^{bk}_{j},\dot{\ob}^k_{j}\ra=\max_{v\in B}\la\psi^{bk}_{j},v\ra.
\end{equation}
\end{corollary}\vspace*{-0.1in}
{\bf Proof}. The local maximization conditions in \eqref{lmp} follow from \eqref{dloc-max} due to the assumed normal-tangent duality. The convexity of the sets $U$, $A$, and $B$ yields the global maximization in \eqref{gmp}, since the limiting normal cone \eqref{nor_con} reduces to the normal cone \eqref{NC} for convex sets. $\h$\vspace*{-0.15in}

\section{Optimality Conditions for Controlled Sweeping Processes}\label{sec:sweep}
\setcounter{equation}{0}\vspace*{-0.1in}

In this section we establish necessary optimality conditions for local minimizers of the original sweeping optimal control problem $(P)$ by passing to the limit as $k\to\infty$ in the optimality conditions for the discrete-time problems $(P^k)$ obtained in Theorem~\ref{Th1}. The limiting procedure employs the well-posedness of discrete approximations presented in Theorem~\ref{str-conver} together with the advanced tools and results of variational analysis and generalized differentiation discussed in Section~\ref{sec:disc} including the second-order calculations taken from Proposition~\ref{Th7}.

Here is the main result of this paper that provides necessary optimality conditions for a given r.i.l.m.\ in problem $(P)$ entirely in terms of its initial data. Taking into account that we do not use below the subdifferentiation of the running cost with respect to the time variable, suppose for simplicity that it does not depend on $t$.\vspace*{-0.1in}

\begin{theorem}[\bf necessary conditions for relaxed intermediate local minimizers]\label{Th2} Let $(\ox(\cdot),\oa(\cdot),\ob(\cdot),\ou(\cdot))$ be an r.i.l.m.\ for problem $(P)$, where the functions $\dot{\ox}(\cdot),\dot{\oa}(\cdot),\dot{\ob}(\cdot)$, and $\ou(\cdot)$ are of bounded variation on $[0,T]$, and where the PLICQ property \eqref{PLICQ} is satisfies along this local minimizer for all $t\in[0,T]$. Suppose in addition to the standing assumptions {\rm(H1)--(H4)} that the terminal cost $\vph$ is locally Lipschitzian around $\ox(T)$, that the perturbation mapping $g$ is ${\cal C}^1$-smooth around $(\ox(t),\ou(t))$ as $t\in[0,T]$, and that the running cost $\ell$ admits the representation
\begin{equation}\label{run_cost}
\ell(x,a,b,u,\dot{x},\dot{a},\dot{b})=\ell_1(x,a,b,u,\dot{x})+\ell_2(\dot{a})+\ell_3(\dot{b}),
\end{equation}
where $\ell$ is locally Lipschitzian with constant $L_\ell$ while $\ell_2$ and $\ell_3$ are ${\cal C}^1$-smooth in $\dot a$ and $\dot b$, respectively, around the given local minimizer. Then there exist $\lm\ge0$, adjoint arcs $p(\cdot)=(p^x(\cdot),p^a(\cdot),p^b(\cdot))\in W^{1,2}([0,T];\R^{n+mn+m})$ and $q(\cdot)=(q^x(\cdot),q^a(\cdot),q^b(\cdot))\in BV([0,T];\R^{n+mn+m})$, signed measures $\gg=(\gg_1,\ldots,\gg_m)\in{\cal C}^*([0,T];\R^m)$ and $\al\in{\cal C}^*([0,T];\R)$, as well as the subgradient functions
\begin{equation*}
\big(w^x(\cdot),w^a(\cdot),w^b(\cdot),w^u(\cdot),v^x(\cdot),v^a(\cdot),v^b(\cdot)\big)\in L^2\big([0,T];\R^{n+mn+m+d}\big)\times L^2\big([0,T];\R^{n+mn+m}\big)
\end{equation*}
satisfying for a.e.\ $t\in[0,T]$ the subdifferential inclusion
\begin{equation}\label{wv}
\big(w^x(t),w^a(t),w^b(t),w^u(t),v^x(t),v^a(t),v^b(t)\big)\in\co\partial\ell\big(\ox(t),\oa(t),\ob(t),\ou(t),\dot{\ox}(t),\dot{\oa}(t),\dot{\ob}(t)\big)
\end{equation}
such that we have the following necessary optimality conditions:
\begin{itemize}
\item[{\bf(i)}] The {\sc primal-dual dynamic relationships} consisting of:\\[1ex]
$\bullet$ The {\sc primal arc representation}
\begin{equation}\label{37}
-\dot{\ox}(t)=\sum_{i=1}^m\eta_i(t)\oa_i(t)-g\big(\ox(t),\ou(t)\big)\;\textrm{ for a.e. }\;t\in[0,T),
\end{equation}
where the functions $\eta_i(\cdot)\in L^2([0,T);\R_+)$ are uniquely determined by \eqref{37} provided the fulfillment of LICQ along the given local minimizer for all $t\in[0,T)$. Furthermore, for a.e.\ $t\in[0,T)$ and all $i=1,\ldots,m$ we have  in the latter case the {\sc dynamic complementary slackness conditions}
\begin{equation}\label{41}
\begin{array}{ll}
\la\oa_i(t),\ox(t)\ra<\ob_i(t)\sr\eta_i(t)=0\;\mbox{ and }\;\eta_i(t)>0\sr\la\oa_i(t),q^x(t)-\lm v^x(t)\ra=0
\end{array}
\end{equation}
$\bullet$ The {\sc adjoint dynamic systems}
\begin{equation}\label{c:6.6}
\left\{\begin{array}{ll}
\qquad\qquad\big(\dot{p}^x(t),\dot{p}^a(t),\dot{p}^b(t)\big)=\lm\big(w^x(t),w^a(t),w^b(t)\big)+\\\\
\Big(-\nabla_x g\big(\ox(t),\ou(t)\big)^*\big(-\lm v^x(t)+q^x(t)\big),\[\eta(t),\rep_m\(-\lm v^{x}(t)+q^{x}(t)\)\],0\Big)
\end{array}\right.
\end{equation}
for a.e.\ $t\in[0,T]$, where the right continuous representative of $q(t)$ is given, for all $t\in[0,T]$ in the form
\begin{equation}\label{c:6.9}
q(t):=p(t)-\int_{[t,T]}\(\sum_{i=1}^m\oa_i(s)d\gg_i(s),2\[\oa(s),d\al(s)\]+\[\rep_m(\ox(s),d\gg(s))\],-d\gg(s)\).
\end{equation}
$\bullet$ The {\sc normal cone adjoint inclusions} for control components: for a.e.\ $t\in[0,T]$ we have
\begin{equation}\label{c:6.6'}
\psi^u(t):=-\lm w^{u}(t)+\nabla_u g\big(\ox(t),\ou(t)\big)^*\big(-\lm v^x(t)+q^x(t)\big)\in
{\rm co}\,N\big(\ou(t);U\big),
\end{equation}
\begin{equation}\label{c:6.6''}
\psi^A(t):=q^a(t)-\lm \nabla\ell_2\big(\dot{\oa}(t)\big)\in{\rm co}\,N\big(\dot{\oa}(t);A\big),\;\mbox{ and }\;\psi^B(t):=q^b(t)-\lm \nabla\ell_3\big(\dot{\ob}(t)\big)\in{\rm co}\,N\big(\dot{\bar{b}}(t);B\big).
\end{equation}
$\bullet$ The {\sc local and global maximum principles}: for a.e.\ $t\in[0,T]$ we have\\[1ex]
---The {\sc local maximization conditions}: assuming that the normal cones in \eqref{c:6.6'}, \eqref{c:6.6''} are tangentially generated as $N(\ou(t);U)=T^*(\ou(t);U)$,
$N(\oa(t);A)=T^*(\oa(t);A)$, and $N(\ob(t);B)=T^*(\ob(t);B)$ yields
\begin{equation}\label{lmax}
\begin{array}{ll}
\la\psi^u(t),\ou(t)\ra=&\disp\max_{u\in T(\ou(t);U)}\la\psi^u(t),u\ra,\;\; \la\psi^A(t),\dot\oa(t)\ra=\disp\max_{v\in T(\dot\oa(t);A)}\la\psi^A(t),v\ra,\\
\la\psi^B(t),\dot\ob(t)\ra=&\disp\max_{v\in T(\dot\ob(t);B)}\la\psi^B(t),v\ra.
\end{array}
\end{equation}
---If the sets $U$, $A$, and $B$ are convex, then the {\sc global maximization conditions} hold:
\begin{equation}\label{max}
\la\psi^u(t),\ou(t)\ra=\max_{u\in U}\la\psi^u(t),u\ra,\;\;\la\psi^A(t),\dot\oa(t)\ra=\max_{v\in A}\la\psi^A(t),v\ra,\;\;\la\psi^B(t),\dot\ob(t)\ra=\max_{v\in B}\la\psi^B(t),v\ra.
\end{equation}
\item[{\bf (ii)}] The {\sc transversality} and {\sc endpoint complementary slackness conditions}: there exists a vector $\eta(T)=(\eta_1(T),\ldots,\eta_m(T))\in\R^m_+$ such that we have
\begin{equation}\label{42x}
-p^{x}(T)\in\lm\partial\vph\big(\ox(T)\big)+\sum_{i\in I(\ox(T),\oa(T),\ob(T))}\eta_i(T)\oa_i(T)\;\mbox{ with}
\end{equation}
\begin{equation}\label{42c}
\sum_{i\in I(\ox(T),\oa(T),\ob(T))}\eta_i(T)\oa_i(T)\in N\big(\ox(T);C(T)\big),
\end{equation}
\begin{equation}\label{42a}
p^{a}(T)+\big[\eta(T),\rep_m(\ox(T))\big]\in\oa(T)\R,
\end{equation}
\begin{equation}\label{42b}
p^{b}(T)=\eta(T),\;\mbox{ and}\;\la\oa_{i}(T),\ox(T)\ra<\ob_i(T)\sr p^b_{i}(T)=0\;\textrm{ for all }\;i=1,\ldots,m,
\end{equation}
where the collection of active constraint indices $I(\ox(T),\oa(T),\ob(T))$ is taken from \eqref{e:AI}.
\item[{\bf (iii)}] The {\sc measure nonatomicity conditions:}\\
{\bf(a)} If $t\in[0,T)$ and $\la\oa_i(t),\ox(t)\ra<\ob_i(t)$ for all $i=1,\ldots,m$, then there exists a neighborhood $V_t$ of $t$ in $[0,T]$ such that $\gg(V)=0$ for all Borel subsets $V$ of $V_t$.\\
{\bf(b)} Take $t\in[0,T]$ with $\|\oa(t)\|=1$. Then there is a neighborhood $W_t$ of $t$ in $[0,T]$ such that $\al(W)=0$ for all the Borel subsets $W$ of $W_t$.
\item[{\bf (iv)}] The {\sc measured coderivative condition:} Considering the $t$-dependent outer limit
$$
\begin{aligned}
\Limsup_{|S|\to 0,\Tilde\gg\to\gg}\Tilde\gg(S)(t)&:=\Big\{y\in\R^m\Big|\,exists\;\;\mbox{sequences}\;\;S_k\subset[0,T]\;\;\mbox{with}\;\;t\in S_k,\,|S_k|\to 0,\\
&\mbox{and measures}\;\;\gg^k_{mes}\;\;\mbox{such that}\;\;\gg^k\st{w^*}\longrightarrow\gg\;\;\mbox{and}\;\;\gg^k_{mes}(S_k)\to y\;\;\mbox{as}\;\;k\to\infty\Big\}
\end{aligned}
$$
over Borel subsets $S\subset[0,T]$ with the Lebesgue measure $|S|$, we have
\begin{equation}\label{mcc}
D^*N\Big(\big(\big\la\oa_i(t),\ox(t)\big\ra-\ob_i(t),\eta_i(t)\big);\R_-\Big)\Big(\big\la\oa_i(t),q^x(t)-\lm v^x(t)\big\ra\Big)\cap\Limsup_{|S|\to 0,\,\Tilde\gg_i\to\gg_i}\Tilde\gg_i(S)(t)\not=\emp,
\end{equation}
for a.e.\ $t\in[0,T]$ and all $i=1,\ldots,m$.
\item[{\bf(v)}] The {\sc nontriviality conditions}. The following hold:\\[1ex]
$\bullet$ We always have the {\sc general nontriviality condition}
\begin{equation}\label{gen-nontr}
\big(\lm,p(\cdot),q(\cdot)\big)\ne 0.
\end{equation}
$\bullet$ If in addition the LICQ assumption is imposed along the local minimizer for all $t\in[0,T)$, and if the dynamic interiority assumption
\begin{equation}\label{inter0}
\big\la\oa_i(t),\ox(t)\big\ra<\ob_i(t)\;\mbox{ for all }\;t\in[0,T)\;\mbox{ and }\;i=1,\ldots,m
\end{equation}
is satisfied, then we have {\sc enhanced nontriviality condition}
\begin{equation}\label{enh0}
\lambda+\|p(T)\|+\sup_{t\in[0,T]}\|\psi^A(t)\|+\sup_{t\in[0,T]}\|\psi^B(t)\|\ne 0.
\end{equation}
$\bullet$ The latter is reduced to the more precise one
\begin{equation}\label{enh1}
\lm+\|p^a(T)\|\ne 0\;\mbox{ whenever }\;p^a(T)\in\oa(T)\R
\end{equation}
provided the fulfillment of the additional interiority assumptions
\begin{equation}\label{inter}
\big\la\oa_i(T),\ox(T)\big\ra<\ob_i(T)\;\mbox{ for all }\;i=1,\ldots,m,
\end{equation}
\begin{equation}\label{psi-con}
\big(\dot\oa(t),\dot\ob(t)\big)\in\mbox{\rm{int}}(A\times B)\;\mbox{ for all }\;t\in [0,T].
\end{equation}
\end{itemize}
\end{theorem}\vspace*{-0.1in}
{\bf Proof}. Theorem~\ref{str-conver} tells us that the prescribed r.i.l.m.\ $(\ox(\cdot),\oa(\cdot),\ob(\cdot),\ou(\cdot))$ for problem $(P)$ is approximated by a sequence of extended optimal solutions $\(\ox^k(\cdot),\oa^k(\cdot),\ob^k(\cdot),\ou^k(\cdot)\)$ to the discrete-time problems $(P^k)$ in the norm topology of $W^{1,2}([0,T];\R^{n+mn+m})\times L^2([0,T];\R^d)$. Then Theorem~\ref{Th1}
gives us necessary optimality conditions that the latter quadruple satisfies for any fixed index $k$. We now pass to the limit as $k\to\infty$ to verify the necessary optimality conditions for the local minimizer $(\ox(\cdot),\oa(\cdot),\ob(\cdot),\ou(\cdot))$ listed in this theorem. For the reader's convenience, we split the proof into the following nine steps.\\[1ex]
{\bf Step~1}: {\em Verification of the subdifferential inclusion}. To justify \eqref{wv}, for each $k\in\N$ denote the piecewise linear extensions of $w^k_j$ and $v^k_j$ to $[0,T]$ by $w^k(\cdot)$ and $v^k(\cdot)$. It follows from the Lipschitz continuity of $\ell$ and relation \eqref{subcol} that  $\(w^k(\cdot),v^k(\cdot)\)$ is bounded in $L^2\big([0,T];\R^{2(n+mn+m)+d}\big)$. Thus we get the weak convergence (without relabeling) of $\(w^k(\cdot),v^k(\cdot)\)$ to some $\(w(\cdot),v(\cdot)\)\in L^2$ due to the weak compactness of bounded sets in this space. Then employing the Mazur weak closure theorem gives us the strong $L^2$-convergence to $\(w(\cdot),v(\cdot)\)$ of a sequence of convex combinations of $\(w^k(\cdot),v^k(\cdot)\)$ and hence its a.e.\ convergence (along a subsequence) on $[0,T]$. This readily justifies \eqref{wv} by taking into account the {\em robustness} (closed-graph property) of the subdifferential mapping \eqref{sub} with respect to perturbations of the initial point.\\[1ex]
{\bf Step~2:} {\em Verification of the primal equation and the first dynamic complementary condition.} Define the functions
\begin{equation*}
\th^{xk}(t):=\dfrac{1}{h^k_j}\th_{j}^{xk},\;\th^{ak}(t):=\dfrac{1}{h^k_j}\th_{j}^{ak},\;\th^{bk}(t):=\dfrac{1}{h^k_j}\th_{j}^{bk},\;\mbox{ and }\;\th^{uk}(t):=\dfrac{1}{h^k_j}\th_{j}^{uk}
\end{equation*}
for $t\in[t^k_j,t_{j+1}^k),\;j=0,\ldots,\nu(k)-1$ and $k\in\N$, where the quadruple $(\th_{j}^{uk},\th_{j}^{xk},\th_{j}^{ak},\th_{j}^{bk})$ is taken from Theorem~\ref{Th1}. This gives us as $k\to\infty$ the estimates
\begin{equation}
\label{th-est}
\begin{aligned}
\int_0^T\n\th^{xk}(t)\en^2dt&=\sum_{j=0}^{\nu(k)-1}\frac{\Big\|\th_{j}^{xk}\Big\|^2}{h^k_j}\le\frac{1}{h^k_j}\sum_{j=0}^{\nu(k)-1}\Big(\int_{t_j^k}^{t_{j+1}^k}\n\dot{\ox}(t)-\frac{\ox_{j+1}^k
-\ox_j^k}{h^k_j}\en dt\Big)^2\\
&\le\sum_{j=0}^{\nu(k)-1}\int_{t_j^k}^{t_{j+1}^k}\Big\|\dot{\ox}(t)-\frac{\ox_{j+1}^k-\ox_j^k}{h^k_j}\Big\|^2dt=\int_0^T\n\dot{\ox}(t)-\dot{\ox}^k(t)\en^2dt\to 0
\end{aligned}
\end{equation}
and shows similarly that $\int_0^T\n\th^{ak}(t)\en^2dt\to 0$ and $\int_0^T\n\th^{bk}(t)\en^2dt\to 0$. Furthermore, we get
\begin{eqnarray*}
\int_0^T\Big\|\th^{uk}(t)\Big\|^2dt&=&\sum_{j=0}^{\nu(k)-1}\frac{\Big\|\th^{uk}_j\Big\|^2}{h^k_j}\le\frac{1}{h^k_j}\sum_{j=0}^{\nu(k)-1}\Big(\int_{t_j^k}^{t_{j+1}^k}\n\ou_j^k-\ou(t)\en dt\Big)^2\\
&\le&\sum_{j=0}^{\nu(k)-1}\int_{t_j^k}^{t_{j+1}^k}\n\ou_j^k-\ou(t)\en^2dt=\int_0^T\n\ou^k(t)-\ou(t)\en^2dt\to 0
\end{eqnarray*}
as $k\to\infty$ due to the strong $W^{1,2}\times W^{1,2}\times W^{1,2}\times L^2$-convergence of $(\ox^k(\cdot),\oa^k(\cdot),\ob^k(\cdot)\ou^k(\cdot))$ to $(\ox(\cdot),\oa(\cdot),\ob(\cdot),\ou(\cdot))$. This verifies the strong convergence of $\(\th^{xk}(\cdot),\th^{ak}(\cdot),\th^{bk}(\cdot),\th^{uk}(\cdot)\)$ to zero and hence its a.e.\ convergence to zero along some subsequence (without relabeling).

Next we define the piecewise constant function $\eta^k(\cdot)$ on $[0,T)$ by $\eta^k(t):=\eta_j^k$ for $t\in[t^k_j,t_{j+1}^k)$, where $\eta_j^k$ are taken from Theorem~\ref{Th1}. It follows from \eqref{87} that
\begin{equation}\label{c:51}
-\dot{\ox}^k(t)=\sum_{i=1}^m\eta_{i}^k(t)\oa_i^k(t)-g\big(\ox^k(t),\ou^k(t)\big)\;\textrm{ whenever }\;t\in(t^k_j,t_{j+1}^k),\quad k\in\N.
\end{equation}
Moreover, the feasibility of $\(\ox(\cdot),\oa(\cdot),\ob(\cdot),\ou(\cdot)\)$ in $(P)$ tells us that $-\dot{\ox}(t)\in G(\ox(t),\oa(t),\ob(t))-g(\ox(t),\ou(t))$ for a.e.\ $t \in [0,T]$, where the normal cone mapping $G(\ox(t),\oa(t),\ob(t)):=N\(\ox(t);C(\oa(t),\ob(t))\)$ is closed-valued and measurable on $[0,T]$ by \cite[Theorem~14.26]{rw}. It follows now from the standard measurable selection results (see, e.g., \cite[Corollary~14.6]{rw}) that there exist nonnegative measurable functions $\eta_i(\cdot)$ on $[0,T]$ as $i=1,\ldots,m$ satisfying the primal arc representation \eqref{37}.

Let us further verify the first dynamic complementarity slackness condition in \eqref{41} under the LICQ condition imposed in this case. Combining $\eqref{c:51}$ and $\eqref{37}$ leads us to the equations
\begin{equation*}
\dot{\ox}(t)-\dot{\ox}^k(t)=\sum_{i=1}^m\big[\eta_i^k(t)\oa_i^k(t)-\eta_i(t)\oa_i(t)\big]+g\big(\ox(t),\ou(t)\big)-g\big(\ox^k(t),\ou^k(t)\big),
\end{equation*}
which are valid for all $t\in(t_j^k,t_{j+1}^k)$ and $j=0,\ldots,k-1$ and hence yield the estimate
\begin{equation}\label{c:52}
\n\sum_{i=1}^m\big[\eta_i^k(t)\oa_i^k(t)-\eta_i(t)\oa_i(t)\big]\en_{L^2}\le\big\|\dot{\ox}^k(t)-\dot{\ox}(t)\big\|_{L^2}+\big\|g\big(\ox(t),\ou(t)\big)-g\big(\ox^k(t),\ou^k(t)\big)\big\|_{L^2}
\end{equation}
on $(t^k_j,t^k_{j+1})$. Invoking then the strong convergence of $\(\ox^k(\cdot),\oa^k(\cdot),\ob(\cdot),\ou^k(\cdot)\)$ to $\(\ox(\cdot),\oa(\cdot),\ob(\cdot),\ou(\cdot)\)$ and taking into account the smoothness of $g$, we get the strong convergence of $\disp\sum_{i=1}^m\left[\eta_i(t)\oa_i(t)-\eta^k_i(t)\oa^k_i(t)\right]$ to zero in $L^2$ and thus its a.e.\ convergence on $[0,T]$ along some subsequence. On the other hand, it follows from \eqref{A5'} that
\begin{equation*}
\begin{aligned}
\int^T_0\[\eta^k_i(t)\]^2dt&\le\dfrac{1}{(1-\ve_k)^2}\int^T_0\sum_{i\in I(\oa^k(t),\ob^k(t),\ox^k(t))}\[\eta^k_i(t)\]^2\|\oa^k_i(t)\|^2dt\\
&\le 4\int^T_0\bigg[\sum_{i\in I(\oa^k(t),\ob^k(t),\ox^k(t))}\eta^k_i(t)\|\oa^k_i(t)\|\bigg]^2dt\\
&\le 4\sigma^2\int^T_0\bigg\|\sum_{i\in I(\oa^k(t),\ob^k(t),\ox^k(t))}\eta^k_i(t)\oa^k_i(t)\bigg\|^2dt\\
&=4\sigma^2\int^T_0\n-\dot\ox^k(t)+g\big(\ox^k(t),\ou^k(t)\big)\en^2dt\\
&\le 8\sigma^2\int^T_0\n\dot\ox^k(t)\en^2dt+\int^T_0\n g\big(\ox^k(t),\ou^k(t)\big)\en^2dt\le M
\end{aligned}
\end{equation*}
for some constant $M>0$, which justifies the boundedness of $\eta^k(\cdot)$ in $L^2$ due to the strong convergence of $\(\ox^k(\cdot),\oa^k(\cdot),\ob^k(\cdot),\ou^k(\cdot)\)$ to $\(\ox(\cdot),\oa(\cdot),\ob(\cdot),\ou(\cdot)\)$. It follows from the weak compactness of bounded sets in $L^2$ that there exists a function $\Tilde\eta(\cdot)$ such that a subsequence of $\{\eta^k(\cdot)\}$ weakly converges to $\Tilde\eta(\cdot)$. Employing the aforementioned Mazur theorem gives us a sequence of convex combinations of the functions from $\{\eta^k(\cdot)\}$, which converges to $\Tilde\eta(\cdot)$ strongly in $L^2$, and hence pointwise for a.e.\ $t\in[0,T)$ along a subsequence. Combining this with the a.e.\ convergence of $\disp\sum^m_{i=1}\eta^k_i(t)\oa^k_i(t)$ to $\disp\sum^m_{i=1}\eta_i(t)\oa_i(t)$ on $[0,T]$ and using the assumed linear independence of the vectors $\nn\oa_i(t)\;\big|\;i\in I(\ox(t),\oa(t),\ob(t))\hnn$ ensure that $\Tilde\eta(t)=\eta(t)$ and that $\eta^k(t)\to\eta(t)$ for a.e.\ $t\in[0,T)$. Invoking finally the strong convergence results from Theorem~\ref{str-conver} and the complementary slackness condition \eqref{eta} for the discrete problems $(P^k)$, we arrive at the first complementary slackness condition in \eqref{41}.\\[1ex]
{\bf Step~3:} {\em Continuous-time extensions of adjoint functions.} Necessary conditions of Theorem~\ref{Th2} for optimal solutions to discrete approximation problems $(P^k)$ involve discrete-time adjoint/dual functions defined at the mesh points. For passing to the limit in these conditions as $k\to\infty$, we need to properly extend the dual functions to the continuous-time interval $[0,T]$ for each $k\in\N$. Let us start with $q^k(\cdot)$ and define $q^k(t)=\(q^{xk}(t),q^{ak}(t),q^{bk}(t)\)$ on $[0,T)$ as the piecewise linear extensions of $q^k(t^k_j):=p^k_j$ as $j=0,\ldots,\nu(k)$. Then extend the functions $\gg^k(\cdot),\psi^k(\cdot)=\(\psi^{uk}(\cdot),\psi^{ak}(\cdot),\psi^{bk}(\cdot)\)$, $\al^{1k}(\cdot)$, and $\al^{2k}(\cdot)$ to $[0,T)$ by
\begin{equation}\label{c:6.25}
\begin{cases}
\gg^k(t):=\gg_j^k,\quad\psi^k(t):=\dfrac{1}{h^k_j}\psi_j^k,\\[1ex]
\al^{1k}(t):=\dfrac{1}{h^k_j}\al^{1k}_j,\quad\al^{2k}(t):=\dfrac{1}{h^k_j}\al^{2k}_j
\end{cases}
\end{equation}
for $t\in[t_j^k,t_{j+1}^k)$ and $j=0,\ldots,\nu(k)-1$ with $\gg^k(T):=0,\;\psi^k(T):=0,\;\al^{1k}(T):=\al^{1k}_{\nu(k)}$, and $\al^{2k}(T):=\al^{2k}_{\nu(k)}$. Defining further the auxiliary functions on the continuous-time interval
\begin{equation*}
\vth^k(t):=\max\big\{t^k_j\big|\;t^k_j\le t,\;0\le j\le\nu(k)\big\}\;\mbox{ for }\;t\in[0,T),
\end{equation*}
we deduce from \eqref{conx}--\eqref{cony} the following relationships for every $t\in(t_j^k,t_{j+1}^k)$ as $j=0,\ldots,\nu(k)-1$:
\begin{equation*}
\begin{aligned}
\dot{q}^{xk}(t)-\lm^k w^{xk}(t)&=-\nabla_xg\big(\ox^k(\vth^k(t)),\ou^k(\vth^k(t))\big)^*\big(-\lm^k\th^{xk}(t)-\lm^k v^{xk}\big(\vth^k(t)\big)+q^{xk}\big(\vth^k(t)+h_k\big)\big)\\
&\,\quad +\sum_{i=1}^m\gg^{k}_i(t)\oa^k_i\big(\vth^k(t)\big),
\end{aligned}
\end{equation*}
\begin{equation*}
\begin{aligned}
\dot{q}^{ak}(t)-\lm^k w^{ak}(t)&=2\big[\al^{1k}(t)+\al^{2k}(t),\oa^k\big(\vth^k(t)\big)\big]+\big[\gg^k(t),\rep_m\ox^k\big(\vth^k(t)\big)\big]\\
&\,\quad +\big[\eta^k_i(t),\rep_m\big(-\lm^k\th^{xk}(t)-\lm^k v^{xk}\big(\vth^k(t)\big)+q^{xk}\big(\vth^k(t)+h^k_j)\big)\big],\\
\dot{q}^{bk}(t)-\lm^k w^{bk}(t)&=-\gg^k(t),
\end{aligned}
\end{equation*}
\begin{equation}\label{cony'}
\begin{array}{ll}
&-\lm^k\th^{uk}(t)-\lm^kw^{uk}\big(\vth^k(t)\big)-\psi^{uk}(t)=\\
&\qquad\;\; -\nabla_ug\big(\ox^k\big(\vth^k(t)\big),\ou^k\big(\vth^k(t)\big)\big)^*\Big(-\lm^k\th^{xk}(t)-\lm^k v^{xk}\big(\vth^k(t)\big)+q^{xk}\big(\vth^k(t)+h_k\big)\Big).
\end{array}
\end{equation}
Next define the adjoint arcs $p^k(t)=\(p^{xk}(t),p^{ak}(t),p^{bk}(t)\)$ on $[0,T]$ by setting
\begin{equation}\label{pqg}
\begin{array}{ll}
p^k(t):=&q^k(t)+\disp\int_{[t,T]}\bigg(\sum_{i=1}^m\gg^{k}_i(s)\oa^k_i\big(\vth^k(s)\big),2\[\al^{1k}(s)+\al^{2k}(s),\oa^k\big(\vth^k(s)\big)\],\\
&\,+\disp\big[\gg^k(s),\rep_m\big(\ox^k(\vth^k(s))\big)\big],-\gg^k(s)\bigg)ds
\end{array}
\end{equation}
for all $t\in[0,T]$. This gives us $p^k(T)=q^k(T)$  and the differential equation
\begin{equation}\label{pk-qk}
\dot p^k(t)=\dot q^k(t)-\bigg(\sum_{i=1}^m\gg^{k}_i(t)\oa^k_i\big(\vth^k(t)\big),2\[\al^{1k}(t)+\al^{2k}(t),\oa^k\big(\vth^k(t)\big)\]+\big[\gg^k(t),\rep_m\big(\ox^k(\vth^k(t))\big)\big],-\gg^k(t)\bigg)
\end{equation}
for a.e.\ $t\in[0,T]$. As the consequence of the above, we arrive at the relationships
\begin{equation}\label{e:px-der}
\dot p^{xk}(t)-\lm^kw^{xk}(t)\in-\nabla_xg\big(\ox^k(\vth^k(t)\big),\ou^k\big(\vth^k(t)\big)\big)^*\(-\lm^k\th^{xk}(t)-\lm^k v^{xk}\big(\vth^k(t)\big)+q^{xk}\big(\vth^k(t)+h_k\big)\),
\end{equation}
\begin{equation}\label{e:pa-der}
\dot p^{ak}(t)-\lm^kw^{ak}(t)=\[\eta^k_i(t),\rep_m\(-\lm^k\th^{xk}(t)-\lm^k v^{xk}\big(\vth^k(t)\big)+q^{xk}\big(\vth^k(t)+h^k_j\big)\)\],
\end{equation}
\begin{equation}\label{e:pb-der}
\dot p^{bk}(t)-\lm^kw^{bk}(t)=0\;\mbox{ for all }\;t\in(t^k_j,t^k_{j+1})\;\mbox{ and }\;j=0,\ldots,\nu(k)-1.
\end{equation}
Next we define the measures $\gg^{k}_{mes}$, $\al^{1k}_{mes}$, and $\al^{2k}_{mes}$ on $[0,T]$ by
\begin{equation}\label{e:discrete-measures}
\left\{
\begin{array}{ll}
\disp\int_Sd\gg^{k}_{mes}:=\int_S\sum^{\nu(k)-1}_{j=0}\dfrac{1}{h^k_j}\chi_{[t^k_j,t^k_{j+1})}\gg^k(t)dt\\[3ex]
\disp\int_Sd\al^{1k}_{mes}:=\int_S\al^{1k}(t)dt,\;\int_Sd\al^{2k}_{mes}:=\int_S\al^{2k}(t)dt
\end{array}\right.
\end{equation}
for any Borel subset $S\subset[0,T]$, where $\chi_Q$ stands the characteristic function of the set $Q$ that is equal to $1$ on $Q$ and $0$ otherwise. Finally, the nontriviality condition \eqref{ntc1} can be equivalently rewritten as
\begin{equation}\label{dis-non}
\lm^{k}+\|p^{k}(T)\|+\sum_{j=0}^{\nu(k)-1}\big(\|\al^{1k}_j|+\|\al^{2k}_j\|\big)+\|\al^{1k}_{\nu(k)}+\al^{2k}_{\nu(k)}\|+\sum_{j=0}^{\nu(k)-1}\sum^m_{i=1}\|\gg^k_{ij}\|=1,
\end{equation}
and hence all the terms in \eqref{dis-non} are uniformly bounded as $k\in\N$.\\[1ex]
{\bf Step~4:} {\em Compactness of extended adjoint functions in discrete approximations}. At this step we verify the compactness of the sequences of extended adjoint functions constructed at Step~3 in appropriate functional spaces. This allows us at the next step to pass to the limit in the discrete optimality conditions and to justify in this way the claimed primal-dual dynamic relationships.

First we get from \eqref{dis-non} that there exists $\lambda\ge 0$ such that $\lm^k\to\lm$ along a subsequence of $k\to\infty$. Our next goal is to justify the uniform boundedness of the sequence $\{p^k_0,\ldots,p^{k}_{\nu(k)}\}$ as $k\in\N$. Indeed, observe by \eqref{conx} that
\begin{equation*}
p_{j+1}^{xk}=p_{j}^{xk}+h^k_j\lm^kw_{j}^{xk}-h^k_j\nabla_xg\(\ox_j^k,\ou_j^k\)^*\Big(-\frac{1}{h^k_j}\lm^k\th_{j}^{xk}-\lm^k v_{j}^{xk}+p_{j+1}^{xk}\Big)+h^k_j\sum_{i=1}^{m}\gg_{ij}^k\oa_{ij}^k
\end{equation*}
for all $j=0,\ldots,\nu(k)-1$, which gives us the estimate
\begin{equation*}
\begin{aligned}
\|p^{xk}_j\|&\le\|p^{xk}_{j+1}\|+h^k_j\lm^k\|w_{j}^{xk}\|+h^k_jL_g\Big\|-\frac{1}{h^k_j}\lm^k\th_{j}^{xk}-\lm^k v_{j}^{xk}+p_{j+1}^{xk}\Big\|+h^k_j\Big\|\sum_{i=1}^{m}\gg_{ij}^k\oa_{ij}^k\Big\|\\
&\le\(1+h^k_jL_g\)\|p^{xk}_{j+1}\|+h^k_j\lm^k\big(\|w_{j}^{xk}\|+L_g\|v^{xk}_j\|\big)+h^k_j\lm^kL_g\|\th^{xk}(t^k_j)\|+h^k_j\Big\|\sum_{i=1}^{m}\gg_{ij}^k\oa_{ij}^k\Big\|
\end{aligned}
\end{equation*}
for all $j=0,\ldots,\nu(k)-1$, where $L_g$ is given in \eqref{e:gL}. Denote further
\begin{equation}\label{mj}
M^k_j:=h^k_j\lm^k\big(\|w_{j}^{xk}\|+L_g\|v^{xk}_j\|\big)+h^k_j\lm^kL_g\|\th^{xk}(t^k_j)\|+h^k_j\Big\|\sum_{i=1}^{m}\gg_{ij}^k\oa_{ij}^k\Big\|.
\end{equation}
It then follows from the $\th^{xk}$-estimate in \eqref{th-est} that
$$
\begin{aligned}
&\sum^{\nu(k)-1}_{j=0}h^k_j\lm^k\|\nabla_x g(\ox_j^k,\ou_j^k)\|\cdot\|\th^{xk}(t^k_j)\|\le\lm^kL_g\sum^{\nu(k)-1}_{j=0}\sqrt{h^k_j\int_{t^k_j}^{t^k_{j+1}}\n\th^{xk}(t)\en^2dt}\\
&\quad \le \lm^kL_g \sqrt{\nu(k)h_k}\sqrt{\int^T_0\n\th^{xk}(t)\en^2dt}\le\lm^kL_g\sqrt{\Tilde\nu}\sqrt{\int^T_0\n\th^{xk}(t)\en^2dt}\to 0\;\;\mbox{as}\;\;k\to\infty.
\end{aligned}
$$
On the other hand, the discrete nontriviality condition \eqref{dis-non} immediately yields
\begin{equation*}
\sum_{j=0}^{\nu(k)-1}h^k_j\Big\|\sum_{i=1}^m\gg^k_{ij}\oa^k_{ij}\Big\|=\int_0^T\sum_{i=1}^m\big\|\gg^k_i(t)\oa^{k}_i(t)\big\|dt=\int_0^T\sum_{i=1}^k\big|\gg^k_i(t)\big|dt\le 1.
\end{equation*}
Furthermore, we get from the subdifferential inclusion \eqref{wv} due to the imposed structure \eqref{run_cost} of the running cost and its Lipschitz continuity with constant $L_\ell$ that
$$
\sum^{\nu(k)-1}_{j=0}\|h^k_jw^{xk}_j\|=\sum^{\nu(k)-1}_{j=0}h^k_j\|w^{xk}(t_j)\|\le\sum^{\nu(k)-1}_{j=0}L_\ell h_k\le Lh_k\nu(k) \le L_\ell\Tilde\nu=:\Tilde L<\infty.
$$
The same arguments give us also the estimate $\disp\sum^{\nu(k)-1}_{j=0}\|h^k_jv^{xk}_j\|\le\Tilde L$. Combining the above and remembering the definition of $M^k_j$ in \eqref{mj} shows that
$\disp\sum^{\nu(k)-1}_{j=0}M^k_j\le\Tilde M$ with some constant $\Tilde M>0$ and implies in turn
\begin{equation*}
\|p^{k}_j\|\le(1+L_gh^k_j)\|p^k_{j+1}\|+M^k_j\;\mbox{ for all }\;j=0,\ldots,\nu(k)-1.
\end{equation*}
Arguing now by induction and using the standard progression estimate lead us to the inequalities
$$
\begin{aligned}
\|p^{k}_j\|&\le(1+L_gh^k_j)^{k-j}\|p^{k}_{\nu(k)}\|+\sum_{i=j}^{\nu(k)-1}M^k_i(1+L_gh^k_j)^{i-j}\\
&\le e^{L_g}+e^{L_g}\sum_{j=0}^{k-1}M^i_m\le e^{L_g}(1+\Tilde M)=:\Tilde M_1,\quad j=2,\ldots,\nu(k)-1,
\end{aligned}
$$
which justifies the boundedness of $\{p^k_{j}\}_{2\le j\le \nu(k)}$ and hence the boundedness of $\{p^k_{j}\}_{0\le j\le \nu(k)}$. Then we have
$$
\begin{aligned}
\sum^{\nu(k)-1}_{j=0}\|q^{xk}(t^k_{j+1})-q^{xk}(t^k_j)\|&=\sum^{\nu(k)-1}_{j=0}\|p^{xk}_{j+1}-p^{xk}_j\|\le\sum^{\nu(k)-1}_{j=0}M^k_j
+\sum^{\nu(k)-1}_{j=0}h^k_j\|\nabla_xg(\ox^k_j,\ou^k_j)^*\|\cdot\|p^{xk}_{j+1}\|\\
&\le\Tilde M +L_gh^k_j\nu(k)\Tilde M_1\le\Tilde M+\Tilde\nu L_g\Tilde M_1,
\end{aligned}
$$
which ensures that the functions $q^{xk}(\cdot)$ are of uniform bounded variation on $[0,T]$. To verify the same property for the sequences of $q^{ak}(\cdot)$ and $q^{bk}(\cdot)$ on $[0,T]$, observe from \eqref{cona} that
$$
\begin{aligned}
\sum^{\nu(k)-1}_{j=0}\|p^{ak}_{j+1}-p^{ak}_j\|&\le\sum^{\nu(k)-1}_{j=0}\Big[\lm^kh^k_j\|w^{ak}_j\|+\|\al^{1k}_j+\al^{2k}_j\|\cdot\|\oa^k_j\|+h^k_j\|\gg^k_j\|\cdot\|\ox^k_j\|\\
&\;+h^k_j\|\eta^k_j\|(\lm^k\|\th^{xk}(t_j)\|+\lm^k\|v^{xk}_j\|+\|p^{xk}_{j+1}\|)\Big]\\
&\le\lm^k\Tilde L+(1+\ve_k)\sum^{\nu(k)-1}_{j=0}\Big[\int^{t^k_{j+1}}_{t^k_j}\|\al^{1k}(t)\|dt+\int^{t^k_{j+1}}_{t^k_j}\|\al^{2k}(t)\|dt\Big]\\
&\;+\sum^{\nu(k)-1}_{j=0}\Big(\int^{t^k_{j+1}}_{t^k_j}\|\gg^k(t)\|dt\Big)(h^k_j\|\ox^k_j\|)\\
&\;+\sum^{\nu(k)-1}_{j=0}h^k_j\|\eta^k_j\|(\lm^k\|\th^{xk}(t_j)\|+\lm^k\|v^{xk}_j\|+\|p^{xk}_{j+1}\|)\Big].
\end{aligned}
$$
The latter implies that there exists a constant $\Tilde M_2$ ensuring the estimates
$$
\|\ox^k_j\|\le\|x_0\|+\int^{t^k_j}_0\|\dot\ox^k(t)\|dt\le\|x_0\|+\int^{T}_0\|\dot\ox^k(t)\|dt\le\|x_0\|+\sqrt T\sqrt{\int^{T}_0\|\dot\ox^k(t)\|^2dt}\le\Tilde M_2\;\mbox{ and }
$$
$$
\|\th^{xk}(t_j)\|\le\int_{t_j^k}^{t_{j+1}^k}\Big\|\dfrac{\ox_{j+1}^k-\ox_j^k}{h^k_j}-\dot{\ox}(t)\Big\|dt\le
\sqrt{h^k_j}\sqrt{\int_{t_j^k}^{t_{j+1}^k}\Big\|\dfrac{\ox_{j+1}^k-\ox_j^k}{h^k_j}-\dot{\ox}(t)\Big\|^2dt}\le\dfrac{\ve}{2}\sqrt{h^k_j}
$$
for all $j=0,\ldots,\nu(k)-1$. Similarly to the above arguments involving $w^{xk}_j$, we get that
$$
\|v^{xk}_j\|\le L_\ell\;\;\mbox{and}\;\;\|p^{xk}_{j+1}\|\le\Tilde M_1\;\mbox{ as }\;j=0,\ldots,\nu(k)-1
$$
and then find, by using the boundedness of $\{\eta^k_i(\cdot)\}$ in $L^2$, a constant $\Tilde M_3>0$ such that
$$
\begin{aligned}
\sum^{\nu(k)-1}_{j=0}\|q^{ak}(t^k_{j+1})-q^{ak}(t^k_j)\|&=\sum^{\nu(k)-1}_{j=0}\|p^{ak}_{j+1}-p^{ak}_j\|\le\lm^k\Tilde L+(1+\ve_k)\[\int^{T}_{0}\|\al^{1k}(t)\|dt+\int^{T}_{0}\|\al^{2k}(t)\|dt\]\\
&\;+h_k\Tilde M_2\int^T_0\|\gg^k(t)\|dt +\(\lm^k\dfrac{\ve}{2}\sqrt{h_k}+\lm^kL+\Tilde M_1\)\sum^{\nu(k)-1}_{j=0}\int^{t^k_{j+1}}_{t^k_j}\|\eta^k(t)\|dt\\
&\le\lm^k\Tilde L+1+\ve_k+h_k\Tilde M_2+\(\lm^k\dfrac{\ve}{2}\sqrt{h_k}+\lm^kL+\Tilde M_1\)\sqrt T\sqrt{\int^T_0\|\eta^k(t)\|^2dt}\\
&\le\Tilde L+2+\Tilde M_2+(\ve+L+\Tilde M_1)\sqrt{Tm\Tilde M_3},
\end{aligned}
$$
which justifies that $q^{ak}(\cdot)$ are of uniform bounded variation on $[0,T]$. For $q^{bk}(\cdot)$ we deduce from \eqref{conb} that
$$
\begin{aligned}
\sum^{\nu(k)-1}_{j=0}\|q^{bk}(t_{j+1})-q^{bk}(t_j)\|&=\sum^{\nu(k)-1}_{j=0}\|p^{bk}_{j+1}-p^{bk}_j\|\le\sum^{\nu(k)-1}_{j=0}h^k_j\lm^k\|w^{bk}_j\|+\sum^{\nu(k)-1}_{j=0}h^k_j\|\gg^k_j\|\\
&\le\lm^k\Tilde L+\int^T_0\|\gg^k(t)\|dt\le\Tilde L+1
\end{aligned}
$$
ensuring the uniform bounded variation of $q^{bk}(\cdot)$ and hence of the entire sequence $\nn q^k(\cdot)\hnn$ on $[0,T]$. Since
$$
\|q^k(t)\|-\|q^k(T)\|\le\|q^k(t)-q^k(0)\|+\|q^k(T)-q^k(t)\|\le{\rm var}(q^k;[0,T])<\infty,
$$
where ${\rm var}(q^k;[0,T])$ denotes the total variation of the function $q^k(\cdot)$ on $[0,T]$, the latter gives us the estimates
$$
\max_{t\in[0,T]}\|q^k(t)\|\le\|q^k(T)\|+{\rm var}(q^k;[0,T])=\|p^k(T)\|+{\rm var}(q^k;[0,T])\le 1+{\rm var}(q^k;[0,T])
$$
and therefore verifies the boundedness of $\{q^k(\cdot)\}$ on $[0,T]$. Employing now Helly's selection theorem, we find a function of bounded variation $q(\cdot)$ such that $q^k(t)\to q(t)$ as $k\to\infty$ for all $t\in[0,T]$. It then follows from \eqref{dis-non} that the sequences of measures $\{\al^{1k}_{mes}\}$, $\{\al^{2k}_{mes}\}$, and $\{\gg^k_{mes}\}$ are bounded in the spaces ${\cal C}^*([0,T];\R_+)$, ${\cal C}^*([0,T];\R_-)$, and ${\cal C}^*([0,T];\R^m)$, respectively. This gives us measures $\gg\in{\cal C}^*([0,T];\R^m)$, $\al^1\in{\cal C}^*([0,T];\R_+)$, and $\al^2\in{\cal C}^*([0,T];\R_-)$ for which the triples $(\gg^k_{mes},\al^{1k}_{mes},\al^{2k}_{mes})$ weak$^*$ converge to $(\gg,\al^1,\al^2)$ along some subsequence. Finally, it follows from \eqref{pqg}, \eqref{dis-non}, and the uniform boundedness of $q^k(\cdot),\;w^k(\cdot)$, and $v^k(\cdot)$ on $[0,T]$ that the sequence $\{p^k(\cdot)\}$ is bounded in $W^{1,2}([0,T];\R^{n+mn+m})$ and hence it is weakly compact in this space.\\[1ex]
{\bf Step~5:} {\em Passing to the limit in the adjoint dynamic relationships}. At this step we employ the compactness results established at the previous step of the proof to furnish the limiting procedures in deriving the remaining primal-dual relationships listed in item (i) of the theorem; namely, the adjoint dynamic systems, the normal cone adjoint inclusions, the local and global maximization conditions, and the second dynamic complementarity slackness condition that are formulated above.

First we use the weak compactness of the sequence $\{p^k(\cdot)\}$ in the space $W^{1,2}([0,T];\R^{2n+mn})$ established above and, invoking once more Mazur's weak closure theorem, find a function $p(\cdot)\in W^{1,2}([0,T];\R^{n+mn+m})$ such that a sequence of convex combinations of $\dot p^k(t)$ converges to $\dot p(t)$ for a.e.\ $t\in[0,T]$. Then the passage to the limit in \eqref{e:px-der} and \eqref{e:pa-der} justifies the claimed representation of $\dot p(\cdot)$ in \eqref{c:6.6}.

Our next goal is to verify the claimed representation \eqref{c:6.9} of the BV-adjoint arc $q(\cdot)$. We proceed similarly to \cite[p.\ 325]{v} to get the convergence
$$
\n\int_{[t,T]}\gg^k_i(s)\oa^k_i(\vth^k(s))ds-\int_{[t,T]}\oa_i(s)d\gg_i(s)\en\to 0\;\;\mbox{as}\;\;k\to\infty
$$
for all $t\in[0,T]$ except a countable subset of $[0,T]$, which tells us in turn that
$$
\int_{[t,T]}\gg^k_i(s)\oa^k_i\big(\vth^k(s)\big)ds\to\int_{[t,T]}\oa_i(s)d\gg_i(s)\;\;\mbox{as}\;\;k\to\infty
$$
for all $i=1,\ldots,m$. Then we have estimate
\begin{equation*}
\begin{aligned}
&\bigg\|\int_{[t,T]}\big[\al^{1k}(s)+\al^{2k}(s),\oa^k_i\big(\vth^k(s)\big)\big]ds-\int_{[t,T]}\oa_i(s)d\al(s)\bigg\|\\
\le&\bigg\|\int_{[t,T]}\big[\al^{1k}(s)+\al^{2k}(s),\oa^k_i\big(\vth^k(s)\big)\big]ds-\int_{[t,T]}\big[\al^{1k}(s)+\al^{2k}(s),\oa_i(s)\big]ds\bigg\|\\
+&\bigg\|\int_{[t,T]}\big[\al^{1k}(s)+\al^{2k}(s),\oa_i(s)\big]ds-\int_{[t,T]}\oa_i(s)d\al(s)\bigg\|\\
=&\bigg\|\int_{[t,T]}\[\al^{1k}(s)+\al^{2k}(s),\oa^k_i\big(\vth^k(s)\big)-\oa_i(s)\]ds\bigg\|\\
+&\bigg\|\int_{[t,T]}\big[\al^{1k}(s)+\al^{2k}(s),\oa_i(s)\big]ds-\int_{[t,T]}\oa_i(s)d\al(s)\bigg\|,
\end{aligned}
\end{equation*}
where $\al:=\al^1+\al^2$. Combining \eqref{dis-non} with the fundamental H\"older integral inequality gives us
$$
\bigg\|\int_{[t,T]}\[\al^{1k}(s)+\al^{2k}(s),\oa^k_i\big(\vth^k(s)\big)-\oa_i(s)\]ds\bigg\|\le\big\|\al^{1k}(\cdot)+\al^{2k}(\cdot)\big\|_{L^2}\|\oa^k_i(\cdot)-\oa_i(\cdot)\|_{L^2}\le \|\oa^k_i(\cdot)-\oa_i(\cdot)\|_{L^2},
$$
which justifies the convergence to zero of the first term in the above estimate. The second term therein also converges to zero for all $t\in[0,T]$ except some countable subset due to the weak$^*$ convergence of $\al^{1k}\to\al^1$ in ${\cal C}^*([0,T];\R_+)$ and $\al^{2k}\to\al^2$ in ${\cal C}^*([0,T];\R_-)$. Thus we get \eqref{c:6.9} by passing to the limit the differential equation \eqref{pk-qk} and using the justified convergence
$$
\int_{[t,T]}\big[\al^{1k}(s)+\al^{2k}(s),\oa^k_i\big(\vth^k(s)\big)\big]ds\to\int_{[t,T]}\oa(s)d\al(s)\;\mbox{ as }\;k\to\infty.
$$

Let us now show that the triple $\psi(\cdot)=(\psi^u(\cdot),\psi^a(\cdot),\psi^b(\cdot))$ defined in \eqref{c:6.6'} and \eqref{c:6.6''} on $[0,T]$ satisfies the normal cone inclusions claimed in those conditions. Indeed, it follows from the construction of $\psi^k(t)$ in \eqref{c:6.25}, from the necessary optimality conditions in \eqref{cony} and \eqref{dloc-max} for the discrete problems $(P^k)$, and from the convergence of all the extended functions  defining $\psi^k(\cdot)$ in \eqref{cony} and \eqref{dloc-max}, which was established in the proof above, that a subsequence $\{\psi^k(\cdot)\}$ weakly converges in $L^2([0,T];\R^{d+mn+m})$. This clearly implies, by using again Mazur's weak closure theorem and passing to the limit in \eqref{c:6.25} for the convexified sequences in both sides therein, that the limiting function $\psi(t)$ satisfies the equations in \eqref{c:6.6'} and \eqref{c:6.6''}. Furthermore, we have the inclusions
\begin{equation}\label{psiAk}
\begin{array}{ll}
&\psi^{uk}(t)\in N\big(u^k(\vth^k(t));U\big),\;\psi^{Ak}(t)\in N\big(\dot a^k(\vth^k(t));A\big),\\
&\psi^{Bk}(t)\in N\big(\dot b^k(\vth^k(t));B\big)\;\mbox{ for a.e. }\;t\in[0,T].
\end{array}
\end{equation}
Passing to the limit in \eqref{psiAk} as $k\to\infty$ with the usage of Mazur's theorem and the robustness of the normal cone \eqref{nor_con} tells us that the limiting function $\psi(t)$ satisfies the convexified inclusions in \eqref{c:6.6'} and \eqref{c:6.6''} for a.e.\ $t\in[0,T]$. The local and global maximization conditions are derived, under the imposed additional assumptions, from the normal cone inclusions in \eqref{c:6.6'} and \eqref{c:6.6''} similarly to the case of discrete-time systems in Corollary~\ref{cor}.

The final segment of this step is verifying the second complementary slackness condition in \eqref{41} under the imposed LICQ assumption. Fix $i\in\{1,\ldots,m\}$ and $t$ from the set of full measure on $[0,T)$ such that $\eta_i(t)>0$. As proved at Step~2, we have $\eta^k_i(\cdot)\to\eta_i(\cdot)$ a.e.\ on $[0,T)$, and hence $\eta^k_i(t)>0$ whenever $k\in\N$ is sufficiently large. This allows us to use the complementarity slackness condition \eqref{96} in $(P^k)$ and conclude that
\begin{equation}\label{compl2}
\big\la\oa^k_i(t),q^{xk}\big(\vth^k(t)+h^k_j\big)-\lm^k\big(\th^{xk}(t)+v^{xk}(t)\big)\big\ra=0
\end{equation}
for all large $k$. Passing to the limit in \eqref{compl2} as $k\to\infty$ and remembering the corresponding convergence of the sequences therein established above give us
\begin{equation*}
\la\oa_i(t),q^x(t)-\lm v^x(t)\ra=0\;\mbox{ for a.e. }\;t\in[0,T),
\end{equation*}
which readily justifies the second complementary slackness condition in \eqref{41}.\\[1ex]
{\bf Step~6:} {\em Verifying the endpoint optimality conditions.} Here we furnish the passage to the limit from the discrete approximations to derive the transversality and complementary slackness conditions at the right endpoint of the controlled sweeping process $(P)$. The developed finite-dimensional limiting procedure is significantly less involved in comparison with its dynamic counterpart accomplished above and allows us, in particular, to verify the endpoint complementary slackness conditions without the additional LICQ assumption.

Define $\eta^k(T):=\eta^k_{\nu(k)}$ for all $k\in\N$, where the latter vectors are taken from Theorem~\ref{Th1}. It follows from the normalized nontriviality condition \eqref{dis-non} for discrete problems $(P^k)$ that the sequence $\{\eta^k(T)\}$ contains a subsequence that converges to some $\eta(T)=(\eta_1,\ldots,\eta_m)\in\R^m_+$ as $k\to\infty$. It follows from \eqref{nmuta} and \eqref{t:7.20} that
\begin{equation}\label{trans1}
\begin{aligned}
&p^{ak}_{\nu(k)}+\big[\eta^k_{\nu(k)},\rep_m(\ox^k_{\nu(k)})\big]=-2\big[\al^{1k}_{\nu(k)}+\al^{2k}_\nu(k),\oa^k_{\nu(k)}\big]\\
&\quad \in-2\oa^k_{\nu(k)}\big(N_{[0,1+\ve_k]}(\|\oa^k_{\nu(k)}\|)+N_{[1-\ve_k,\infty)}(\|\oa^k_{\nu(k)}\|)\big).
\end{aligned}
\end{equation}
Passing now to the limit in \eqref{nmuta} and \eqref{trans1} as $k\to\infty$ and $\ve_k\dn 0$ with taking into account the convergence of $\{\al^{1k}_{\nu(k)}+\al^{2k}_{\nu(k)}\}$ due to \eqref{dis-non}, we verify the transversality conditions in \eqref{42a}. The passage to the limit in \eqref{nmutx} and \eqref{nmutb} by using the justified convergence of $\{p^k_{\nu(k)}\}$ to $p(T)$ gives us \eqref{42b}. To verify finally the remaining transversality condition \eqref{42x} accompanied by the inclusion in \eqref{42c}, we pass to the limit in the discrete counterpart \eqref{nmutx} as $k\to\infty$ with the usage of representation \eqref{F-rep} for the normal cone \eqref{NC} to the polyhedral convex set \eqref{SP} as well as the robustness of the subdifferential mapping \eqref{sub}.\\[1ex]
{\bf Step~7:} {\em Verifying the measure nonatomicity conditions.} Starting with condition (a), pick any $t\in[0,T]$ with $\la\oa_i(t),\ox(t)\ra<\ob_i(t)$ for all $i=1,\ldots,m$ and by the continuity of $\ox(\cdot)$ find a neighborhood $V_t$ of $t$ such that $\la\oa_i(s),\ox(s)\ra<\ob_i(s)$ for these indices whenever $s\in V_t$. The obtained convergence of the discrete optimal solutions tells us that $\la\oa^k_i(t^k_j),\ox^k(t^k_j)\ra<\ob^k_i(t^k_j)$ if $t^k_j\in V_t$ for all $j=0,\ldots,\nu(k)-1$, $i=1,\ldots,m$, and large $k\in\N$. It follows from \eqref{eta} and \eqref{94} that $\eta^k_{ij}=0$ and $\gg^k_{ij}=0$ for all such indices. Thus we deduce from \eqref{e:discrete-measures} that
\begin{equation*}
\|\gg^k\|(V)=\disp\int_Vd\|\gg^k\|=\int_V\|\gg^k(t)\|dt=0.
\end{equation*}
Passing now to the limit as $k\to\infty$ and using the measure convergence established in Step~3 give us $\|\gg\|(V)=0$, which verifies the claimed condition (a). The measure nonatomicity condition (b) is justified similarly.\\[1ex]
{\bf Step~8:} {\em Verifying the measured coderivative condition.} To prove \eqref{mcc}, we first rewrite \eqref{congg} in the form
\begin{equation*}
\gg^k_i(t)\in D^*N_{\R_-}\Big(\big\la\oa^k_i(t),\ox^k(t)\big\ra-\ob^k_i(t),\eta^k_i(t)\Big)\Big(\big\la\oa^k_i(t),-\lm^k\th^{xk}(t)-\lm^kv^{xk}(t)+q^{xk}\big(\vth^k(t)+h^k_j\big)\big\ra\Big)
\end{equation*}
for $t\in(t^k_j,t^k_{j+1})$ as $j=0,\ldots,\nu(k)-1$ and $i=1,\ldots,m$. Pick $t\in[0,T]$ and for any $k\in\N$ find $j\in\{0,\ldots,\nu(k)-1\}$ such that $t\in[t^k_j,t^k_{j+1}]=:S_k$. It follows from the construction of $\gg^k_{mes}$ that
\begin{equation*}
\gg^k(t)=\dfrac{1}{h^k_j}\int_{S_k}\gg^k(s)ds=\int_{S_k}d\gg^k_{mes}=\gg^k_{mes}(S_k).
\end{equation*}
Furthermore, the sequence $\{\gg^k(t)\}$ is bounded by \eqref{dis-non}, which allows us to select a subsequence of $\{\gg^k(t)\}$ (without relabeling) and a vector $y\in\R^n$ such that $\gg^k(t)\to y$ as $k\to\infty$. Combining the latter with the coderivative robustness, we arrive at the inclusion
\begin{equation*}
y\in D^*N_{\R_-}\Big(\big\la\oa_i(t),\ox(t)\big\ra-\ob_i(t),\eta(t)\Big)\Big(\big\la\oa_i(t),q^x(t)-\lm v^x(t)\big\ra\Big)\cap\Limsup_{|S|\to 0,\;\Tilde\gg_i\to\gg_i}\Tilde\gg_i(S)(t)
\end{equation*}
for $i=1,\ldots,m$, which thus verifies the claimed measured coderivative condition \eqref{mcc}.\\[1ex]
{\bf Step~9:} {\em Verifying the nontriviality conditions.} Let us first show that
\begin{equation}\label{e:83}
\lm+\|p(T)\|+\|\al\|_{TV}+\|\gg\|_{TV} \ne 0,
\end{equation}
where $\|\cdot\|_{TV}$ denotes the measure total variation. Arguing by contraposition, suppose that \eqref{e:83} fails, i.e., $\lm=0$, $p(T)=0$. This yields $\lm^k\to 0$ and $p^k(T)\to 0$, which implies in turn that $p^{ak}_{\nu(k)}\to 0$ and $p^{bk}_{\nu(k)}\to 0$ as $k\to\infty$. Next we justify the following limiting relationships while considering the two cases:
\begin{equation}\label{al-mes}
\disp\int_{[0,T]}|\al^{1k}(t)|dt\to 0\;\mbox{ and }\;\disp\int_{[0,T]}|\al^{2k}(t)|dt\to 0\;\mbox{ as }\;k\to\infty.
\end{equation}\\
$\bullet$ $0\notin Q$. In this case we clearly have the first convergence condition in \eqref{al-mes}.\\[1ex]
$\bullet$ $0\in Q$, The countability of $Q$ allows us to find in this case a real number $\delta\in\big(0,(h^k_0)^2\big)$ with $\delta\notin Q$. Then we have
$\disp\lim_{k\to\infty}\int_{[\delta,T]}|\al^{1k}(t)|dt=0$, which implies therefore the estimates
\begin{equation*}
\begin{aligned}
\disp0&\le\int_0^T|\al^{1k}(t)|dt=\int_\delta^T|\al^{1k}(t)|dt+\int_0^\delta|\al^{1k}(t)|dt\le\int_\delta^T|\al^{1k}(t)|dt+\int_0^{(h_{0}^k)^2}\dfrac{\al^{1k}_0}{h^k_0}dt\\
&\le\int_\delta^T|\al^{1k}(t)|dt+\al^{1k}_0h^k_0\to 0\;\;\mbox{as}\;\;k\to\infty
\end{aligned}
\end{equation*}
since $0\le\al^{1k}_0\le 1$ due to \eqref{dis-non}. This readily verifies the first convergence condition in \eqref{al-mes}. The second conditions therein is justified similarly. Hence
\begin{equation*}
\disp\sum_{j=0}^{\nu(k)-1}\big(|\al^{1k}_j|+|\al^{2k}_j|\big)\to 0\;\mbox{ as }\;k\to\infty.
\end{equation*}

Next we aim at checking the endpoint convergence
\begin{equation}\label{al-end}
\al^{1k}_{\nu(k)}+\al^{2k}_{\nu(k)}\to 0\;\mbox{ as }\;k\to\infty.
\end{equation}
Indeed, it follows from \eqref{nmutb} and the above convergence of $p^k(T)\to 0$ that $\eta^k_{\nu(k)}\to 0$ as $k\to\infty$. Then \eqref{nmuta} yields $(\al^{1k}_{\nu(k)}+\al^{2k}_{\nu(k)})\oa^k_{\nu(k)}\to 0$ as $k\to\infty$, which gives us \eqref{al-end} since $\oa^k_{\nu(k)}\not=0$ for all large $k$ by \eqref{const_a}.

Now let us proceed with the $\gg$-measure part. By \eqref{dis-non} it follows from the above that
\begin{equation*}
\disp\sum_{j=0}^{\nu(k)-1}\sum^m_{i=1}|\gg^k_{ij}|\to 1\;\mbox{ as }\;k\to\infty.
\end{equation*}
Define the sequence of measurable vector functions $\be^k\colon[0,T]\to\R^m$ by
\begin{equation}\label{bek}
\be^k_i(t):=\left\{\begin{array}{cccc}
\dfrac{\gg^k_i(t)}{|\gg^k_i(t)|}&\mbox{if }\;\gg^k_i(t)\not=0,\\
0&\mbox{if }\;\gg^k_i(t)=0
\end{array}\right.
\end{equation}
for $i=1,\ldots,m$ and all $t\in[0,T]$. The Jordan measure decomposition gives us the representations $\gg^k_{mes}=(\gg^k)^+-(\gg^k)^-$ and $\gg=\gg^+-\gg^-$. Since the space ${\cal C}^*([0,T];\R^m)$ is separable and the measure sequence $\{\gg^k_{mes}\}$ is bounded, we find a subsequence of $\{\gg^k\}$ (without relabeling) such that $\{(\gg^k)^+\}$ and $\{(\gg^k)^-\}$ weak$^*$ converge in ${\cal C}^*([0,T];\R^m)$ to $\gg^+$ and $\gg^-$, respectively. It follows from \eqref{bek} that the sequence $\{\be^k_i(\cdot)\}$ is bounded on $[0,T]$. Applying now the convergence result from  \cite[Proposition~9.2.1]{v} (with $A=[0,1]^m$ therein) gives us Borel measurable functions $\be^+,\be^-\colon[0,T]\to\R^m$ such that some subsequence of $\{\be^k(\gg^k)^+\}$ and $\{\be^k(\gg^k)^-\}$ weak$^*$ converge to $\be^+\gg^+$ and $\be^-\gg^-$, respectively. This yields the relationships
\begin{equation*}
\begin{aligned}
&\disp\bigg\|\int_{[0,T]\backslash Q}\be^+(t)d\gg^+(t)-\int_{[0,T]\backslash Q}\be^-(t)d\gg^-(t)\bigg\|=\disp\lim_{k\to\infty}\bigg\|\int_{[0,T]\backslash Q}\be^k(t)d(\gg^k)^+(t)-\int_{[0,T]\backslash Q}\be^k(t)d(\gg^k)^-(t)\bigg\|\\
&=\disp\lim_{k\to\infty}\bigg\|\int_{[0,T]\backslash Q}\be^k(t)d\gg^k_{mes}(t)\bigg\|=\disp\lim_{k\to\infty}\bigg\|\int_{[0,T]\backslash Q}\(\be^k_1(t)d\gg^k_{1,mes}(t),\ldots,\be^k_m(t)d\gg^k_{m,mes}(t) \)\bigg\|\\
&=\disp\lim_{k\to\infty}\bigg\|\bigg(\sum^{\nu(k)-1}_{j=0}|\gg^k_{1j}|,\ldots,\sum^{\nu(k)-1}_{j=0}\|\gg^k_{mj}\|\bigg)\bigg\|=\disp\lim_{k\to\infty}\sqrt{\sum^m_{i=1}\bigg[\sum^{\nu(k)-1}_{j=0}\|
\gg^k_{ij}|\bigg]^2}\\
&\ge\disp\lim_{k\to\infty}\dfrac{1}{\sqrt m}\disp\sum^m_{i=1}\sum^{\nu(k)-1}_{j=0}|\gg^k_{ij}|=\disp\dfrac{1}{\sqrt m}>0
\end{aligned}
\end{equation*}
for some countable set $Q\subset[0,T]$. Furthermore, we have the estimates
\begin{equation*}
\begin{aligned}
&\bigg\|\int_{[0,T]\backslash Q}\be^+(t)d\gg^+(t)-\int_{[0,T]\backslash Q}\be^-(t)d\gg^-(t)\bigg\|\le\bigg\|\int_{[0,T]\backslash Q}\be^+(t)d\gg^+(t)\bigg\|+\bigg\|\int_{[0,T]\backslash Q}\be^-(t)d\gg^-(t)\bigg\|\\
&=\bigg\|\bigg(\int_{[0,T]\backslash Q}\be^+_1(t)d\gg^+_1(t),\ldots,\int_{[0,T]\backslash Q}\be^+_m(t)d\gg^+_m(t)\bigg)\bigg\|+\bigg\|\bigg(\int_{[0,T]\backslash Q}\be^-_1(t)d\gg^-1(t),\ldots,\int_{[0,T]\backslash Q}\be^-_m(t)d\gg^-_m(t)\bigg)\bigg\|\\
&=\sqrt{\bigg[\int_{[0,T]\backslash Q}\be^+_1(t)d\gg^+_1(t)\bigg]^2+\ldots+\bigg[\int_{[0,T]\backslash Q}\be^+_m(t)d\gg^+_m(t)\bigg]^2}\\
&+\sqrt{\bigg[\int_{[0,T]\backslash Q}\be^-_1(t)d\gg^-_1(t)\bigg]^2+\ldots+\bigg[\int_{[0,T]\backslash Q}\be^-_m(t)d\gg^-_m(t)\bigg]^2}\\
&\le\sqrt{\bigg[\int_{[0,T]\backslash Q}d\gg^+_1(t)\bigg]^2+\ldots+\bigg[\int_{[0,T]\backslash Q}d\gg^+_m(t)\bigg]^2}+\sqrt{\bigg[\int_{[0,T]\backslash Q}d\gg^-_1(t)\bigg]^2+\ldots+\bigg[\int_{[0,T]\backslash Q}d\gg^-_m(t)\bigg]^2}\\
&=\big\|\big(\gg^+_1([0,T]\backslash Q),\ldots,\gg^+_m([0,T]\backslash Q)\big)\big\|+\big\|\big(\gg^-_1([0,T]\backslash Q),\ldots,\gg^-_m([0,T]\backslash Q)\big)\big\|\\
&=\int_{[0,T]\backslash Q}d\|\gg^+(t)\|+\int_{[0,T]\backslash Q}d\|\gg^-(t)\|\le\|\gg^+\|_{TV}+\|\gg^-\|_{TV}=\|\gg\|_{TV}.
\end{aligned}
\end{equation*}
Combining the above relationships shows that $\|\gg\|_{TV}>0$. Thus we arrive at a contradiction and verify therefore the claimed condition \eqref{e:83}.

Next we check that \eqref{e:83} yields the measure-free general nontriviality condition in \eqref{gen-nontr}. Arguing by contraposition, suppose that $(\lm,p,q)=0$ and get from the third component in the representation of $q(t)$ in \eqref{c:6.9} that $\gg(\cdot)=0$. Using then the second component of \eqref{c:6.9} tells us that $[\oa(s),d\al(s)]=0$, which implies by \eqref{a} that $\|\al\|_{TV}=0$, and thus we come to a contradiction with \eqref{e:83}.

To establish now the enhanced nontriviality condition \eqref{enh0} under the imposed additional assumptions, suppose on the contrary that $\lm=0$, $p(T)=0$, $\psi^A(t)=0$, and $\psi^B(t)=0$ for all $t\in[0,T]$. It immediately follows from \eqref{c:6.6''} that $q^a(t)=0$ and $q^b(t)=0$ for all $t\in[0,T]$. Then we deduce from the third component of the equality in \eqref{c:6.6} that $p^b(t)=0$ for all $t\in[0,T]$. Furthermore, the first complementarity slackness condition in \eqref{41}, which is valid under the imposed LICQ, tells us that $\eta(t)=0$ for a.e.\ $t\in[0,T)$ due to the interiority condition \eqref{inter0}. Thus we get from the second component of \eqref{c:6.6} that $p^a(t)=0$ for all $t\in[0,T]$. Looking again at \eqref{c:6.9} and taking into account that $q(\cdot)$ is right continuous lead us to a contradiction with the general nontriviality condition \eqref{gen-nontr}. The obtained contradiction verifies the enhanced nontriviality condition in \eqref{enh0}.

It remains to show that the enhanced nontriviality condition \eqref{enh0} is reduced to the simplified form \eqref{enh1} under the fulfillment of the endpoint interiority assumptions \eqref{inter} and \eqref{psi-con}. Indeed, the former assumption ensures by the implication in \eqref{42b} that the violation of \eqref{enh1} yields $p^b(T)=0$, and hence $\eta(T)=0$ by the equality in \eqref{42b}. In addition, it follows from \eqref{psi-con} and \eqref{c:6.6''} that $\psi^A(t)=0$ and $\psi^B(t)=0$ for all $t\in[0,T]$. Finally, we deduce from \eqref{42x} that $p^x(T)=0$, and thus $p(T)=0$, which contradicts the enhanced nontriviality condition \eqref{enh0}. This completes therefore the entire proof of the theorem. $\h$\vspace*{-0.15in}

\section{Examples}\label{sec:exa}\vspace*{-0.1in}

In this section we present explicit examples that illustrate our main optimality conditions obtained for the sweeping control problem $(P)$. Three different cases will be considered in which the {\em three types of controls} appear separately. Namely, controls in the {\em dynamics}, controls in {\em generating normals} of moving polyhedra, and control by {\em shifting} polyhedron edges.
\vspace*{-0.1in}

\begin{example}[controls in the dynamics]\label{ex_1} Consider first problem $(P)$ with the same data as in our preceding paper \cite{ccmn19a}, where the discrete version of this problem was treated.
Here the $a$-components and $b$-components of controls are fixed, i.e., the set $C(t)$ is given, and only the $u$-components are used for optimization.
Keeping the notation of formulas from \eqref{minimize} to \eqref{ab}, our data read as follows:
\begin{equation}\label{data}
\left\{
\begin{array}{ll}
n=2,\;m=1,\;T=1,\;x_0=\(\frac{3}{2},1\),\;a=\(-\frac{1}{\sqrt5},-\frac{2}{\sqrt 5
}\),\;b=-\frac{2}{\sqrt5},\\[1ex]
g(x,u)=g(x_1,x_2,u_1,u_2):=(u_1,u_2),\;\;\vph(x):=x_1+x_2,\;\;\ell(t,x,a,b,u,\dot x,\dot a,\dot b):=\frac{1}{2}u^2_1+u^2_2,\\[1ex]
U:=[-1,1]\times[-1,1].
\end{array}\right.
\end{equation}
The set $C(t)$ in the sweeping inclusion \eqref{SPC} is independent of time and is described by
$$
C(t)=C:=\nn(x_1,x_2)\in\R^2\big|\;x_1+2x_2\ge 2\hnn\;\mbox{ for all }\;t\in[0,1].
$$
In other words, we wish to minimize the cost functional
$$
J(x,u)=x_1(1)+x_2(1)+\int_0^1\Big(\frac12 u_1^2(t)+u_2^2(t)\Big)dt
$$
subject to the controlled sweeping dynamics
$$
\dot{x}\in-N(x;C)+u,\; x(0)=\Big(\frac32,1\Big)
$$
with $|u_1|,|u_2|\le 1$. In what follows we are going to show that applying the necessary optimality conditions of Theorem~\ref{Th2}, together with some heuristics, allows us to find optimal solutions to the original control problem $(P)$. Since the initial condition lies in the interior of $C$, the first time of hitting the boundary $t^\ast$ is positive. Moreover, in the case where $t^\ast<1$, the object is expected to slide along the boundary of $C$ for the whole interval $[t^\ast,1]$. It is easy to see that all the assumptions of Theorem~\ref{Th2} are satisfied for \eqref{data}, and hence we can employ the obtained necessary optimality conditions to find some dual elements satisfying the following relationships, where we $\bar{x}=(\bar{x}_1,\bar{x}_2)$ and $\bar u=(\bar{u}_1,\bar{u}_2)$ are the optimal trajectory and the optimal control, respectively:
\begin{enumerate}
\item $(w^x(t),w^u(t))=\big(0,\ou_1(t),2\ou_2(t)\big)$ a.e.\ $t\in[0,1]$.
\item $v^x(t)=0$ a.e.\ $t\in[0,1]$.
\item $\dot\ox(t)=\ou(t)+\dfrac{\eta(t)}{\sqrt5}(1,2)$ a.e.\ $t\in[0,1]$.
\item $p^x(t)\equiv p^x(1)$.
\item $q^x(t)=p^x(1)+\disp\int_{[t,1]}\l(\dfrac{1}{\sqrt5},\dfrac{2}{\sqrt5}\r)d\gg(s)=p^x(1)+\l(\dfrac{1}{\sqrt5},\dfrac{2}{\sqrt5}\r)\gg([t,1])$.
\item $\psi^u(t)=-\lm\big(\ou_1(t),2\ou_2(t)\big)+q^x(t)$, a.e.\ $t\in[0,1]$.
\item $\psi^u(t)\in N(\ou(t);U)$, which implies that
$$
\psi^u_1(t)\ou_1(t)+\psi^u_2(t)\ou_2(t)=\max_{(u_1,u_2)\in[-1,1]\times[-1,1]}\big\{\psi^u_1(t)u_1+\psi^u_2(t)u_2\big\}\;\mbox{ for a.e. }\;t\in[0,1].
$$
\item $\ox_1(t)+2\ox_2(t)>2\sr\gg(V)=0$ whenever $t\in[0,1]$ for all the Borel subsets $V$ of $V_t$ and $\eta(t)=0$, where $V_t$ is a neighborhood of $t$.
\item $\eta(t)>0\sr\big\la(-1,-2),\;q^x(t)-\lm v^x(t)\big\ra=\big\la(-1,-2),\;q^x(t)\big\ra=0$ for all $t\in[0,1]$.
\item $-p^x(1)=\lm(1,1)+\disp\eta(1)\Big(-\frac{1}{\sqrt 5},-\frac{2}{\sqrt5}\Big)$.
\item $\big(\lm,p^x(\cdot),q^x(\cdot)\big)\ne(0,0,0)$ with $\lm\ge 0$.
\end{enumerate}
Combining items~5 and 10 gives us the expression
$$
q^x(t)=-\lm(1,1)+\dfrac{(1,2)}{\sqrt5}\big(\eta(1)+\gg([t,1])\big)\;\mbox{ on }\;[0,1].
$$
Since we start from an interior point of the polyhedron, the measure $\gg$ is zero in the interval $[0,t^\ast)$, where
$$
t^\ast=\inf\big\{t\in[0,1]\;\big|\;x(t)\in{\rm bd}\,C\big\}.
$$
Observe that the value $t^\ast=1$ tells us that either $\ox$ does not hit ${\bf}\,C$, or $\ox$ hits ${\rm bd}\,C$ only at the final time. Let us investigate behavior of the trajectory before and after the hitting time $t=t^\ast$.

Since the trajectory does not hit the boundary before $t=t^\ast$, it follows from item~8 that $\gg([t,1])\equiv\gg([t^\ast,1])$ for all $t\in[0,t^\ast)$. As a consequence, we get
\begin{equation}\label{qx}
q^x(t)\equiv\lm(1,1)+\dfrac{(1,2)}{\sqrt 5}\big(\eta(1)+\gg([t^\ast,1])\big)=:q^\ast=(q^\ast_1,q^\ast_2)
\end{equation}
for all $t\in[0,t^\ast)$, and hence item~6 implies that
$$
\psi^u_1(t)=-\lm\ou_1(t)+q^\ast_1\;\;\mbox{and}\;\;\psi^u_2(t)=-2\lm\ou_2(t)+q^\ast_2\;\mbox{ for a.e. }\;t\in[0,t^\ast).
$$
Thus the maximization condition in item~7 tells us that
$$
-\lm\big(\ou^2_1(t)+2\ou^2_2(t)\big)+q^\ast_1\ou_1(t)+q^\ast_2\ou_2(t)=\max_{|u_1|\le 1,|u_2|\le 1}\big\{-\lm\big(\ou_1(t)u_1+2\ou_2(t)u_2\big)+q^\ast_1u_1+q^\ast_2u_2\big\}
$$
for a.e.\ $t\in[0,t^\ast)$, which can be written by \eqref{qx} as
\begin{equation}\label{max_con}
\begin{split}
-&\lm\big(\ou_1^2(t)+2\ou_2^2(t)\big)-\lm\ou_1(t)-\lm\ou_2(t)+\frac{\eta(1)+\gg({[t^\ast,1]})}{\sqrt{5}}\big(\ou_1(t)+2\ou_2(t)\big)\\
=&\underset{|u_1|,|u_2|\le 1}\max\Big\{\lm\big(u_1(-1-\ou_1(t))+u_2(-1-2\ou_2(t))\big)+\frac{\eta(1)+\gg([t^\ast,1])}{\sqrt{5}}(u_1+2u_2)\Big\}
\end{split}
\end{equation}
with $\eta(1)\ge 0$. In order to maximize the function
$$
\phi(u_1,u_2):=u_1\Big(\frac{\eta(1)+\gg([t^\ast,1])}{\sqrt{5}}-\lm\big(1+\ou_1(t)\big)\Big)+u_2\Big(\frac{2}{\sqrt{5}}\big(\eta(1)+\gg([t^\ast,1])\big)-\lm\big(1+2\ou_2(t)\big)\Big)
$$
with respect to $u_1\in [-1,1]$ and $u_2\in[-1,1]$, we observe that if the optimal value for $u_1$, i.e.,
$\bar{u}_1(t)$ is in the interior of $[-1,1]$, then $\frac{\partial \phi}{\partial u_1}(\bar{u}_1(t))=0$, and the same occurs for
$u_2$. In other words,
\begin{itemize}
\item if $|\ou_1(t)|<1$, then $\eta(1)+\gg([t^\ast,1])=\lm(1+\ou_1(t))\sqrt{5}$, and
\item if $|\ou_2(t)|<1$, then $\eta(1)+\gg([t^\ast,1])=\lm (1+2\ou_2(t))\frac{\sqrt{5}}{2}$,
\end{itemize}
while if both the above cases occur, then
$$
2\lm\big(1+\ou_1(t)\big)=\lm\big(1+2\ou_2(t)\big)=\eta(1)+\gg([t^\ast,1]).
$$
The last case gives us the relationship
$$
\ou_1(t)-\ou_2(t)=-\frac{1}{2}
$$
and tells us that the controls $\ou_1(\cdot)$ and $\ou_2(\cdot)$ are constant on the interval $[0,t^\ast)$ provided that $\lm>0$. (It is shown below that the case where $\lm=0$ cannot occur).
Moreover, in this situation we get that $\eta(1)+\gg([t^\ast,1])>0$, which implies that the optimal trajectory eventually hits the boundary. Indeed, $\eta >0$ or $\gamma [t^\ast,1]>0$ imply (by \eqref{42c} and by the measure nonatomicity conditions) that the constraint is eventually active. The sliding condition then implies that the constraint remains active until the final time. Since the optimal performance on the interval $[0,t^\ast]$ can be obviously obtained with $\ou_1(t),\ou_2(t)\le 0$ a.e., it follows that on $[0,t^\ast]$ the optimal control pair is constant. Thus is natural to assume that the optimal control pair is constant also on $(t^\ast,1]$ together with the normal vector $\eta(t)$, i.e., we conclude that
$$
\ou(t)=\l\{\begin{array}{ll}
\ou_0=(\ou_{10},\ou_{20})&\mbox{ if }\;t\in[0,t^\ast),\\[1ex]
\ou_\ast=(\ou_{1\ast},\ou_{2\ast})&\mbox{ if }\;t\in[t^\ast,1],
\end{array}
\right.
$$
$$
\eta(t)=\begin{cases}
\eta_0=0\;\;\mbox{if}\;\;t<t^\ast,\\
\eta^\ast\ge 0\;\;\mbox{if}\;\;t\ge t^\ast.
\end{cases}
$$
Supposing that $\ox$ slides along the boundary ${\bd}\,C$ on $[t^\ast,1]$ gives us the expression
$$
\ox(t)=\nn
\begin{array}{ll}
\disp\(\frac{3}{2}+t\ou_{10},1+t\ou_{20}\)&\mbox{ if }\;t\in[0,t^\ast),\\[1ex]
\disp\(\frac{3}{2}+t^\ast\ou_{10}+\(t-t^\ast\)\(\ou_{1\ast}+\frac{\eta^\ast}{\sqrt5}\),1+t^\ast\ou_{20}+\(t-t^\ast\)\(\ou_{2\ast}+\frac{2\eta^\ast}{\sqrt5}\)\)&\mbox{ if }\;
t\in[t^\ast,1].
\end{array}
\right.
$$
Consequently, we arrive at the equation for the hitting time
$$
\frac{7}{2}+t^\ast(\ou_{10}+2\ou_{20})+\(t-t^\ast\)\(\ou_{1\ast}+2\ou_{2\ast}+\sqrt 5\eta^\ast\)=2\;\mbox{ if }\;t\ge t^\ast,
$$
which results in the formula
\begin{equation}
\label{htnew}
t^\ast=-\dfrac{3}{2(\ou_{10}+2\ou_{20})}.
\end{equation}
Since $0\le t^\ast\le 1$, this in turn implies the condition
\begin{equation}\label{u-con1}
\ou_{10}+2\ou_{20}\le-\frac{3}{2}.
\end{equation}
Observe further that the sliding condition on $[t^\ast,1]$ implies that $\ou_{1\ast}+2\ou_{2\ast}+\sqrt5\eta^\ast=0$, from which we get
\begin{equation}\label{u-eta}
\eta^\ast=-\dfrac{\ou_{1\ast}+2\ou_{2\ast}}{\sqrt5}.
\end{equation}
\begin{figure}[ht]
\centering
\includegraphics[scale=0.4]{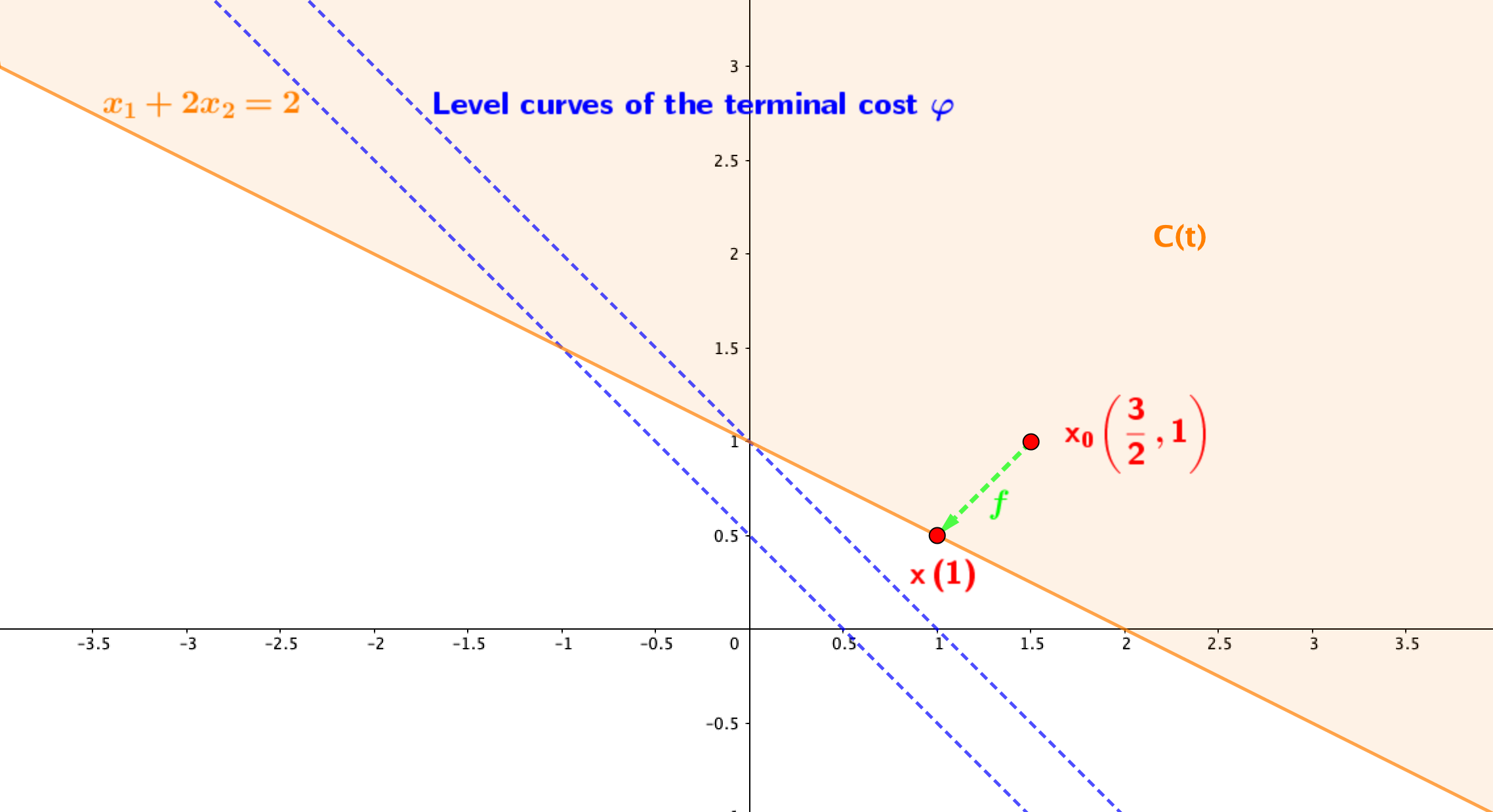}
\caption{The optimal trajectory for Example~\ref{ex_1}.}
\label{opt_ex_1}
\end{figure}
\end{example}\vspace*{-0.1in}
The cost functional is calculated therefore by
$$
J[\ox,\ou]= t^\ast\(\dfrac{\ou^2_{10}}{2}+\ou^2_{20}+\ou_{10}+\ou_{20}\)+(1-t^\ast)\(\dfrac{\ou^2_{1\ast}}{2}
+\ou^2_{2\ast}+\ou_{1\ast}+\ou_{2\ast}+\dfrac{3\eta^\ast}{\sqrt5}\)+\frac{5}{2},
$$
where the numbers $t^\ast$ and $\eta^\ast$ are taken from \eqref{htnew} and \eqref{u-eta}, respectively. The unique candidate for optimal control is now determined via direct calculations as follows:
$$
\ou_{10}=-\frac{5}{6}\;\mbox{ and }\;\ou_{20}=-\frac{1}{3}.
$$
Since in this case the hitting time is $t^\ast=1$, the controls $\bar{u}_{1\ast},\bar{u}_{2\ast}$ are irrelevant, and the optimal cost is $J[\ox,\ou]=\frac{43}{24}$. Observe finally that the case of $\lambda=0$ cannot occur. Indeed, supposing the contrary and taking into account that $\ox_1(t)+2\ox_2(t)>2$ for all $t<1$ tells us that
$\bar{u}_{10},\bar{u}_{20}\in(-1,1)$. It then follows from item~6 that $q^x(t)=0$. Invoking further items~8 and 5 yields $p^x(1)$=0, which clearly contradicts the enhanced nontriviality condition.

The two subsequent examples concern sweeping control problems with {\em controlled moving sets}. First we consider the case where controls are acting in the normal vectors generating the moving polyhedron.\vspace*{-0.1in}

\begin{example}[\bf controls in generating normals of moving polyhedra]\label{ex_2} Let $n=2$, $m=1$, and $T=1$. Consider control functions $a(t):=(\cos\vartheta(t),\sin\vartheta(t))$ acting in the moving set with $\dot{\vartheta}(t)\in[-\frac{\pi}{2},\frac{\pi}{2}]$ and $\vartheta(0)=\frac{\pi}{2}$. Then we form the controlled moving polyhedron by
$$
C(t):=\big\{x\in\mathbb{R}^2\;\big|\;\la a(t),x\ra\le 0\big\},
$$
which means in the notation of problem $(P)$ that $\dot{a}(t)\in A:=\{v\in\R^2\;|\;\|v\|\le\frac{\pi}{2}\big\}$ and $b(t)\equiv 0$. Set also $g(x,u)\equiv(0,3)$ and $x_0=(-1,-1)$ meaning that the sweeping dynamics reads as
\begin{equation}\label{ex_xdot}
\dot{x}(t)\in-N\big(x(t);C(t)\big)+(0,3)=-\eta(t)a(t)+(0,3)\;\mbox{ for a.e. }\;t\in[0,1]\;\mbox{ and }\;x(0)=(-1,-1)
\end{equation}
with a suitable $\eta(t)\ge 0$ a.e.\ on $[0,1]$ including $t=1$. Our goal is to minimize the cost functional
$$
J[x,a]:=\frac12\Big(\|x(1)\|^2 +\int_0^1\|\dot{a}(t)\|^2(t)\,dt\Big).
$$
Observing that $\dot{a}(t)=\dot{\vartheta}(t)(-\sin\vartheta(t),\cos\vartheta(t))$, the cost functional reduces to
$$
J[x,\vartheta]=\frac12\Big(\|x(1)\|^2+\int_0^1\dot{\vartheta}^2\,dt\Big).
$$
Since $x(0)\in\textrm{int}\,C(0)$, any trajectory of \eqref{ex_xdot} remains in the interior of $C(t)$ on the interval $[0,t^\ast)$, where
$$
t^\ast=\inf\big\{t>0\;\big|\;x(t)\in{\rm bd}\,C(t)\big\}.
$$
The main necessary conditions of Theorem~\ref{Th2} give us the following for the problem under consideration:
\begin{enumerate}
\item $p^x(t)\equiv p^x=\lambda \bar{x}(1)-\eta(1)\bar{a}(1)$.
\item $q^x(t)=p^x-\disp\int_{[t,1]}\bar{a}(s)\,d\gamma(s)$, and so $q^x\equiv p^x-\disp\int_{[t^\ast,1]}\bar{a}(s)\,d\gamma(s)$ is constant on $[0,t^\ast)$.
\item $\dot{p}^a(t)=\eta(t)q^x(t)=\eta (t)\big(p^x-\int_{[t,1]}\bar{a}(s)d\gamma(s)\big)$, which implies in turn that $\dot{p}^a\equiv 0$ a.e.\ on $[0,t^\ast)$ since $\eta\equiv 0$ a.e.\ on $[0,t^\ast)$, $p^a(1)+\eta(1)\bar{x}(1)\in a(1)\mathbb{R}$, $q^a(t)=p^a(t)-\disp\int_{[t,1]}\big(2\bar{a}(s)d\al(s)+\bar{x}(s)\,d\gamma(s)\big)$.
\item $\psi^a(t)=q^{a}(t)-\lm\dot{\bar a}(t)\in N(\dot{\bar a}(t);A)\;\textrm{ a.e. }\;t\in[0,1]$.
\item $\eta(1)\bar{a}(1)\in N(\bar{x}(1);\bar{C}(1))$.
\item If $\eta(t)>0$, then it follows from \eqref{41} that $\langle q^x(t),\bar{a}(t)\rangle=0$ a.e.
\item $(\lm,p(1),q(\cdot))\ne(0,0,0)$ with $\lm\ge 0$.
\end{enumerate}
Item~3 implies that $q^a$ is constant in $[0,t^\ast)$, i.e., $q^a(t)\equiv\bar{q}^a$ for all $t\in[0,t^\ast)$, since the measures $\alpha$ and $\gamma$ are not supported on intervals where the constraint is not active. If $\lambda>0$, then item~4 tells us that if $|\dot{\bar a}(t)|<\frac{\pi}{2}$ on a set $E$ of positive measure in $[0,t^\ast)$, then $\dot{\bar a}(t)$ has the constant value $q^a$, which yields $\dot{\bar a}(t)=0$ on that set. Observe that $\dot{\vartheta}<0$ on a set of positive measure means a counterclockwise rotation on that time set, which
evidently means wasting energy, as the cost of the final position is not made smaller. Therefore, the case where $\lambda >0$ ensures that either $\dot{\bar{\vartheta}}(t)=0$, or $\dot{\bar{\vartheta}}(t)=\frac{\pi}{2}$ a.e.\ on $[0,t^\ast)$. Set now $\dot{\vartheta}_0:=\frac{1}{t^\ast}\int_0^{t^\ast}\dot{\bar{\vartheta}}(s)\, ds$, and hence $\bar{\vartheta}(t^\ast)=\frac{\pi}{2}+t^\ast \dot{\vartheta}_0$. Moreover, by the strict convexity of the integrand yields $\int_0^{t^\ast} \dot{\bar{\vartheta}}^2(t)\,dt>t^\ast(\dot{\vartheta}_0)^2$. This tells us that on the interval $[0,t^\ast)$ there is no loss of generality in considering a constant angular velocity $\dot{\vartheta}_0\in[0,\frac{\pi}{2}]$.

Recalling further that $\dot{\oa}(t)=\dot{\vartheta}(t)(-\sin\vartheta(t),\cos\vartheta(t))$ and that $\dot{\bar{a}}$ has bounded variation (hence it is a.e.\ differentiable) on $[0,1]$ gives us the equality
$$
\ddot{\oa}(t)=\ddot{\vartheta}(t)\big(-\sin\vartheta(t),\cos\vartheta(t)\big)-\dot{\vartheta}^2(t) \big(\cos\vartheta(t),\sin\vartheta(t)\big),
$$
and thus we arrive at the following representation, which is valid on a time set of full measure:
\begin{eqnarray*}
\big\la\ddot{\oa}(t),\dot{\oa}(t)\big\ra &=&\big\la\big(-\ddot{\vartheta}(t)\sin\vartheta(t)-\dot{\vartheta}^2(t)\cos\vartheta(t),\ddot{\vartheta}(t)\cos\vartheta(t)-\dot{\vartheta}^2(t)
\sin\vartheta(t)\big),\;\dot{\vartheta}(t)\big(-\sin\vartheta(t),\cos\vartheta(t)\big)\big\ra\\
&=&\ddot{\vartheta}(t)\dot{\vartheta}(t)\sin^2\vartheta(t) + \dot{\vartheta}^3(t)\cos\vartheta(t)
\sin\vartheta(t)-\ddot{\vartheta}(t)\dot{\vartheta}(t)\cos^2\vartheta(t)-\dot{\vartheta}^3(t)\sin\vartheta(t)\cos\vartheta(t)=\dot{\vartheta}(t)\ddot{\vartheta}(t).
\end{eqnarray*}
Since the drift $g\equiv(0,3)$ is pushing on the boundary of $C(t)$, it follows that $\eta (t)>0$ for a.e.\ $[t^\ast,1]$. Then item~6 implies that
$\langle q^x(t),\bar{a}(t)\rangle =0$. By item~4 we have that $q^a(t)$ is a nonnegative multiple of $\dot{\bar{a}}(t)$, and hence it is expressed similarly to $q^x(t)$. This allows us to
make for simplicity the extra assumption that also $\langle q^a(t)\bar{a}(t)\rangle=0$. Under this assumption, the above representation implies that
$\dot{\vartheta}\ddot{\vartheta}=0$ a.e.. This means that the function $\dot{\vartheta}(t)$ is piecewise constant up to a possible Cantor part, which we neglect for simplicity.

It follows from $\vartheta(0)=\dfrac{\pi}{2}$ that $t^\ast$ solves the equation
$$
-1+3t=\mbox{cotan}\Big(\dot{\vartheta}_0t+\dfrac{\pi}{2}\Big)\Big(=-\tan\(\dot{\vartheta}_0t\)\Big).
$$
This tells us that, e.g., that if $\dot{\vartheta}_0=0$, then $t^\ast:=t^\ast_0=\frac{1}{3}$, while in the case where $\dot{\vartheta}_0=\dfrac{\pi}{2}$ we get that $t^\ast:=t^\ast_1\approx 0.215803$.
By arguing as above on the interval $[0,t^\ast)$, we can assume that $\dot{\vartheta}$ is constant also on $(t^\ast,1]$; say, $\dot{\vartheta}\equiv\dot{\vartheta}_\ast$. Observe that
\begin{equation}\label{(a)}
\la \ox(t),\oa(t)\ra=0\;\mbox{ whenever }\;t\in[t^\ast,1]
\end{equation}
meaning that on this interval the sweeping trajectory $\ox(t)$ slides on the rotating line perpendicular to $\oa(t)$. Differentiating \eqref{(a)} and recalling that $\ox(\cdot)$ is a solution to \eqref{ex_xdot}, i.e.,
\begin{equation}\label{x_sol}
\dot{\ox}(t)=-\eta(t)\big(\cos\vartheta(t),\sin\vartheta (t)\big)+(0,3)
\end{equation}
for a.e.\ $t\in[t^\ast,1]$, we have the equation
$$
-\eta(t)\cos^2\vartheta (t)+\big(3-\eta(t)\sin\vartheta(t)\big)\sin\vartheta(t)-\dot{\vartheta}_\ast x_1(t)\sin\vartheta(t)+\dot{\vartheta}_\ast x_2(t)\cos\vartheta(t)=0
$$
the solution of which gives us the formula
\begin{equation}\label{eta(t)}
\eta(t)=3\sin\vartheta(t)-\dot{\vartheta}_\ast\big(x_1(t)\sin\vartheta(t)-x_2(t)\cos\vartheta(t)\big).
\end{equation}
Observe first that $\eta(t)\ge 0$ as long as $\vartheta(t)\le\pi$. Moreover, \eqref{eta(t)} yields the following three facts:\\[1ex]
$\bullet$ on the interval $[t^\ast,1]$ we get that $\eta(\cdot)$ has the same regularity as $\vartheta(\cdot)$; in particular, it is smooth.\\[0.05in]
$\bullet$ $\eta^\ast:=\eta(t^{\ast})=3\sin\vartheta^\ast-\dot{\vartheta}_\ast(x_1^\ast\sin\vartheta^\ast-x_2\cos\vartheta^\ast)$, where $\vartheta^\ast=\vartheta (t^\ast)$ and $(x_1^\ast,x_2^\ast)=(x_1(t^\ast),x_2(t^\ast))$.\\[0.05in]
$\bullet$ If $\dot{\vartheta}_\ast=0$, then $\eta(t)$ is constant with $\eta(t)\equiv 3\sin\vartheta_\ast$.

Differentiating \eqref{eta(t)} tells us that for a.e.\ $t\in[t^\ast,1]$ we have the equality
$$
0=-\dot{\eta}(t)+3\dot{\vartheta}_\ast\cos\vartheta(t)-\dot{\vartheta}_\ast\big(\dot{x}_1(t)\sin\vartheta (t)+\dot{\vartheta}_\ast x_1(t)\cos\vartheta(t)-\dot{x}_2(t)\cos\vartheta(t)+\dot{\vartheta}_\ast x_2(t)\sin\vartheta(t)\big),
$$
which yields by using \eqref{(a)} and \eqref{x_sol} the differential equation
$$
\dot{\eta}(t)=6\dot{\vartheta}_\ast\cos\vartheta(t)\;\mbox{ for a.e.}\;t\in[t^\ast,1].
$$
The integration of the above equation gives us the exact formula
$$
\eta (t)=\eta^\ast+6\big(\sin(\dot{\vartheta}_\ast(t-t^\ast)+\vartheta^\ast)-\sin\vartheta^\ast\big)\;\mbox{ for all }\;t\in[t^\ast,1].
$$
Consequently, the sweeping trajectory $\ox(\cdot)$ can be computed explicitly on $[t^\ast,1]$, and thus the corresponding cost value $J$ can be explicitly computed as well, although we do not write down their (long) expressions for the sake of brevity. Note that $J$ can be written as a function of two one-dimensional variables; namely, of the switching time $t^\ast\in[t^\ast_1,t^\ast_0]$
and of the control $\dot{\vartheta}_\ast\in[-\frac{\pi}{2},\frac{\pi}{2}]$. By using computer calculations, all the cases can be analyzed. This shows that the best performance is obtained by choosing $\bar{t}^\ast\approx 0.270266$, which gives us $\dot{\bar{\vartheta}}_0\approx 0.691889$, $\dot{\bar{\vartheta}}_\ast\approx 0.350021$, and the minimal cost value $\bar{J}\approx 0.167854$.
The optimal trajectory (with a final time slightly larger than $1$ in order to show that it is meaningful and approaches the origin) is depicted in Figure~\ref{opt_ex_2}.
\begin{figure}[ht]
\centering
\includegraphics[scale=0.35]{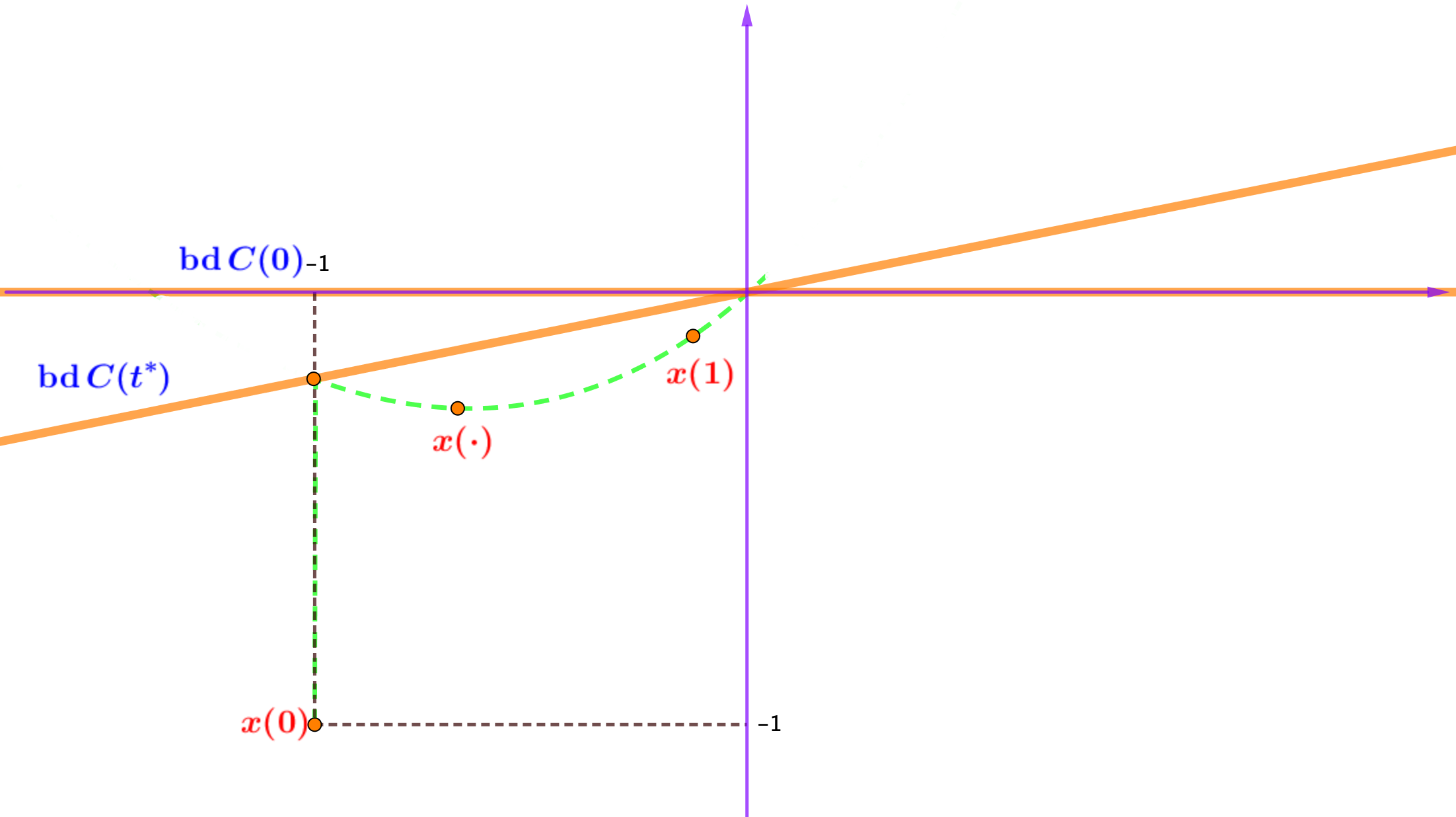}
\caption{The optimal trajectory for Example~\ref{ex_2}.}
\label{opt_ex_2}
\end{figure}
\end{example}\vspace*{-0.1in}

The last example addresses the setting where controls are use only to shift polyhedron edges.\vspace*{-0.1in}

\begin{example}[control by shifting polyhedron edges]\label{ex_3} Consider the sets
\begin{eqnarray*}
C_1(t):&=&\big\{x=(x_1,x_2)\in\R^2\;\big|\;\la x,a_1\ra\le 1\big\},\\
C_2(t):&=&\big\{x_1,x_2)\in\R^2\;\big|\;\la x,a_2\ra\le b(t)\big\},\\
C(t):&=&C_1(t)\cap C_2(t)\;\mbox{ for all }\;t\in[0,1],
\end{eqnarray*}
where $a_1:=\frac{(1,1)}{\sqrt{2}}$, $a_2:=(0,1)$,
$b(0):=\frac{3}{2}$, and $\dot{b}(t)\in[-1,1]$ for a.e.\ $t\in [0,1]$. Define the Cauchy problem for the controlled sweeping dynamics
\begin{equation}\label{ex_CP}
\dot{x}(t)\in-N\big(x(t);C(t)\big)+(0,2)\;\mbox{ a.e. }\;t\in[0,1],\;x(0)=(0,1),
\end{equation}
with controls $b(t)$ acting in the edges of the moving polyhedron $C(t)$. The sweeping dynamics \eqref{ex_CP} can be equivalently described as follows:
$$
\dot{x}(t)=-\eta_1(t)\frac{(1,1)}{\sqrt{2}}-\eta_2(t)(0,1)+(0,2),
$$
where $\eta_1(t),\eta_2(t)\ge 0$ and $\eta_i(t)>0\sr x(t)\in{\bd C}_i(t)$ as $i=1,2$ for a.e.\ $t\in[0,1]$.

We aim at minimizing the cost functional
$$
J[x,b]:=\frac12\Big(\|x(1)\|^2+\int_0^1\dot{b}^2(t)dt\Big)
$$
over all the feasible solutions to \eqref{ex_CP} as formulated in Theorem~\ref{Th2}. The main necessary optimality conditions obtained in that theorem ensure
the existence of $\lm\ge 0$, $\eta=\(\eta_1(t),\eta_2(t)\)\in\R^2_+$, measures $\gg_1,\,\gg_2$, absolutely continuous functions $p^x,p^b$, and BV functions $q^x,q^b$ such that
\begin{enumerate}
\item $p^x(t)\equiv p^x(1)=q^x(1)=-\lambda x(1)-\eta_1(1)\disp\frac{(1,1)}{\sqrt{2}}-\eta_2(1)(0,1)$.
\item $q^x(t)=p^x(1)-\disp\frac{(1,1)}{\sqrt{2}}\gamma_1([t,1])-(0,1)\gamma_2([t,1])$ with $\eta_i(t)>0\sr\la a_i,q^x(t)\ra=0$, $i=1,2$.
\item $\dot{p}^b(t)=0$ a.e.\ $t\in[0,T]$, $p^b\equiv\eta_2(1)\ge 0$.
\item $q^b(t)=p^b-\gamma_2([t,1])$ a.e.\ $t\in[0,T]$.
\item $\psi^b (t)=q^b(t)-\lambda\dot{\ob}(t)\in N(\dot{b}(t);[-1,1])$ a.e.\ $t\in[0,T]$.
\item $\big(\lambda,p^x(1),p^b(1),q^1(\cdot),q^2(\cdot)\big)\ne(0,0,0,0,0)$.
\end{enumerate}
Since the initial point $x(0)$ belongs to the interior of $C_2(0)$, we have $t^\ast:=\inf\{t>0\;|\;x(t)\in{\bd C}_2(t)\}>0$. The structure of the dynamics implies that $x(t)\in{\bd C}_1(t)$ for all $t\in[0,t^\ast]$ while $x(t)$ belongs to the interior of $C_1(t)$ on $(t^\ast,1]$. It follows from item~2 that
$$
q^x(t)=p^x(1)-\gg_1([t,t^\ast])\frac{(1,1)}{\sqrt{2}}-\gg_2([t^\ast,1])(0,1)\;\textrm{ a.e. in }\;[0,t^\ast],
$$
$$
q^x(t)=p^x(1)-\gg_2([t,1])(0,1)\;\textrm{ a.e. in }\;(t^\ast,1].
$$
Moreover, we deduce from item~4 that $q^b(t)\equiv q^b=p^b-\gg_2([t^\ast,1])$ is a constant nonnegative multiple of $\dot{\ob}(t)$ for a.e.\ $t\in[0,t^\ast]$. If
therefore $|\dot{\ob}(t)|<1$ on this interval, then $\psi^b(t)=0$ and thus $\lambda \dot{\ob}(t)=q^b$ for a.e.\ $t\in[0,t^\ast)$. Now it is easy to observe from the problem structure that $\dot{b}(t)\le 0$ for a.e.\ $t\in[0,t^\ast]$. Indeed, $\dot{\ob}(t)>0$ means spending energy without pushing the trajectory towards the origin, i.e., without affecting the cost of the final position. This tells us that if $\lm>0$, then either $-1<\dot{b}(t)\le 0$ (and in this case $\dot{\ob}\equiv b_0$, a constant value), or $\dot{b}(t)=-1$. Hence $\dot{\ob}(t)$ can get only two values. Assuming for simplicity that $\dot{b}$ is piecewise continuous implies that it is piecewise constant. This leads us to treating controls $\dot{b}$ as piecewise constant functions with finitely many switchings. The latter is equivalent to considering a constant control $\dot{b}(t)\equiv\dot{b}_0$ on the interval $[0,t^\ast]$. Indeed, arguing as in Example~\ref{ex_2} confirms that the same position on the boundary of $C_2(t)$ can be obtained by using a constant control, and that the strict convexity of the integrand implies that a constant control ensures a better performance. Observe also that the problem dynamics \eqref{ex_CP} implies that $\eta_1(t)>0$ and $\eta_2\equiv 0$ for a.e.\ $t\in [0,t^\ast]$. Thus the usage of the second part of item~2 ensures that
$$
0=\la a_1,q^x(t)\ra=\la a_1,p^x(1)-a_1\gg_1([t,1])-a_2\gg_2\([t^\ast,1]\)\ra\;\mbox{a.e. on }\;[0,t^\ast),
$$
which yields in turn the relationship
$$
p^{x_1}-\sqrt{2}\gg_1([t,t^\ast])+p^{x_2}-\gamma_2([t^\ast,1])=0\;\mbox{ a.e. on }\;[0,t^\ast),
$$
where $p^x(1)=(p^{x_1},p^{x_2})$. This tells us that $q^x$ is constant on $[0,t^\ast)$. Furthermore, for a.e.\ $t\in [0,t^\ast]$ we get
$$
(\dot{x}_1,\dot{x}_2)=(0,2)-\Big\langle(0,2),\frac{(1,1)}{\sqrt{2}}\Big\rangle\frac{(1,1)}{\sqrt{2}}=(-1,1),
$$
and hence $x_1(t)=-t,\,x_2(t)=1+t$. Taking $b(t)=\frac32 +\dot{b}_0t$ gives us
\begin{equation}\label{tstar}
t^\ast=t^\ast(\dot{b}_0)=\frac12\frac{1}{1-\dot{b}_0}.
\end{equation}
On the interval $[t^\ast,1]$ we have $\la(1,1),x(t)\ra<1$ and $\eta_2(t)>0$ while $\eta_1\equiv 0$. Then using item~2 gives us
$$
0=\la(0,1),q^x(t)\ra=\sqrt{2}p^{x_2}-\gamma_1({t^\ast})-\sqrt{2}\gg_2([t,1]),
$$
which implies, in particular, that $\gg_2([t,1])$ is constant on $[t^\ast,1]$. Assuming for simplicity that the function $\dot{b}$ is piecewise continuous on $[t^\ast,1]$, we obtain from item~4 that $\dot{b}$ is piecewise constant on $[t^\ast,1]$. Then arguing as in the above case of $[0,t^\ast)$, we get that $\dot{b}$ is constant on $[t^\ast,1]$, i.e., $\dot{b}(t)\equiv\dot{b}^\ast$. This tells us that
$$
x_1(t)\equiv x_1(t^\ast)\;\mbox{ and }\;x_2(t)\equiv x_2(t^\ast)+\dot{b}^\ast(t-t^\ast)\;\mbox{ on }\;[t^\ast,1].
$$
Thus the corresponding value of the cost functional is
$$
J(\dot{b}_0,\dot{b}_\ast)=\dfrac12\[(t^\ast)^2+\big(1+t^\ast+\dot{b}^\ast(1-t^\ast)\big)^2+t^\ast\dot{b}_0^2+(1-t^\ast)\dot{b}_\ast^2\].
$$
Substituting into the above expression the value of $t^*$ from \eqref{tstar} gives us
$$
J(\dot{b}_0,\dot{b}_\ast)=\frac{10-2\dot{b}_0+6\dot{b}_\ast+3\dot{b}_\ast^2-2\dot{b}_0(6
+8\dot{b}_\ast+5\dot{b}_\ast^2)+\dot{b}_0^2(6+8\dot{b}_\ast+8\dot{b}_\ast^2)}{8(1-\dot{b}_0)^2}.
$$
Employing  finally computer calculations, we arrive at the optimal control $(b_0\approx-0.877931,\,b_\ast\approx-0.730354)$ with the optimal cost value $J\approx 0.600458$. Figure~\ref{opt_ex_3} illustrates behavior of the optimal solution.
\begin{figure}[ht]
\centering
\includegraphics[scale=0.4]{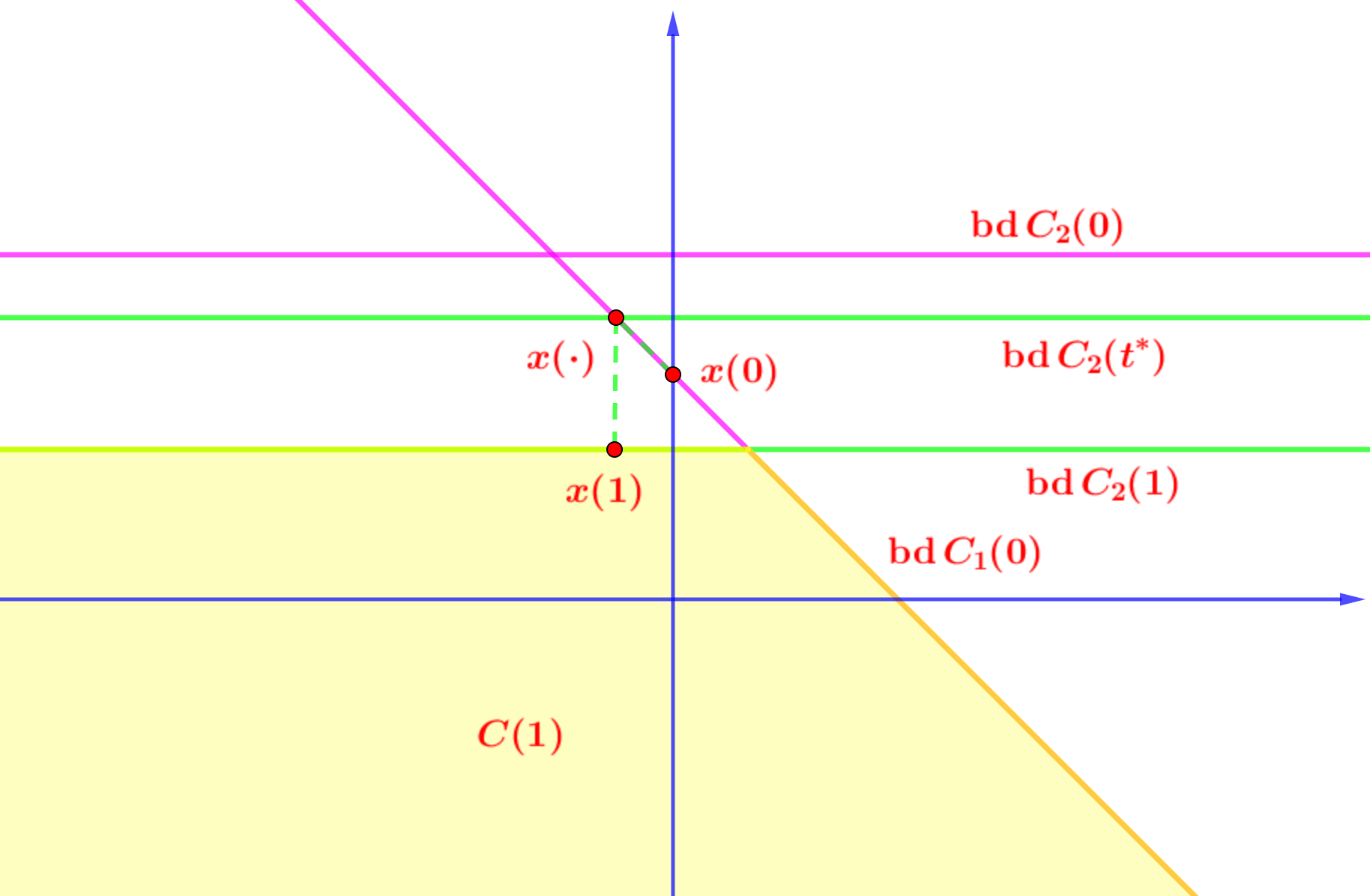}
\caption{The optimal trajectory for Example \ref{ex_3}.}
\label{opt_ex_3}
\end{figure}
\end{example}\vspace*{-0.1in}
{\bf Acknowledgements}. The authors are very grateful to both anonymous referees for their helpful remarks that allowed us to improve the original presentation.
\end{document}